\documentclass[final]{elsarticle}

\makeatletter
\def\ps@pprintTitle{%
    \let\@oddhead\@empty
    \let\@evenhead\@empty
    \let\@oddfoot\@empty
    \let\@evenfoot\@empty
    }
\makeatother

\usepackage[utf8]{inputenc}
\usepackage{geometry}
\usepackage{array}
\usepackage{amsmath,bm,times}
\usepackage{graphicx}
\usepackage{natbib}
\usepackage{caption}
\usepackage{hyperref}
\usepackage{cleveref}
\usepackage{subcaption}
\usepackage{bbm}
\usepackage{bbold}
\usepackage{enumitem}
\usepackage{multicol}
\usepackage{multirow}
\usepackage{comment}
\usepackage{mathtools}
\usepackage{textcomp}
\usepackage{amsthm}
\usepackage{soul}
\usepackage{xcolor}
\usepackage{algorithm}
\usepackage{algpseudocode}
\usepackage{siunitx}
\usepackage{tabularx}
\usepackage{booktabs}
\usepackage{tikz} 
\usetikzlibrary{shapes.geometric, positioning, arrows.meta}

\newcolumntype{L}[1]{>{\raggedleft\let\newline\\\arraybackslash\hspace{0pt}}p{#1}}
\newcommand{\rowgroup}[1]{\hspace{-1em}#1}

\setcitestyle{authoryear,round}

\begin{document}

\begin{frontmatter}
\journal{Water Research}
\title{Optimal design-for-control of self-cleaning water distribution networks using a convex multi-start algorithm}

\author[mymainaddress]{Bradley Jenks}

\author[mymainaddress]{Filippo Pecci}
\author[mymainaddress]{Ivan Stoianov}

\address[mymainaddress]{Department of Civil and Environmental Engineering, Imperial College London, London SW7 2BB, United Kingdom}

\begin{abstract}
The provision of self-cleaning velocities has been shown to reduce the risk of discolouration in water distribution networks (WDNs). Despite these findings, control implementations continue to be focused primarily on pressure and leakage management. This paper considers the control of diurnal flow velocities to maximize the self-cleaning capacity (SCC) of WDNs. We formulate a new optimal design-for-control problem where locations and operational settings of pressure control and automatic flushing valves are jointly optimized. The problem formulation includes a nonconvex objective function, nonconvex hydraulic conservation law constraints, and binary variables for modelling valve placement, resulting in a nonconvex mixed integer nonlinear programming (MINLP) optimization problem. Considering the challenges with solving nonconvex MINLP problems, we propose a heuristic algorithm which combines convex relaxations (with domain reduction), a randomization technique, and a multi-start strategy to compute feasible solutions. We evaluate the proposed algorithm on case study networks with varying size and degrees of complexity, including a large-scale operational network in the UK. The convex multi-start algorithm is shown to be a more robust solution method compared to an off-the-shelf genetic algorithm, finding good-quality feasible solutions to all design-for-control numerical experiments. Moreover, we demonstrate the implemented multi-start strategy to be a fast and scalable method for computing feasible solutions to the nonlinear SCC control problem. The proposed method extends the control capabilities and benefits of dynamically adaptive networks to improve water quality in WDNs.

\vspace{0.25cm}
\end{abstract}

\end{frontmatter}


\section{Introduction}
\label{sec:intro}
The management of water quality in water distribution networks (WDNs) presents a complex operational challenge. As a direct consequence of ageing and deteriorating infrastructure, the mitigation of discolouration incidents is becoming one of the key operational challenges within water quality management programmes. In addition to discolouration being the largest source of customer complaints \citep{VREEBURG2007,HUSBAND2011,ARMAND2017}, there is a growing body of evidence suggesting its occurrence harbours increased microbial activity \citep{LIU2013,LIU2014,VANDERWIELEN2016}. These conditions can accelerate biofilm growth and drastically reduce the efficacy of disinfectant residuals in protecting against waterborne illness and contaminants intrusion. Moreover, WDNs in the UK are highly sectorized and operated with fixed topology for purposes of leakage management. This type of network configuration, referred to as district metered areas (DMAs), has been demonstrated to exacerbate water quality deterioration and increase the risk of discolouration incidents \citep{MACHELL2014,ARMAND2018}. With progressively stringent water quality regulations, water companies are seeking effective and cost-efficient operational control strategies to reduce the risk of discolouration.  

Discolouration is primarily a consequence of resuspended material accumulated within WDNs \citep{VREEBURG2007}. It can materialize from the cumulative impact of the following processes \citep{BOXALL2005}: (i) the ingress and/or development of particulate matter; (ii) the accumulation of particulates at the pipe invert and/or formation of cohesive layers at the pipe wall; and (iii) a hydraulic disturbance (i.e. trigger event), which mobilizes loose particulates and generates sufficient shear stress to overcome cohesive forces at the pipe wall. Such hydraulic disturbances can be generated from different phenomena, including pressure transients during unsteady hydraulic conditions \citep{AISOPOU2012}. Apart from their origin, the physical pathways of discolouration are intrinsically connected to network hydraulics. In a recent study focusing on the impact of network sectorization on water quality, \citet{ARMAND2018} proposed a set of surrogate hydraulic variables for discolouration risk assessment. Central to their findings was the role of diurnal flow velocities on particle transport and fate. This connection between discolouration and hydrodynamic conditions has been supported by numerous experimental and theoretical studies; see \citet{VANSUMMEREN2017} and \citet{ARMAND2018} for reviews on the topic. These studies have mainly focused on the development of predictive tools for modelling particle transport and accumulation processes. Most notably, \citet{BOXALL2001} developed the Prediction of Discolouration in Distribution Systems (PODDS) model, an empirically-based numerical tool which aims to characterize cohesive layer strength at the pipe wall. The PODDS model was later updated to account for material regeneration in \citet{FURNASS2014}, where both erosion and regeneration processes require calibration using continuous flow and turbidity data. Because such tools require extensive field testing and are generally limited to pipe-level assessments, their use in practice has not yet been widespread. Recognizing this limitation, \citet{VANSUMMEREN2017} presented a theoretical particle transport model, combining the effects of gravitational settling, hydraulic shear stresses, and bed-load transport. To complement this, several laboratory-based experimental studies have emerged to better understand the complex interactions between particle properties and pipe hydraulics (e.g. \citeauthor{SHARPE2019}, \citeyear{SHARPE2019}; \citeauthor{BRAGA2020}, \citeyear{BRAGA2020}).

In addition to predictive modelling, research has also focused on reducing the severity and frequency of discolouration incidents through network design, maintenance, and control. Water companies in the Netherlands have been conducting experimental research on the design and implementation of controls for self-cleaning networks. The self-cleaning capacity (SCC) of a WDN is defined as the ability for pipes to experience peak daily flow velocities above a threshold required to routinely re-suspend particles and thus prevent accumulation \citep{VREEBURG2009}. Previous experimental programmes have suggested resuspension velocities on the order of 0.2 \si{\meter/\second} to 0.25 \si{\meter/\second} in distribution pipes \citep{RYAN2008,BLOKKER2010}. This has been corroborated with a recent field study monitoring turbidity under various flow rates, where an increase in turbidity levels were observed at flow velocities greater than 0.2 \si{\meter/\second} \citep{PREST2021}. Water companies in the Netherlands have demonstrated successful self-cleaning implementations by redesigning looped, oversized networks to branched layouts with smaller diameter pipes \citep{VREEBURG2009}. A recent study has also investigated the trade-off between self-cleaning velocities and fire flow capacity in North American WDNs \citep{GIBSON2019}. However, since the redesign of WDN infrastructure becomes cost-prohibitive at scale, there have been recent forays in the reconfiguration of existing network topology to promote self-cleaning networks \citep{BLOKKER2012,ABRAHAM2016,ABRAHAM2018}.   

Combining UK and Dutch experience, \citet{ABRAHAM2016,ABRAHAM2018} formulated an optimization problem for increasing SCC by redistributing flow through changes in network topology. More specifically, the optimization problem aimed to maximize the number of pipes with flow velocities above a self-cleaning threshold through two separate strategies: (i) optimal closure of isolation valves and (ii) optimal operational settings of existing pressure control valves \citep{ABRAHAM2016,ABRAHAM2018}. \citet{ABRAHAM2018} solved the problem of optimizing valve closures using a linear graph analysis tool. Following \citet{SCHAUB2014}, a line-outage distribution factor (LODF) matrix was computed to estimate the flow redistribution resulting from an outage (closure) of an (or multiple) edge-to-edge relation(s). For the optimal control problem, \citet{ABRAHAM2016} computed a local solution by approximating the nonsmooth objective function as a continuous nonlinear function, followed by application of a tailored sequential convex programming algorithm. While the benefits of the LODF solution method for optimal valve closures were demonstrated numerically using an operational network in the Netherlands, results from the control problem were limited to a small-scale theoretical network. Moreover, decision variables were restricted to the control of existing unidirectional pressure reducing valves (PRVs). Building on the SCC optimization problem posed in \citet{ABRAHAM2016,ABRAHAM2018}, this manuscript considers both control and design-for-control problem formulations. The latter involves the simultaneous optimization of valve placement and operational settings for both existing and new control valves. In addition to unidirectional PRVs, this work also considers bidirectional dynamic boundary valves (DBVs) and automatic flushing valves (AFVs) as dynamic hydraulic controls. These hydraulic controls were developed to facilitate the novel operational framework of dynamically adaptive networks \citep{WRIGHT2014,ULUSOY2022b}. The resulting optimization problem is formulated as a nonconvex mixed integer nonlinear program (MINLP).

Both mathematical optimization and heuristic methods have been used to solve design and control problems in WDNs (see literature review in \citeauthor{MALA2017}, \citeyear{MALA2017}). For mathematical optimization methods, scalability is recognized as a current limitation in solving MINLP problems to global optimality \citep{KOCH2012,SAHINIDIS2019}; that is, the implementation of global solvers become impractical for large problem cases. Consequently, heuristic approaches are often employed to compute satisfactory feasible solutions. A common heuristic method used for WDN optimization problems is the genetic algorithm (GA). While GAs have been successfully applied to design problems, the computational effort required to find solutions sufficiently close to the global optimum grows rapidly with problem size \citep{MAIER2014}. In this manuscript, we develop a heuristic algorithm based on convex optimization and a multi-start scheme to compute feasible solutions to the considered MINLP problem. To handle integer variables, we first formulate a convex subproblem through polyhedral relaxations of nonconvex terms and the continuous relaxation of binary variables. We subsequently employ a randomization heuristic to sample $N$ candidate valve configurations from the set of fractional values generated from the convex subproblem. We then fix binary variables for each sampled valve configuration and compute (locally) optimal operational settings from a nonlinear programming (NLP) control problem. This follows the heuristic algorithm presented in \citet{PECCI2022}, extending its application to the SCC design-for-control problem and to the nonsmooth Hazen-Williams friction model. Since the degree of nonlinearity of the SCC problem is higher than the problem investigated in \citet{PECCI2022}, we include a multi-start strategy and a feasibility restoration problem for selecting starting points. This step aims to minimize the risk of poor local optima as well as ensure hydraulic feasibility of the NLP control problem. Finally, the best feasible solution is selected from the set of sampled valve configurations. The proposed heuristic algorithm further increases the benefits from the implementation of dynamically adaptive networks as it expands their control capabilities to enhance water quality in WDNs. 

This manuscript is organized as follows. In \Cref{sec:prob_form}, we formulate the design-for-control problem of maximizing the network SCC through dynamic hydraulic controls. We then present the proposed heuristic algorithm in \Cref{sec:sol_methods}. Finally, in \Cref{sec:results}, we demonstrate the performance of the developed heuristic algorithm using three case study networks with varying size and degrees of complexity. To facilitate a broader discussion on heuristic approaches for the design and control of WDNs, we compare the results with an off-the-shelf GA implementation, which is a common approach used in the literature.

\section{Problem formulation}
\label{sec:prob_form}
We investigate a design-for-control problem to maximize the length of network pipes experiencing flow velocities above a given self-cleaning capacity (SCC) threshold. This is achieved by installing new valves and/or controlling their operational settings. For this purpose, our problem formulation considers three valve types as pressure and connectivity control actuators. First, pressure reducing valves (PRVs), which are modelled having unidirectional flow. Second, bidirectional dynamic boundary valves (DBVs), for which flow is permitted in both directions across discrete model time steps. Here, DBVs represent the operation of remote-controlled isolation valves, which modulate flow and pressure between adjacent zones. Third, automatic flushing valves (AFVs), whose flushing rate is bounded by a set maximum value. Throughout this manuscript, we refer to either PRVs or DBVs as \textit{control valves}, as both have the capability of controlling pressure and are modelled at network links. On the other hand, AFVs are simply referred to as \textit{flushing valves} and are modelled at network nodes. We consider operational scenarios, for which PRV locations have been fixed to minimize average zone pressure (AZP), and thus decision variables include only their operational settings. In comparison, both locations and operational settings of DBVs and AFVs are considered as decision variables. The operational settings of valves are modelled as continuous variables, whereas their placement (location) are modelled through binary variables. All network links and nodes are considered as potential locations of DBVs and AFVs, respectively. As the current stage of this work focuses on the self-cleaning capacity of pipes at the DMA or distribution level, we do not consider storage tanks or pumping activity as forms of hydraulic control. Therefore, we assume discrete and hydraulically independent model time steps. Finally, it is noted that this work relies on the availability of a calibrated hydraulic model.

\subsection{Hydraulic variables and constraints}
The problem considers a water distribution network (WDN) with $n_p$ links, $n_n$ demand nodes and $n_0$ known head nodes (e.g. water sources, reservoirs). The network is modelled as a directed graph with $n_p$ edges (links) and $n_n$ + $n_0$ vertices (nodes). A demand-driven hydraulic analysis is used to simulate steady-state network hydraulics over $n_t$ discrete time steps. For each time step $t \in \{1,\ldots,n_t\}$, known hydraulic conditions are given by vectors of nodal demands $d_t \in \mathbb{R}^{n_{n}}$ and source hydraulic heads $h^0_t \in \mathbb{R}^{n_{0}}$. Moreover, vectors $\eta_t \in \mathbb{R}^{n_p}$ and $\alpha_t \in \mathbb{R}^{n_n}$ are included to model local losses introduced by the action of control valves and operational demands at flushing valves, respectively. Unique vectors of hydraulic states $q_t \in \mathbb{R}^{n_p}$ and $h_t \in \mathbb{R}^{n_n}$ are computed by solving the following steady-state energy \eqref{eq:hydraulic_conservation_a} and mass \eqref{eq:hydraulic_conservation_b} conservation equations governing pipe flow:
\begin{subequations}
\label{eq:hydraulic_conservation}
\begin{align}
\begin{split}
\label{eq:hydraulic_conservation_a}
A_{12}h_t + A_{10}h^0_t + \phi(q_t) + \eta_t = 0
\end{split} \\
\begin{split}
\label{eq:hydraulic_conservation_b}
A_{12}^T q_t - d_t -\alpha_t = 0,
\end{split}
\end{align}
\end{subequations}
where $A_{12} \in \mathbbm{R}^{n_p \times n_n}$ and $A_{10} \in \mathbbm{R}^{n_p \times n_0}$ are the link-node incidence matrices for demand and known head nodes, respectively; and the vector $\phi(q_t) = [\phi_1(q_{1,t}) \dots \phi_{n_p}(q_{n_{p},t})]^T$ models frictional head losses associated with flows $q_t$. Omitting time index $t$, $\phi_j(q_j)$ is defined in general form for flow conveyed across link $j$ as: 
\begin{equation}
    \label{eq:head_loss_model}
    \phi_j(q_j) = r_j|q_j|^{n_j-1}q_j, \quad \forall j \in \{1,\ldots,n_p\},
\end{equation}
where the resistance coefficient $r_j$ and exponent $n_j$, both independent of time $t$, take different values depending on the link type (e.g. pipe or valve) and on the frictional head loss model. For valve links, $n_j = 2$ and $r_j = \frac{8K_j}{g\pi^2D_j^4}$, with $K_j$ and $D_j$ representing the valve loss coefficient and diameter, respectively \citep{LAROCK1999}. In this work, we apply the Hazen-Williams (HW) model to characterize frictional head losses across pipe links. The HW model is an explicit and empirical relationship between pipe flow and frictional head loss, with $n_j = 1.852$ and $r_j$ is defined as follows for all $j \in \{1,\ldots,n_p\}$: 
\begin{equation}
    \label{eq:HW_model}
    r_{j} = \frac{10.67L_j}{C_j^{1.852}D_j^{4.871}},
\end{equation}
where $C$ is the HW coefficient, a dimensionless number representing frictional characteristics; $L$ is pipe length in meters; and $D$ is pipe diameter in meters \citep{LAROCK1999}. Similarly, explicit approximations of the Darcy-Weisbach formula (e.g. \citeauthor{VALIANTZAS2008}, \citeyear{VALIANTZAS2008}) could be used to model frictional head losses. Finally, it is convenient to isolate the nonlinear term $\phi(q_t)$ in \eqref{eq:hydraulic_conservation_a}. Here, we introduce a vector of auxiliary variables $\theta_t \in \mathbb{R}^{n_{p}}$, which separates the energy conservation constraint into its linear and nonlinear components, as follows:
\begin{subequations}
\label{eq:energy_conservation}
\begin{align}
\begin{split}
\label{eq:hydraulic_conservation_c}
A_{12}h_t + A_{10}h^0_t + \theta_t + \eta_t = 0
\end{split} \\
\begin{split}
\label{eq:hydraulic_conservation_d}
\theta_t - \phi(q_t) = 0.
\end{split}
\end{align}
\end{subequations}

Valve placement and operation are modelled as follows. For each time step $t \in \{1,\ldots,n_t\}$, the continuous variable $\eta_t \in \mathbb{R}^{n_p}$ presented in \eqref{eq:hydraulic_conservation_a} models the local losses introduced by the action of control valves and the continuous variable $\alpha_t \in \mathbb{R}^{n_n}$ presented in \eqref{eq:hydraulic_conservation_b} models the flow emitted at flushing valves. Moreover, binary variables $z \in \{0,1\}^{n_p}$ are included to model PRV and DBV placement, and $v^+_t \in \{0,1\}^{n_p}$ and $v^-_t \in \{0,1\}^{n_p}$ to assign their control capabilities in the positive or negative flow direction, respectively, across each time step $t$. Thus, for all links $j \in \{1,\ldots,n_p\}$ and time steps $t \in \{1,\ldots,n_t\}$, binary variables $z_j$, $v^+_{j,t}$, and $v^-_{j,t}$ are set as
\begin{equation} 
\label{eq:valves_direction}
\begin{alignedat}{3}
    &z_j = \;\, \begin{cases} 1 &\text{control valve on link $j$} \\[-3pt]
    0 & \text{no valve} \end{cases} \\[5pt]
    &v^+_{j,t} = \begin{cases} 1 &\text{control valve on link $j$ in positive direction} \\[-3pt]
    0 & \text{no valve} \end{cases} \\[5pt]
    &v^-_{j,t} = \begin{cases} 1 &\text{control valve on link $j$ in negative direction} \\[-3pt]
    0 & \text{no valve} \end{cases}
\end{alignedat}
\end{equation}
Analogously, the placement of AFVs at network nodes is modelled using binary variables $y \in \{0,1\}^{n_n}$, defined as
\begin{equation} 
\label{eq:flush_binary}
\begin{alignedat}{3}
    &y_i = \begin{cases} 1 &\text{flushing valve placed at node $i$} \\[-3pt]
    0 & \text{no valve} \end{cases}
\end{alignedat}
\end{equation}
These binary variables are subject to the following physical and economical constraints, which limit pressure control capabilities in a single direction at each time step $t$ \eqref{eq:valves_physical_a} and enforce a maximum number of control valves $n_v$ and flushing valves $n_f$ considered for installation \eqref{eq:valves_physical_b}-\eqref{eq:valves_physical_c}:
\begin{subequations}
\label{eq:valves_physical}
\begin{align}
    &v^+_{j,t} + v^-_{j,t} \leq z_j, \quad \forall j \in \{1,\ldots,n_p\},\; \forall t \in \{1,\ldots,n_t\} 
    \label{eq:valves_physical_a}\\[3pt]
    &\sum_{j=1}^{n_p}z_j = n_v
    \label{eq:valves_physical_b}\\[3pt]
    &\sum_{i=1}^{n_n}y_i = n_f.
    \label{eq:valves_physical_c}
\end{align}
\end{subequations}
Since we assume that existing PRVs have fixed location and unidirectional flow, $z_j$, $v^+_{j,t}$, and $v^-_{j,t}$ are set \textit{a priori} for all time steps $t \in \{1,\ldots,n_t\}$ at the known set of PRV links $N_\text{PRV} \subseteq \{1,\ldots,n_p\}$.

We introduce constant vectors to bound the continuous hydraulic variables and formulate \mbox{big-M} constraints for modelling the operation of control and flushing valves. For a given vector of maximum allowed velocities $u^{\max} \in \mathbb{R}^{n_p}$ and vector of link cross-sectional areas $A \in \mathbb{R}^{n_p}$, let $q^L_t = -Au^{\max}$ and $q^U_t = Au^{\max}$ be the vectors of lower and upper bound flows across network links at time step $t$, respectively. Bounds on the auxiliary head loss vector $\theta_t$ are set as $\theta^L_t := \phi(q^L_t)$ and $\theta^U_t := \phi(q^U_t)$. Moreover, let $h^{\min}_t$ and $h^{\max}_t \in \mathbb{R}^{n_n}$ specify minimum and maximum heads at network nodes, respectively. The minimum head is set to a minimum regulatory pressure plus the node elevation and the maximum head is set to the largest available known source head. Bounds on $\eta_t$ for $j \in \{1,\dots,n_p\}$ are then defined as follows:
\begin{subequations}
\label{eq:eta_bounds}
\begin{align}
    &(\eta^L_t)_j := (h^{\min}_t)_i - (h^{\max}_t)_k, \quad \forall i \xrightarrow{j} k\\[3pt]
    &(\eta^U_t)_j := (h^{\max}_t)_i - (h^{\min}_t)_k, \quad \forall i \xrightarrow{j} k
\end{align}
\end{subequations}
We formulate \mbox{big-M} constraints to model valve placement and enforce energy conservation at control valve links, ensuring $\eta_t$ and $q_t$ act in the same direction. These constraints are written as follows:
\begin{subequations}
\label{eq:valve_bigM}
\begin{align}
&\eta_t - \text{diag}(\eta^U_t)v^+_t \leq 0, \quad \forall t \in \{1,\ldots,n_t\} \label{eq:valve_bigM_a}\\
&-\eta_t + \text{diag}(\eta^L_t)v^-_t \leq 0, \quad \forall t \in \{1,\ldots,n_t\} \label{eq:valve_bigM_b}\\
&-q_t - \text{diag}(q^L_t)v^+_t \leq -q^L_t, \quad \forall t \in \{1,\ldots,n_t\} \label{eq:valve_bigM_c}\\
&q_t + \text{diag}(q^U_t)v^-_t \leq q^U_t, \quad \forall t \in \{1,\ldots,n_t\} \label{eq:valve_bigM_d}\\
&-\theta_t - \text{diag}(\theta^L_t)v^+_t \leq -\theta^L_t, \quad \forall t \in \{1,\ldots,n_t\} \label{eq:valve_bigM_e}\\
&\theta_t + \text{diag}(\theta^U_t)v^-_t \leq \theta^U_t, \quad \forall t \in \{1,\ldots,n_t\}. \label{eq:valve_bigM_f}
\end{align}
\end{subequations}
Additionally, let $\alpha^U_t$ be a known upper bound on the flushing rate at AFVs. The following \mbox{big-M} constraint is then included to model flushing valve placement and enforce bounds on the continuous variable $\alpha_t$:
\begin{equation}
    \alpha_t - \text{diag}(\alpha^U_t)y \leq 0, \quad \forall t \in \{1,\ldots,n_t\}. \label{eq:valve_bigM_g}
\end{equation}
Finally, lower and upper bounds on hydraulic variables $h_t$, $q_t$, $\eta_t$ and $\theta_t$ are set to define the feasible solution space,
\begin{subequations}
\label{eq:hyd_bounds}
\begin{align}
    &q^L_t \leq q_t \leq q^U_t, \quad \forall t \in \{1,\ldots,n_t\}  \label{eq:hyd_bounds_a}\\
    &h^{\min}_t \leq h_t \leq h^{\max}_t, \quad \forall t \in \{1,\ldots,n_t\}  \label{eq:hyd_bounds_b}\\
    &\eta^L_t \leq \eta_t \leq \eta^U_t, \quad \forall t \in \{1,\ldots,n_t\}  \label{eq:hyd_bounds_c}\\
    &0 \leq \alpha_t \leq \alpha^U_t, \quad \forall t \in \{1,\ldots,n_t\} \label{eq:hyd_bounds_d}\\ 
    &\theta^L_t \leq \theta_t \leq \theta^U_t, \quad \forall t \in \{1,\ldots,n_t\}. \label{eq:hyd_bounds_e}
\end{align}
\end{subequations}

\subsection{Self-cleaning capacity objective function}
The objective of this study is to maximize the length of network pipes satisfying the self-cleaning capacity (SCC) threshold. The SCC objective function is defined as the following length-weighted sum over all pipes $n_p$ and hydraulic time steps $n_t$ \citep{ABRAHAM2018}:
\begin{equation}
\label{eq:SCC_objective}
\begin{aligned}
& & f_\text{SCC} := \frac{1}{n_t} \sum_{t=1}^{n_t} \sum_{j=1}^{n_p} w_j \kappa_j\bigg(\frac{q_{j,t}}{A_j}\bigg), \\
\end{aligned}
\end{equation}
where $A$ is the link cross-sectional area; and $\kappa_j\left(\cdot\right)$ is an indicator function which models the state of pipe velocities with reference to a minimum threshold. The indicator function is described for link $j$ as follows
\begin{equation}
\label{eq:SCC_indicator_function}
    \kappa_j(u) = \begin{cases} 1 &\text{if $|u|$ $>$ $u^{\min}_j$} \\
    0 & \text{otherwise},
    \end{cases}
\end{equation}
with $u^{\min}_j$ representing the threshold flow velocity at link $j$,  defined \textit{a priori} (see presented literature in \Cref{sec:intro}). Moreover, a weighting is included to normalize the length of link $j$ to the entire network, $w_{j} = \frac{L_j}{\sum_{k=1}^{n_p}L_k}$, where $L \in \mathbbm{R}^{n_p}$ is the vector of pipe lengths.

\begin{figure}[!h] 
    \begin{center}
    \captionsetup{justification=centering}
    \includegraphics[width=0.5\linewidth]{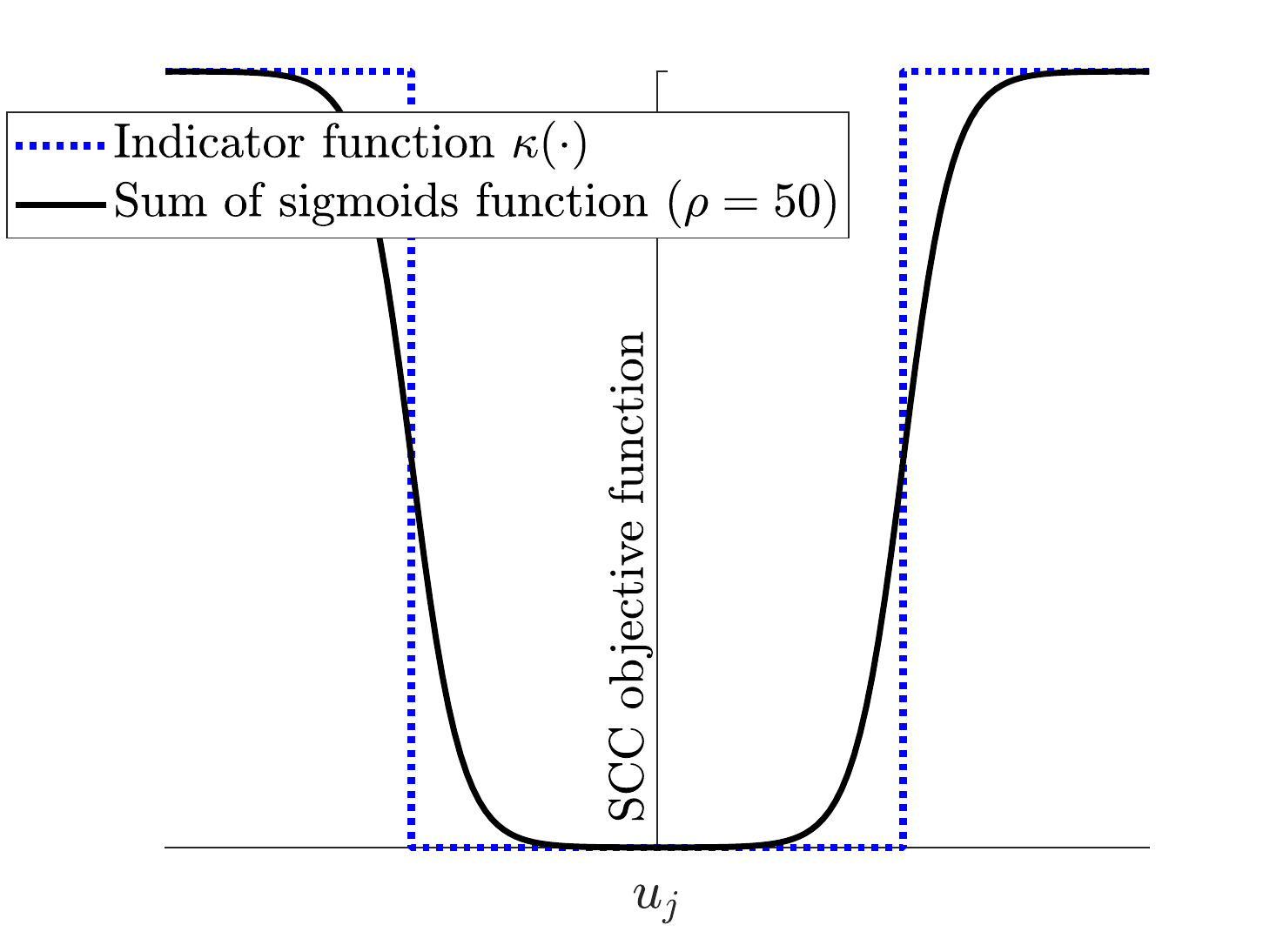}
    \caption{SCC indicator function $\kappa(\cdot)$ and its continuous sum of sigmoids approximation $f_{\widetilde{\text{SCC}}}$ for link $j$}
    \label{fig:SCC_obj}
    \end{center}
\end{figure}

The SCC objective function $f_{\text{SCC}}$ is nonsmooth at $\pm u^{\min}$, resulting in unbounded gradients. Therefore, in order to employ gradient-based optimization methods, the threshold function $\kappa(\cdot)$ is approximated with a continuous sum of sigmoids (or logistic) function, as proposed in \citet{ABRAHAM2016}. The sigmoidal function has positive and negative components, defined by $\psi^+_j(u) := \left(1+e^{-\rho(u-u^{\min}_j)}\right)^{-1}$ and $\psi^-_j(u) := \left(1+e^{-\rho(-u-u^{\min}_j)}\right)^{-1}$, respectively, where $\rho$ is a parameter which sets the sigmoid function curvature. The following expression combines these sigmoid functions to approximate $f_{\text{SCC}}$ posed in \eqref{eq:SCC_objective}:
\begin{equation}
\label{eq:SCC_objective_sigmoid}
    f_{\widetilde{\text{SCC}}} := \frac{1}{n_t} \sum_{t=1}^{n_t} \sum_{j=1}^{n_p}w_j \left(\psi^+_{j}\bigg(\frac{q_{j,t}}{A_j}\bigg) + \psi^-_j\bigg(\frac{q_{j,t}}{A_j}\bigg)\right).
    \vspace{0.1cm}
\end{equation}
where velocity is defined as $u_{j,t} = \big(\frac{q_{j,t}}{A_j}\big)$. In accordance with that reported in \citet{ABRAHAM2016}, we found $\rho \leq 100$ provided a step-like objective function, whilst still being sufficiently smooth at the threshold boundaries for gradients to exist. An example of the SCC indicator function $\kappa(\cdot)$ and its continuous sum of sigmoids approximation $f_{\widetilde{\text{SCC}}}$ are shown in \Cref{fig:SCC_obj}.

\subsection{Mixed-integer nonlinear program}
The SCC design-for-control problem aims to maximize \eqref{eq:SCC_objective_sigmoid}, subject to hydraulic conservation laws and physical and economical valve constraints. The problem formulation includes continuous variables denoted by $\bm{x}:=[\bm{q} \, \bm{h} \, \bm{\eta} \, \bm {\theta} \, \bm {\alpha}]^T$ and binary variables denoted by $\bm{v}:=[\bm{v^+} \, \bm{v^-}]^T$, $\bm{y}$ and $\bm{z}$. Here, the objective function is replaced with its additive inverse to result in a minimization optimization problem. The resulting mixed integer nonlinear program (MINLP) is summarized by the following problem formulation.
\begin{equation} \tag{$\text{MINLP}$}
\label{eq:MINLP_problem}
\begin{aligned}
& \underset{\substack{\bm{x}, \, \bm{v}, \, \bm{y}, \, \bm{z}}}{\text{minimize}}
& & -f_{\widetilde{\text{SCC}}}\;\: \eqref{eq:SCC_objective_sigmoid} \\
& \text{subject to}
& & \text{linear hydraulic conservation constraints}\;\: \eqref{eq:hydraulic_conservation_b}\ \text{and}\  \eqref{eq:hydraulic_conservation_c}\\
& & & \text{nonconvex HW head loss model constraints}\;\: \eqref{eq:hydraulic_conservation_d} \\
& & & \text{big-M constraints for control valve operation}\;\: \eqref{eq:valve_bigM_a} -  \eqref{eq:valve_bigM_f}\\
& & & \text{big-M constraints for flushing valve operation}\;\: \eqref{eq:valve_bigM_g} \\
& & & \text{physical and economical valve constraints}\;\: \eqref{eq:valves_physical_a} -  \eqref{eq:valves_physical_c}\\
& & & \bm{x} \in \bm{Q}\;\: \eqref{eq:hyd_bounds_a} -  \eqref{eq:hyd_bounds_e}\\
& & & \bm{v} \in \{0,1\}^{2n_pn_t}, \; \bm{y} \in \{0,1\}^{n_n}, \; \bm{z} \in \{0,1\}^{n_p}\\
\end{aligned}
\end{equation}
where $\bm{Q}$ is a rectangle representing upper and lower bounds for the continuous decision variables. The continuous decision variables are defined as: $\bm{q} := (q_t)_{t=1,\dots,n_t}$, $\bm{h} := (h_t)_{t=1,\dots,n_t}$, $\bm{\eta} := (\eta_t)_{t=1,\dots,n_t}$, $\bm{\theta} := (\theta_t)_{t=1,\dots,n_t}$ and $\bm{\alpha} := (\alpha_t)_{t=1,\dots,n_t}$. The binary decision variables varying with time step $t$ are defined as: $\bm{v^+} := (v^+_t)_{t=1,\dots,n_t}$ and $\bm{v^-} := (v^-_t)_{t=1,\dots,n_t}$.

Problem \eqref{eq:MINLP_problem} has $n_t\left(3n_p+2n_n\right)$ continuous variables, $2n_tn_p+n_p+nn$ binary variables and $2n_tn_p$ nonconvex terms. Observe that the problem grows rapidly with the size of the considered WDN (see \Cref{table:prob_data}), making it a difficult nonconvex MINLP problem to solve. To overcome these challenges, we develop a convex heuristic to compute feasible solutions to Problem \eqref{eq:MINLP_problem}. The following section describes the solution algorithm and its implementation details.

\section{Solution method}
\label{sec:sol_methods}
\begin{figure}[!h] 
\begin{center}
\captionsetup{justification=centering}
\includegraphics[width=0.65\linewidth]{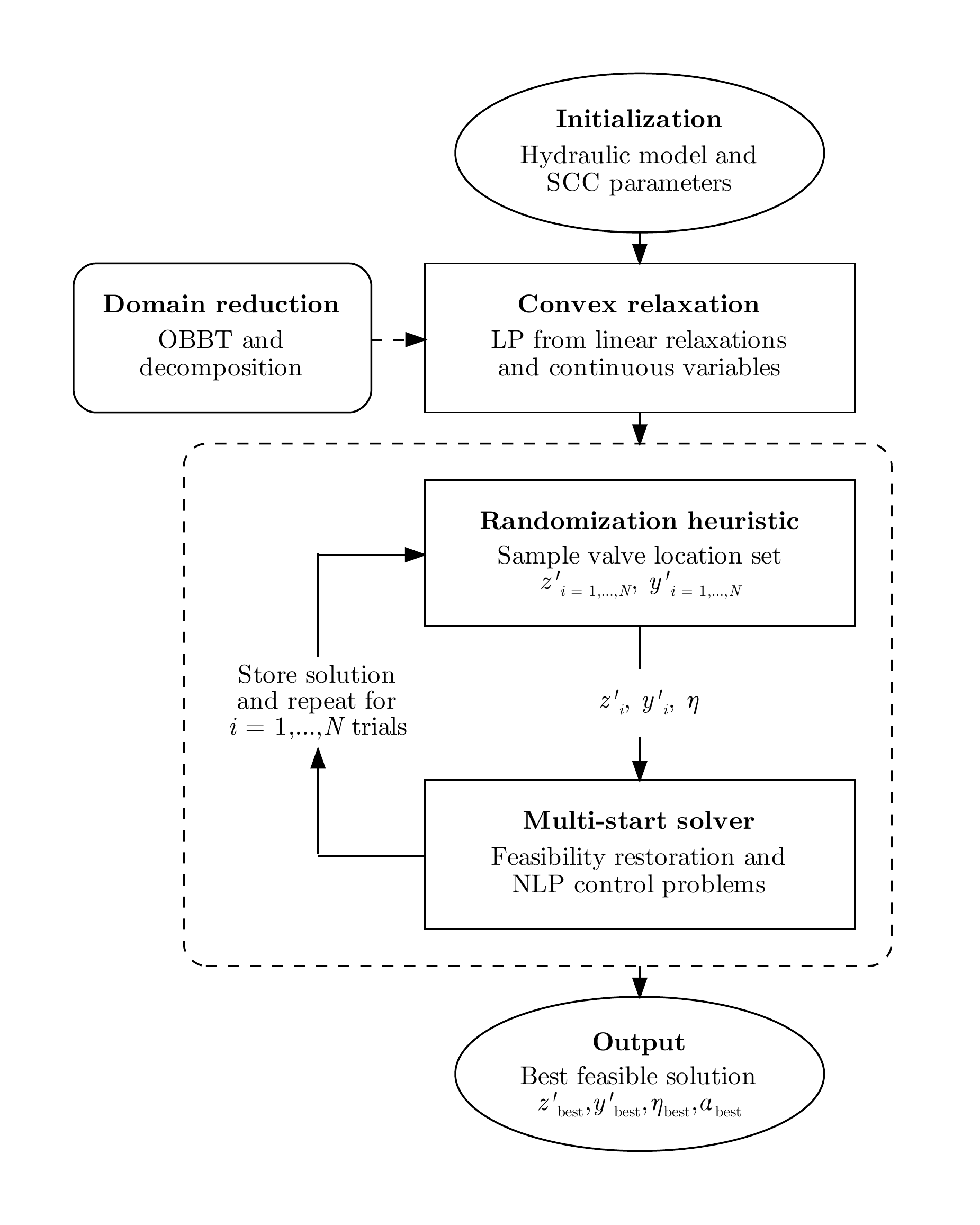}
\caption{Proposed convex multi-start (CMS) algorithm}
\label{fig:solution_algorithm_chart}
\end{center}
\end{figure}

The proposed solution algorithm combines convex relaxations with a randomization heuristic and multi-start solver to compute feasible solutions to Problem \eqref{eq:MINLP_problem}. First, we formulate convex relaxations of Problem \eqref{eq:MINLP_problem}, which yield a linear programming (LP) subproblem. Here, the nonconvex SCC objective function \eqref{eq:SCC_objective} and nonconvex energy conservation constraints \eqref{eq:hydraulic_conservation_d} are relaxed using polyhedral envelopes, and a continuous relaxation is applied to binary decision variables. We also include a domain reduction step where the resulting convex subproblem is tightened using both model decomposition and optimization-based bound tightening (OBBT) techniques. Then, a randomization heuristic uses the fractional valve placement values from the convex subproblem solution to form nonlinear programming (NLP) control problems. Local solutions to these NLP problems are computed using a strictly feasible sequential convex programming (SFSCP) solver. We also implement a multi-start strategy, which includes an optimization-based feasibility restoration problem to ensure hydraulic feasibility of the starting points. \Cref{fig:solution_algorithm_chart} offers a detailed overview of the solution method, referred to as the convex multi-start (CMS) algorithm. The algorithm steps and overall implementation details are provided in the following subsections.

\subsection{Convex relaxation}
\label{sec:convex_relax}
There are two sources of nonlinear nonconvexity that make Problem \eqref{eq:MINLP_problem} difficult to solve. These are the SCC objective function $f_{\widetilde{\text{SCC}}}$ \eqref{eq:SCC_objective_sigmoid} and the HW head loss model $\phi(\cdot)$ in equality constraint \eqref{eq:hydraulic_conservation_d}. One approach to overcome nonconvexity is through polyhedral relaxations, which can be formulated as linear constraints and thus efficiently handled by state-of-the-art linear solvers. Since the resulting mixed integer linear program (MILP) may still have a large number of binary variables, it is convenient to apply a continuous relaxation to binary variables. In addition to the computational advantages, a continuous relaxation increases the search space for the set of optimal binary variables. Here, we implement the aforementioned relaxation techniques to formulate a convex subproblem of Problem \eqref{eq:MINLP_problem}, with its solution forming the basis of the heuristic algorithm.

\begin{figure}[h!]
    \centering
    \captionsetup{justification=centering}
    \subfloat[\label{fig:convex_relax_a}\text{Sum of sigmoids objective function $\psi(u_j)$}]{
        \includegraphics[width=0.44\textwidth]{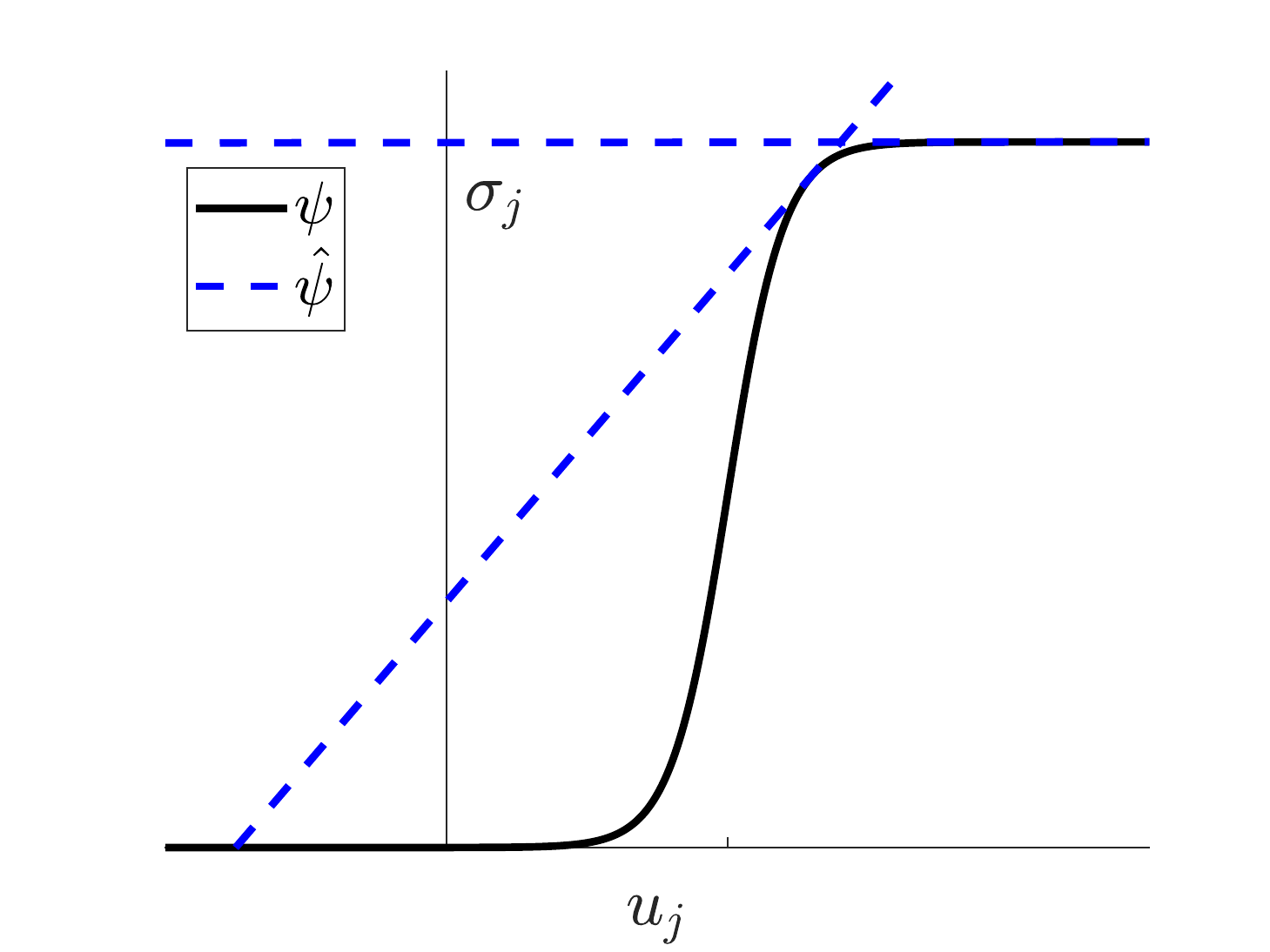}}
        \hspace{-0.15cm}
    \subfloat[\label{fig:convex_relax_b}\text{HW head loss model $\phi(q_j)$}]{
        \includegraphics[width=0.44\textwidth]{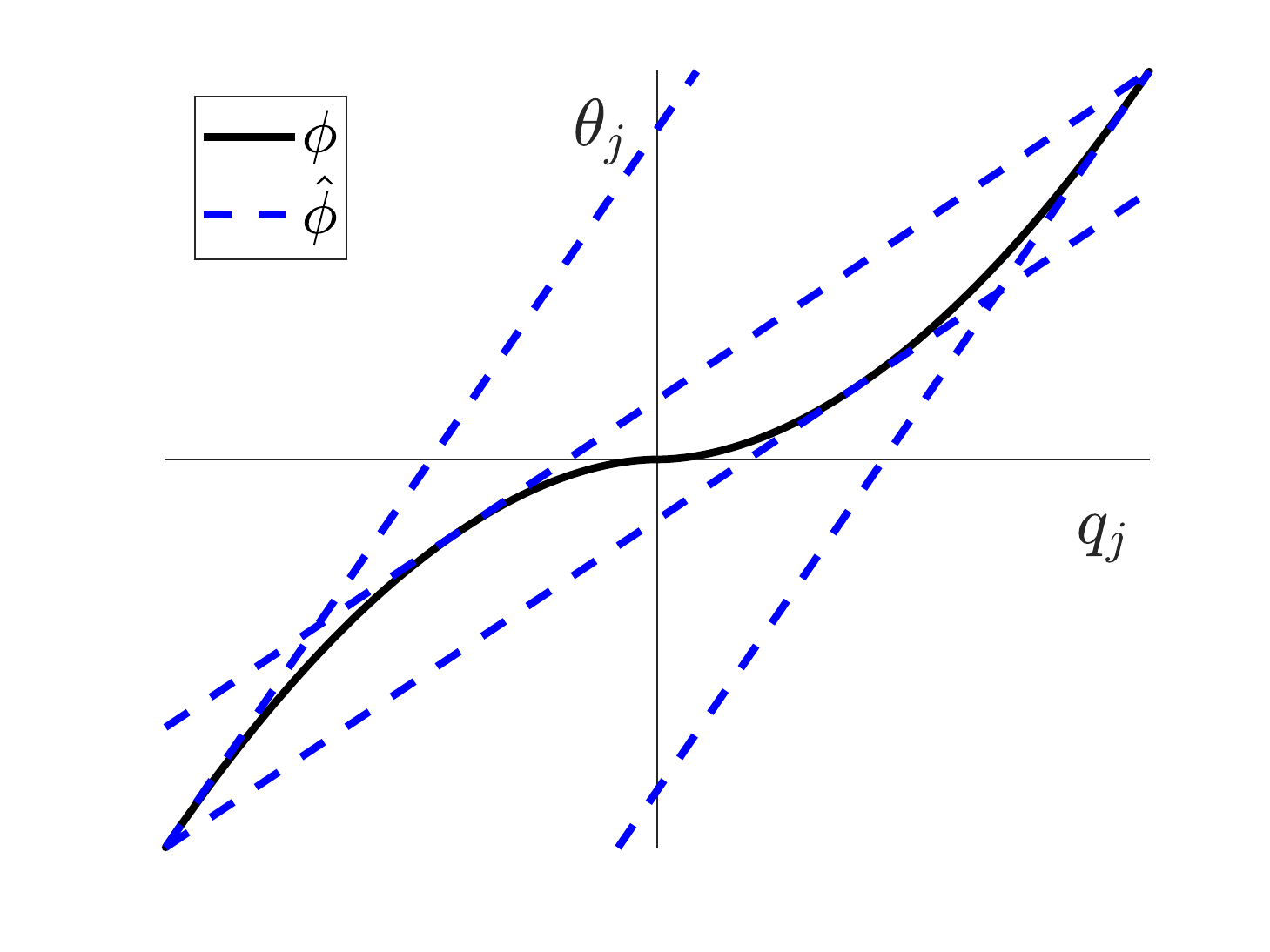}}
    \\ \vspace{-0.1cm}
    \caption{Examples of implemented linear relaxations for link $j$}
    \label{fig:convex_relax}
\end{figure}

We first relax the nonconvex objective function in \eqref{eq:SCC_objective_sigmoid} with a linear outer approximation. Let $\sigma_t^+ \in \mathbb{R}^{n_{p}}$ and $\sigma_t^- \in \mathbb{R}^{n_{p}}$ be vectors of auxiliary variables introduced to model the positive $\psi^+$ and negative $\psi^-$ sigmoid functions, respectively. The objective function $f_{\widetilde{\text{SCC}}}$ is then reformulated as the following set of inequality constraints:
\begin{subequations}
    \label{eq:SCC_obj_relax_a}
    \begin{align}
        \label{eq:SCC_obj_relax_a1}
        &\sigma^+_t \leq \psi^+\bigg(\frac{q_{t}}{A}\bigg), \quad t \in \{1,\ldots,n_t\}\\[3pt]
        \label{eq:SCC_obj_relax_a2}
        & \sigma^-_t \leq \psi^-\bigg(\frac{q_{t}}{A}\bigg), \quad t \in \{1,\ldots,n_t\}.
    \end{align}
\end{subequations}
To ensure equivalence with \eqref{eq:SCC_objective_sigmoid}, the SCC objective function becomes
\begin{equation}
    \label{eq:SCC_obj_relax_b}
    f_{\widetilde{\text{SCC}}} := \frac{1}{n_t} \sum_{t=1}^{n_t} \sum_{j=1}^{n_p}w_j \left(\sigma^+_{j,t} + \sigma^-_{j,t}\right).
\end{equation}
We construct concave envelopes for the positive $\psi^+$ and negative $\psi^-$  components of $f_{\widetilde{\text{SCC}}}$. This follows the methodology presented in \citet{UDELL2014,UDELL2016} for sigmoidal functions, resulting in piecewise linear relaxations. These relaxations are written as the following constraint:
\begin{equation}
    \label{eq:SCC_obj_relax_c}
    S_tq_t + T_t\sigma_t \le s_t, \quad \forall t \in \{1,\ldots,n_t\} 
\end{equation}
where matrices $S_t := [S_t^+ \, S_t^-]^T$ and $T_t := [T_t^+ \, T_t^-]^T$ and vector $s_t := [s_t^+ \, s_t^-]^T$ depend on flow velocity bounds $u_t^L := \big(\frac{q_t^L}{A}\big)$ and $u_t^U := \big(\frac{q_t^U}{A}\big)$ as well as the sigmoid function parameters. A detailed derivation of these relaxations is provided in Appendix A.1. Moreover, an example of the implemented relaxation is illustrated in \Cref{fig:convex_relax_a}, with $\hat{\psi}$ denoting the set of linear relaxations.

We then implement polyhedral relaxations for the HW head loss model $\phi(\cdot)$ in equality constraint \eqref{eq:hydraulic_conservation_d}. This builds on the relaxation methods for monomials of odd degree introduced by \citet{LIBERTI2003} by formulating polyhedral relaxations for the HW head loss model. Similar to \eqref{eq:SCC_obj_relax_c}, the formulated relaxations are written as the following linear constraint:
\begin{equation}
    \label{eq:HW_model_relax}
    R_tq_t + E_t\theta_t \le r_t, \quad \forall t \in \{1,\ldots,n_t\} 
\end{equation}
where matrices $R_t$ and $E_t$ and vector $r_t$ are derived from flow bounds $q_t^L$ and $q_t^U$ and the HW model parameters. Further details are presented in Appendix A.2 and an example of the implemented relaxation is illustrated in \Cref{fig:convex_relax_b}, noting that the polyhedral relaxation is denoted by $\hat{\phi}$.

Lastly, we implement the continuous relaxation of binary variables $\bm{v}$, $\bm{y}$, and $\bm{z}$. These are combined with continuous decision variables $\bm{x}$, defined in Problem \eqref{eq:MINLP_problem}, and the continuous auxiliary variables associated with the $f_{\widetilde{\text{SCC}}}$ inequality constraints \eqref{eq:SCC_obj_relax_a1} and \eqref{eq:SCC_obj_relax_a2}, denoted by $\bm{\sigma}:=[\bm{\sigma^+}\bm{\sigma^-}]$. The resulting convex relaxation of Problem \eqref{eq:MINLP_problem} is represented by the following LP subproblem.
\begin{equation} \tag{$\text{LP}$}
\label{eq:LP_problem}
\begin{aligned}
& \underset{\substack{\bm{x}, \, \bm{v}, \, \bm{y}, \, \bm{z}, \, \bm{\sigma}}}{\text{minimize}}
& & -f_{\widetilde{\text{SCC}}}\;\: \text{reformulated in}\ \eqref{eq:SCC_obj_relax_b} \\
& \text{subject to}
& & \text{linear relaxations of HW head loss model constraints}\;\: \eqref{eq:HW_model_relax} \\
& & & \text{linear relaxations of $\psi^+$ and $\psi^-$ sigmoid functions}\;\: \eqref{eq:SCC_obj_relax_c}\\
& & & \text{linear hydraulic conservation constraints}\;\: \eqref{eq:hydraulic_conservation_b}\ \text{and}\  \eqref{eq:hydraulic_conservation_c}\\
& & & \text{big-M constraints for control valve operation}\;\: \eqref{eq:valve_bigM_a} -  \eqref{eq:valve_bigM_f}\\
& & & \text{big-M constraints for flushing valve operation}\;\: \eqref{eq:valve_bigM_g} \\
& & & \text{physical and economical valve constraints}\;\: \eqref{eq:valves_physical_a} -  \eqref{eq:valves_physical_c}\\
& & & \bm{x} \in \bm{Q}\;\: \eqref{eq:hyd_bounds_a} -  \eqref{eq:hyd_bounds_e}, \; \bm{\sigma} \in [0,1]^{2n_pn_t}\\
& & & \bm{v} \in [0,1]^{2n_pn_t}, \; \bm{y} \in [0,1]^{n_n}, \; \bm{z} \in [0,1]^{n_p}\\
\end{aligned}
\end{equation}
The optimal value to Subproblem~\eqref{eq:LP_problem} yields a lower bound to the original nonconvex Problem \eqref{eq:MINLP_problem}. In particular, vectors $y \in \mathbbm{R}^{n_n}$ and $z \in \mathbbm{R}^{n_p}$ of continuous variables for valve placement can be interpreted as the probability distribution from which random samples are drawn to solve the NLP control problem. Moreover, the vector $\eta \in \mathbbm{R}^{n_p}$ from Subproblem \eqref{eq:LP_problem} is used as one of $M$ starting points in the multi-start strategy. These steps are discussed in the subsequent sections.

Furthermore, we implement a domain reduction procedure to reduce the bound intervals and thus strengthen the convex relaxations formulated for Subproblem \eqref{eq:LP_problem}. The main procedure implemented is an optimization-based bound tightening (OBBT) algorithm \citep{BELOTTI2009}. The OBBT algorithm solves a series of LP optimization problems, with objective functions set to both maximize and minimize link flow. This is augmented with a forest-core decomposition scheme to reduce the number of flow variables whose bounds are tightened \citep{SIMPSON2014}. The forest of a WDN comprises the disjoint union of all outer branch (or tree) components \citep{DEUERLEIN2008}. The core represents the set of looped (or block) graphs, which contain the roots of all forest trees. Following \citet{PECCI2019}, we only perform bound tightening on the set of core links. We include psuedocode for the OBBT algorithm in Appendix B.

\subsection{Randomization heuristic}
In searching for a good quality local solution to Problem \eqref{eq:MINLP_problem}, we employ a randomization heuristic to sample candidate valve configurations. Here, the fractional values of $y \in \mathbbm{R}^{n_n}$ and $z \in \mathbbm{R}^{n_p}$ yielded from the solution to Subproblem \eqref{eq:LP_problem} define a discrete probability distribution for the index sets $\{1,\dots,n_n\}$ and $\{1,\dots,n_p\}$, respectively. We propose to randomly sample candidate pressure control $n_v$ and flushing $n_f$ valve configurations from the respective probability distributions over $N$ sampling trials. This creates sampled vectors of binary variables $y'_{i=1,\dots,N} \in \{0,1\}^{n_n}$ and $z'_{i=1,\dots,N} \in \{0,1\}^{n_p}$, from which a local solution to the NLP control problem is obtained for each trial $i \in \{1,\dots,N\}$ (see \Cref{sec:multistart_solver}). To avoid redundant valve configurations, we store binary values from each trial in the set $\mathcal{P} \in \mathbbm{R}^{N \times (n_v+n_f)}$, checking its intersection with $\{y'_i \cup z'_i\}$ before fixing valve placement values and proceeding to the NLP control problem. Psuedocode for the implemented randomization heuristic is detailed in \Cref{alg:randomization heuristic}.

\begin{algorithm}
    \caption{Randomization heuristic}
    \label{alg:randomization heuristic}
    \begin{algorithmic}[1]
    \State \textbf{Input:} vectors of fractional values $y \in [0,1]^{n_n}$ and $z \in [0,1]^{n_p}$ \Comment{Problem \eqref{eq:LP_problem}}
    \State \textbf{Output:} sampled vectors of binary variables $y'_{i=1,\dots,N} \in \{0,1\}^{n_n}$ and $z'_{i=1,\dots,N} \in \{0,1\}^{n_p}$
    \State Initialize $\mathcal{P} \;\leftarrow\; \emptyset$
    \For{$i = 1,\dots,N$}
    \While{$\{y'_i \cup z'_i \} = \emptyset$ \textbf{or} $\{y'_i \cup z'_i\} \in \mathcal{P}$}
    \State Sample $y'_i$ from vector $y$ with probability weights equal to fractional values
    \State Sample $z'_i$ from vector $z$ with probability weights equal to fractional values
    \EndWhile
    \State Update visited valve locations: $\mathcal{P}_{(i)} \;\leftarrow\; \{y'_i \cup z'_i\}$
    \State Vectors $y'_i$ and $z'_i$ \Comment{Input to multi-start solver}
    \EndFor
    \end{algorithmic}
\end{algorithm}

\subsection{Multi-start solver}
\label{sec:multistart_solver}
\begin{figure}[!h] 
    \vspace{0.25cm}
    \captionsetup{justification=centering}
    \includegraphics[width=01.0\linewidth]{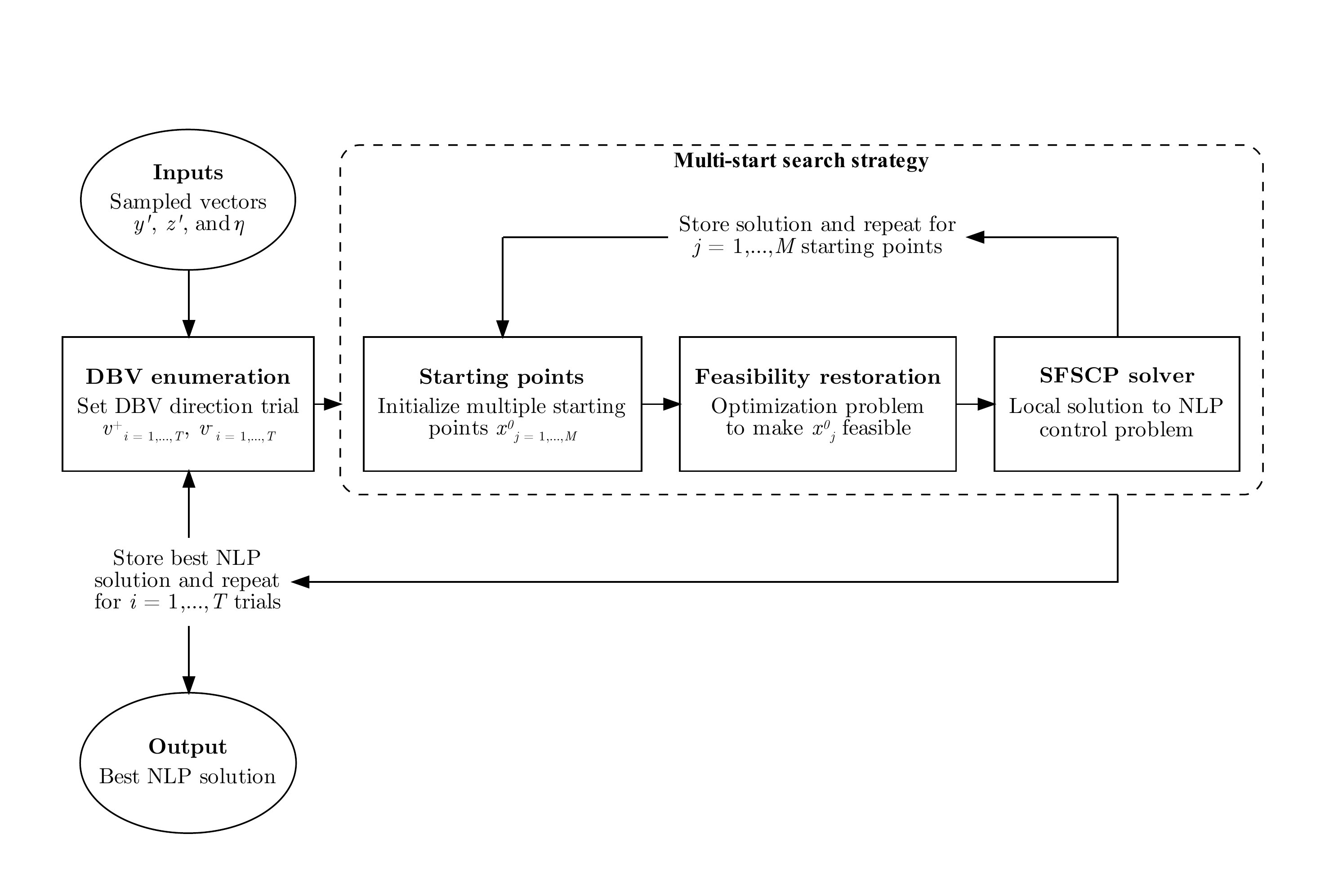}
    \caption{Process diagram of multi-start solver}
    \label{fig:solution_multistart}
\end{figure}

The final component of the CMS solution algorithm concerns the optimization of valve settings. This control problem is formulated to compute (locally) optimal operational settings for each candidate valve configuration generated from the randomization heuristic. Because the objective function in Problem \eqref{eq:MINLP_problem} is highly nonlinear, we implement a strategy to use multiple starting points in order to avoid getting trapped in poor local optima. This includes a feasibility restoration subproblem to ensure hydraulically feasible starting points are passed to the NLP solver. We refer to the overall solution process as the multi-start solver, which has steps highlighted by the process diagram shown in \Cref{fig:solution_multistart}.

\subsubsection{Strictly feasible sequential convex programming solver}
With the placement of DBVs and AFVs from the set of sampled configurations $\{1,\dots,N\}$, we then fix binary flow direction variables $v^+_t$ and $v^-_t$ at DBV locations for each time step $t \in \{1,\dots,n_t\}$. This yields a nonconvex NLP problem with continuous decision variables denoted by $\bm{x}:=[\bm{q} \, \bm{h} \, \bm{\eta} \, \bm {\alpha}]^T$. Moreover, big-M constraints for pressure control and flushing valve operations found in Problem \eqref{eq:MINLP_problem} are reformulated as variable bounds. These bounds are derived on the basis of the fixed binary variables. The corresponding NLP control problem is formulated as follows.
\begin{equation} \tag{$\text{NLP}$}
\label{eq:NLP_problem}
\begin{aligned}
& \underset{\substack{\bm{x}}}{\text{minimize}}
& & -f_{\widetilde{\text{SCC}}}\;\: \eqref{eq:SCC_objective_sigmoid} \\
& \text{subject to}
& & \text{hydraulic conservation constraints}\;\: \eqref{eq:hydraulic_conservation_a}\ \text{and}\  \eqref{eq:hydraulic_conservation_b}\\
& & & \bm{x} \in \bm{\tilde{Q}}\\
\end{aligned}
\end{equation}
where $\bm{\tilde{Q}}$ represents the upper and lower variable bounds modified to account for directional constraints enforced by the fixed binary variables. Since the focus of this study is at the DMA- or distribution-level, which in practice limits the number of installed DBVs for consideration, we use an enumeration approach to test all DBV flow direction combinations for $v^+_t$ and $v^-_t$. This approach was deemed appropriate on the basis of the relatively small computational (CPU) times recorded in \Cref{sec:results}. For example, the numerical experiment corresponding to the maximum number of installed DBVs on the largest network tested in this study (see \Cref{sec:results_setup} for details) resulted in a CPU time of less than five minutes. This CPU time is well within an assumed one hour limit for application in control strategies. We acknowledge, however, that this assumption warrants re-evaluation for system-wide studies which inevitably consider a larger set of control valves. Here, a total of $T = 2^{n_{\text{DBV}}}$ control problems are solved to find the best DBV direction combination at each time step $t \in \{1,\dots,n_t\}$ (see \Cref{fig:solution_multistart}).   
In this work, we use the strictly feasible sequential convex programming (SFSCP) method outlined in \citet{WRIGHT2015} for solving a local optima to Problem \eqref{eq:NLP_problem}. The SFSCP solver, also used to optimize SCC control valve settings in \citet{ABRAHAM2016}, is implemented because of its reliable convergence properties, namely its strict feasibility requirement for each iterate of the optimization process. The SFSCP solver linearizes the objective function and nonlinear hydraulic conservation constraints in Problem \eqref{eq:NLP_problem}, from which a sequence of linear programs can be efficiently solved. Linear approximations are formulated using first-order Taylor expansions around the point $x_k := [q_k\, h_k\, \eta_k\, \alpha_k]^T$. Let $f := f_{\widetilde{\text{SCC}}}$ and the hydraulic conservation equality constraints \eqref{eq:hydraulic_conservation_a} and \eqref{eq:hydraulic_conservation_b} be denoted by $g(\cdot)$. It follows that
\begin{equation}
    \label{eq:SCC_obj_linearize}
    f(x) \approx f(x_k) + \frac{\partial{f}}{\partial{q}}(x_k)(q-q_k),
\end{equation}
and
\begin{equation}
    \label{eq:hydraulic_linearize}
    g(x) \approx g(x_k) + \frac{\partial{g}}{\partial{q}}(x_k)(q-q_k) + \frac{\partial{g}}{\partial{\eta}}(x_k)(\eta-\eta_k) + \frac{\partial{g}}{\partial{\alpha}}(x_k)(\alpha-\alpha_k).
\end{equation}
are the linearized expressions used to formulate a convex subproblem for each iterate $k$. This subproblem computes step directions $d\eta_k$ and $d\alpha_k$ at point $x_k$ for pressure control and flushing valve operational variables, respectively. Following \citet{WRIGHT2015}, a strictly feasible line search method is implemented to guarantee an improvement in objective function whilst ensuring hydraulic feasibility of the computed local optima. We define the hydraulically feasible region by the set
\begin{equation}
    \label{eq:feas_region}
    \mathcal{F} := \bigg\{x \in \mathbbm{R}^{n_p+n_n+n_v+n_f} \;\big| \; x \; \text{satisfies constraints } \eqref{eq:hydraulic_conservation_a}, \eqref{eq:hydraulic_conservation_b}, \text{ and variable bounds} \; \tilde{Q}\bigg\}.
\end{equation}
The null space solver proposed in \citet{ABRAHAM2015} is used to solve hydraulic states $q_{k+1}$ and $h_{k+1}$ for the optimal controls obtained at iterate $k$. If $x_k \in \mathcal{F}$, then the optimization process moves to the next iterate $ k := k+1$. Otherwise, a backtracking procedure is implemented to reduce $d\eta_k$ and $d\alpha_k$ until hydraulic feasibility is achieved. Termination of the SFSCP solver is triggered when relative improvements in the objective function are below a set tolerance $\epsilon_{\text{tol}}$. The SFSCP solver steps are detailed in psuedocode presented in Appendix C \citep{WRIGHT2015}. Moreover, the subsequent section describes a multi-start strategy for selecting feasible starting points $x^0$ to be passed to the SFSCP solver. 

For completeness, we provide a brief discussion on the performance comparison with state-of-the-art nonlinear solver IPOPT. This is included as many large-scale WDN optimization studies have used IPOPT to compute locally optimal solutions for network controls (e.g. \citeauthor{PECCI2019}, \citeyear{PECCI2019}; \citeauthor{ULUSOY2020}, \citeyear{ULUSOY2020}; \citeauthor{NERANTZIS2022}, \citeyear{NERANTZIS2022}). Most of these studies, however, have applied a quadratic approximation (QA) to model the head loss function $\psi(\cdot)$ as the nonsmooth HW formula presents challenges for IPOPT, namely that its Hessian is unbounded at the origin. Preliminary efforts to introduce a QA model for the current SCC problem resulted in large inaccuracies compared to the HW model. This was expected as the current study aims to maximize pipe flow velocities, which increases the QA model error since a larger maximum flow value is required (\citeauthor{PECCI2017a}, \citeyear{PECCI2017a}; Equation 21). Therefore, we tested IPOPT using the HW model and relied on its feasibility restoration scheme to mitigate potential errors. Infeasibility errors were observed for the larger case study network when using the IPOPT solver, ultimately leading to the SFSCP solver as the more reliable and robust option for the current SCC control problem. The SFSCP and IPOPT solver comparison is detailed in Appendix D.

\subsubsection{Multi-start strategy}
As highlighted in \Cref{fig:solution_multistart}, a multi-start strategy is implemented to find good quality local optima to Problem \eqref{eq:NLP_problem}. Here, we initialize starting points $x^0$ for multiple initial conditions $j \in \{1,\dots,M\}$, which are sequentially passed to the SFSCP solver. Starting points are formed by simulating hydraulic states for different initial control valve settings $\eta^0$. For initial condition $j = 1$, we set $\eta^0_{j=1}$ equal to the control valve settings resulting from Subproblem \eqref{eq:LP_problem}. Otherwise, $\eta^0_{j=2,\dots,M}$ is randomly generated between the previously set lower and upper $\eta$ bounds. Note that initial flushing operations are set as $\alpha^0 = \{0\}$ for all initial conditions $j \in \{1,\dots,M\}$. 
Since the initial control conditions do not guarantee a solution in the hydraulically feasible region, we also include a feasibility restoration step. That is, if starting point $x^0_j \notin \mathcal{F}$ for all $j \in \{1,\dots,M\}$, feasibility is restored through the solution to the following NLP optimization problem.
\begin{equation} \tag{$\text{FR}$}
\label{eq:FR_problem}
\begin{aligned}
& \underset{\substack{\bm{x}}}{\text{minimize}}
& & \| \eta - \eta^0\|_2^2\\
& \text{subject to}
& & \text{hydraulic conservation constraints}\;\: \eqref{eq:hydraulic_conservation_a}\ \text{and}\  \eqref{eq:hydraulic_conservation_b}\\
& & & \bm{x} \in \bm{\tilde{Q}}\\
\end{aligned}
\end{equation}
where the objective function is set to minimize the squared $\ell_{2}$-norm between the vector of initial settings $\eta^0$ and settings $\eta$ which produce a hydraulically feasible starting point. Problem \eqref{eq:FR_problem} is solved using the nonlinear IPOPT solver. Apart from cases where the fixed DBV directions strictly cannot facilitate a hydraulically feasible solution, we did not encounter infeasibility errors when using IPOPT as the NLP solver for Problem \eqref{eq:FR_problem}. The discussion in Appendix D suggests that the nonconvexity of the SCC objective function in Problem \eqref{eq:NLP_problem} is a key factor for the observed infeasibility issues.

\section{Numerical experiments and discussion}
\label{sec:results}
We evaluate the SCC design-for-control problem using the proposed CMS solution algorithm. Since CMS is a heuristic, we also compare its performance with an off-the-shelf GA implementation, which is a common heuristic method for solving WDN design problems. This section is organized as follows. First, in \Cref{sec:results_setup} we describe the problem parameters, case study networks, and computational resources used to carry out the numerical experiments. Next, \Cref{sec:results_control} discusses the results obtained when solving the control-only subproblem to Problem \eqref{eq:MINLP_problem}, where valves locations are fixed. In this case, operational settings of existing PRV and DBV (if applicable) valve configurations are optimized to maximize self-cleaning conditions. We then investigate the solution to the design-for-control problem in \Cref{sec:results_dfc}, where the placement and operation of both pressure control and flushing valves are considered to improve SCC performance.   

\subsection{Computational setup}
\label{sec:results_setup}
All numerical experiments were performed in MATLAB R2021b (64-bit) for Microsoft Windows 11, installed on a 2.50-GHz Intel(R) Core(TM) i9-11900H CPU with 8 cores and 32.0 GB of memory (RAM). The LP problems involved in the convex relaxation, OBBT algorithm, and SFSCP solver were solved using Gurobi (v9.5.0; \citeauthor{GUROBI2022}, \citeyear{GUROBI2022}), accessed via its MATLAB interface. With the exception of large problem cases, which had solver parameters (e.g. NumericFocus, Crossover, and Presolve) modified to deal with numerical issues, Gurobi's default parameters were applied. Local solutions to the NLP feasibility restoration problems were solved using the state-of-the-art interior point solver IPOPT (v3.12.19; \citeauthor{WACHTER2006}, \citeyear{WACHTER2006}). IPOPT was accessed in MATLAB via the OPTI toolbox interface \citep{CURRIE2012} and was implemented passing the exact Jacobian and Hessian matrices as inputs to the solver. The single objective GA was implemented using MATLAB's off-the-shelf Global Optimization Toolbox with its default solver parameters (e.g. population size $=$ 100; crossover fraction $=$ 0.8; function tolerance $=$ $1\mathrm{e}{-6}$). For each GA fitness evaluation, EPANET2.2 \citep{EPANET2.2} hydraulic solver was accessed via the EPANET-MATLAB Toolkit interface \citep{ELIADES2016}. We set maximum time limits for the GA implementation associated with control and design-for-control problems of 6 and 12 hours, respectively. Furthermore, MATLAB's parallel computing toolbox was employed (with eight local workers) to speed up code in the OBBT algorithm and multi-start solver. Finally, network hydraulics were computed within CMS using the null space method proposed in \citet{ABRAHAM2015}.

\begin{figure}[h!]
    \centering
    \captionsetup{justification=centering}
    \subfloat[\label{fig:netLayouts_a}\text{\texttt{Pescara}}]{
        \includegraphics[width=0.43\textwidth]{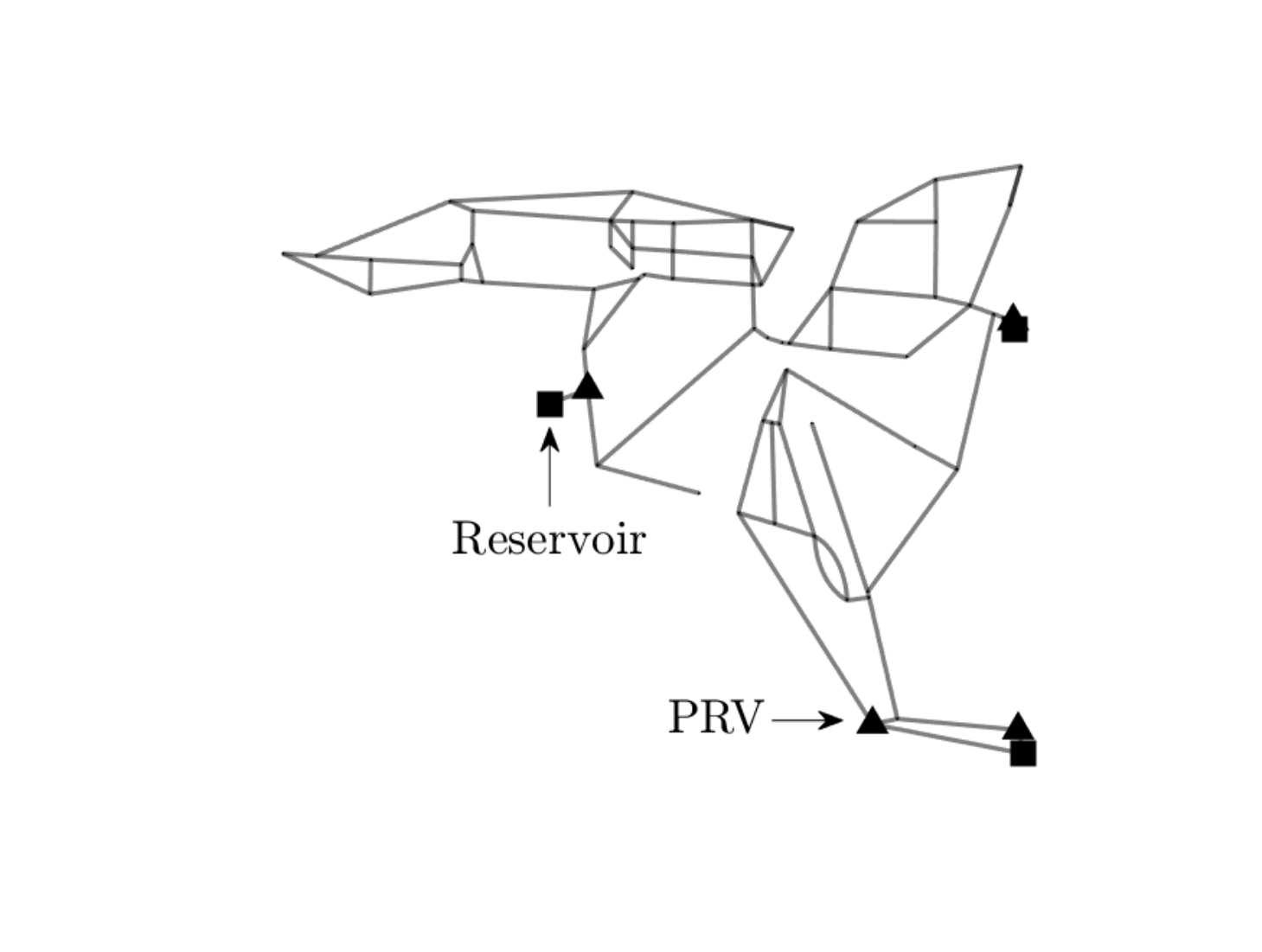}}
    \hspace{0.4cm}
    \subfloat[\label{fig:netLayouts_b}\text{\texttt{Modena}}]{
        \includegraphics[width=0.43\textwidth]{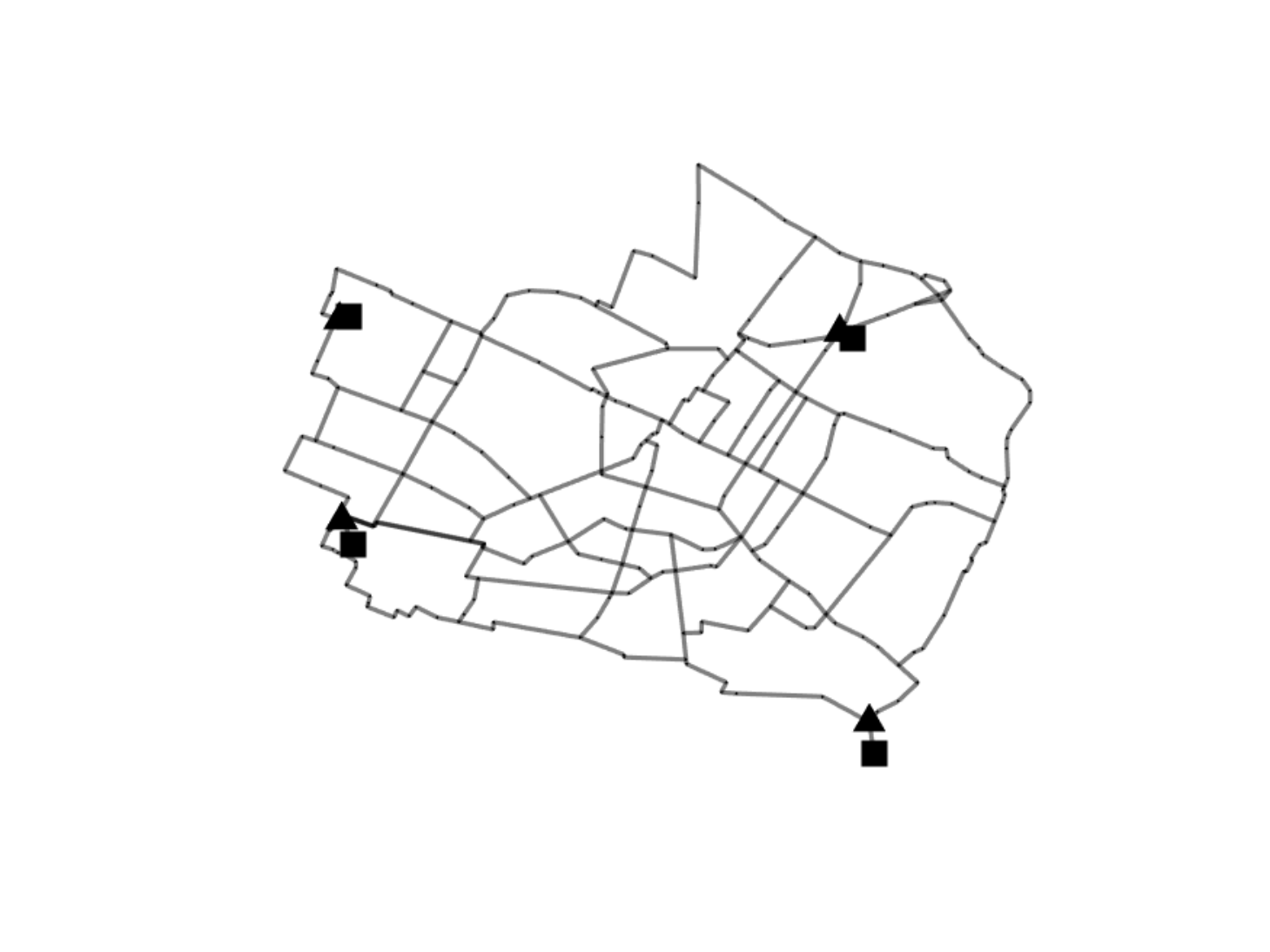}}
    \vspace{0.4cm}    
    \subfloat[\label{fig:netLayouts_c}\text{\texttt{BWFLnet}}]{
        \includegraphics[width=0.45\textwidth]{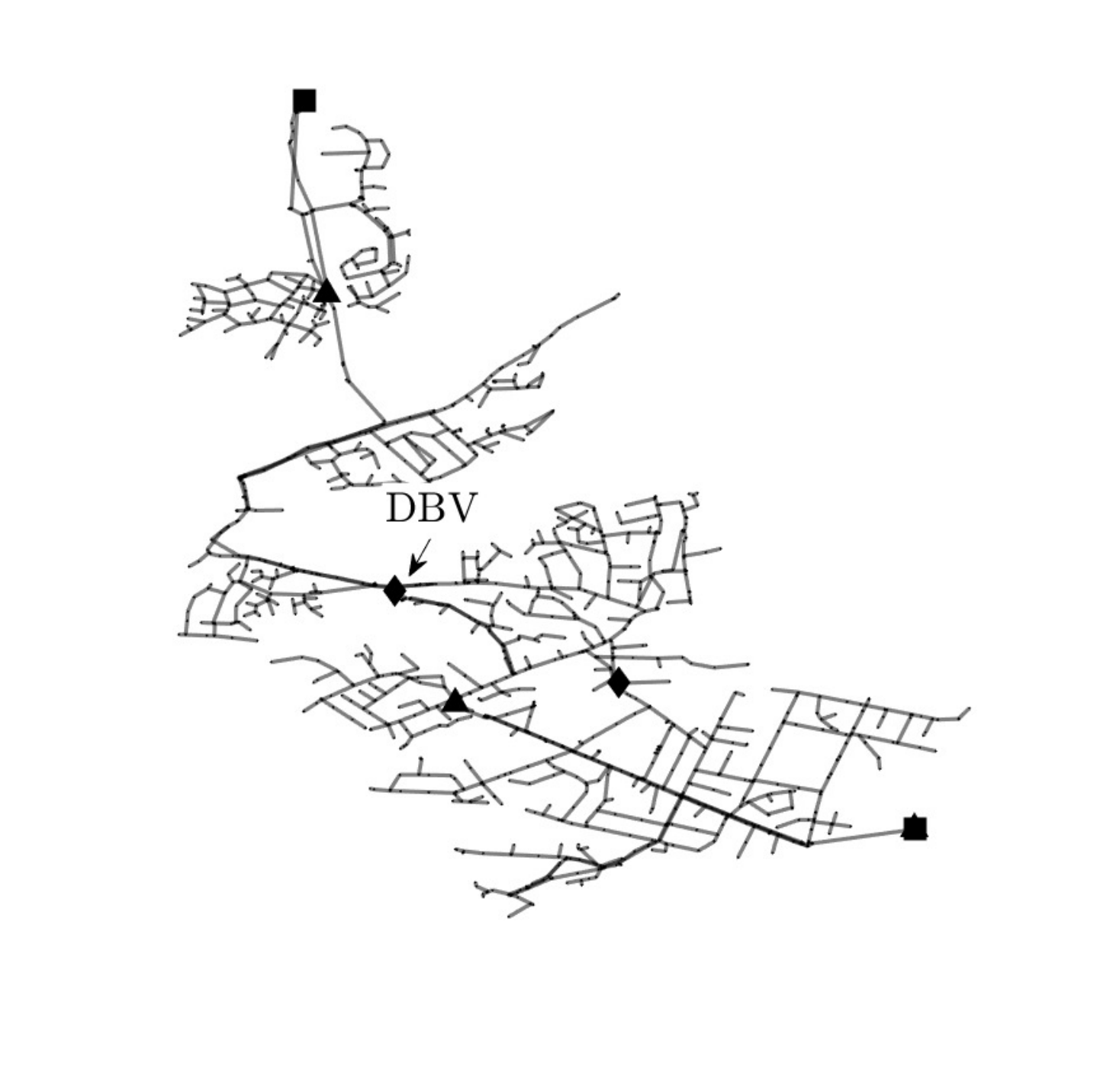}}

    \caption{Case study network layouts}
    \label{fig:netLayouts}
\end{figure}

We evaluated the performance of the solution methods using three case study networks: \texttt{Pescara}, \texttt{Modena}, and \texttt{BWFLnet}. Problem data for each network, including the size of the resulting MINLP optimization problems, are summarized in \Cref{table:prob_data}. Moreover, network layouts are shown in \Cref{fig:netLayouts}. \texttt{Pescara} and \texttt{Modena} are skeletonized versions of WDNs in medium-sized Italian cities, published by \citet{BRAGALLI2012} for benchmarking purposes. Since these are theoretical networks, unidirectional PRVs were placed at the outlet of each reservoir to represent existing valve configurations. \texttt{BWFLnet} is a large-scale operational network representing an area of the City of Bristol's WDN, located in southwest England. It currently operates with dynamically controlled PRV and bidirectional DBV valves, which modulate flow and pressure between adjacent zones to minimize AZP. \texttt{BWFLnet}'s current control valve configuration was applied to model existing conditions - for more details on \texttt{BWFLnet} see \citet{WRIGHT2014} and \citet{WALDRON2020}. Moreover, existing kept-shut boundary valves were opened in order to test the full capabilities of the proposed solution algorithm. Observe that the evaluated case study networks vary by size, network connectivity, and demand scenarios, enabling the proposed CMS algorithm to be tested on a range of network conditions. \texttt{Pescara} and \texttt{Modena} are relatively small-scale looped networks with theoretical (and quite large) demands. \texttt{BWFLnet}, on the other hand, represents a larger-scale operational network having a branched structure and demands representing water consumption in a typical urban area in the UK. For reference, we note that the size of the resulting MINLP optimization problems for \texttt{BWFLnet} are at the upper end of the range reported in the MINLPLib database \citep{VIGERSKE2022}; MINLPLib is a widely used library of benchmarking MINLP problem instances.

\begin{table}[h]
    \centering
    \setlength{\tabcolsep}{16pt}
    \caption{Case study network problem data}
    \label{table:prob_data}
        \begin{tabular}{lrrr}
            \toprule
            Problem data & \texttt{Pescara} & \texttt{Modena} & \texttt{BWFLnet}\\
            \midrule
            $n_n$ & 67 & 268 & 2,745\\
            $n_0$ & 3 & 4 & 2\\
            $n_p$ & 98 & 317 & 2,816\\
            Existing PRVs ($n_{\text{PRV}}$) & 4 & 4 & 3\\
            Existing DBVs ($n_{\text{DBV}}$) & 0 & 0 & 2\\
            No. free cont. variables & 660 & 2,340 & 22,224\\
            No. total cont. variables & 1,712 & 5,948 & 55,752\\
            No. binary variables & 949 & 3,121 & 28,089\\
            No. nonconvex terms & 784 & 2,536 & 22,528\\
            \bottomrule
        \end{tabular}
    \end{table} 

The known input parameters for Problem \eqref{eq:MINLP_problem} were selected as follows. For the SCC objective function $f_{\widetilde{\text{SCC}}}$ \eqref{eq:SCC_objective_sigmoid}, we set the sigmoid function curvature parameter $\rho$ to 50, unless stated otherwise. On the basis of the reviewed literature in \Cref{sec:intro}, the self-cleaning velocity threshold $u_j^{\min}$ was set to 0.2 \si{\meter/\second} for all $ j \in \{1,\ldots,n_p\}$. Since we are concerned with the redistribution of peak flows, we define four discrete hydraulic time steps $n_t$ to represent peak demand scenarios within a typical diurnal demand profile. These were split between the morning and evening high demand periods and were unique to each case study. In relation to hydraulic feasibility constraints \eqref{eq:hyd_bounds_b}, we applied a minimum regulatory pressure head of 15 m (UK regulations) at all nonzero demand nodes, for all case study networks. Moreover, the upper bound for flushing valve demands $\alpha^U$ shown in constraint \eqref{eq:hyd_bounds_d} was set to 25 L/s. This value was selected with reference to the recommended fire flow capacity at fire hydrants in the UK \citep{FIRE2007}. For the design-for-control problem, we evaluated the placement and operation of $n_v = 1,\dots,3$ bidirectional DBVs and $n_f = 0,\dots,3$ AFVs. Hence, for each network, a total of 12 valve configurations were considered. The number of randomization trials $N$ in \Cref{alg:randomization heuristic} was set to the minimum of valve configuration combinations generated from the continuous relaxation in Subproblem \eqref{eq:LP_problem} or a maximum trial count; maximum trial counts were set to $N = 50$ for \texttt{Pescara} and \texttt{Modena} and $N = 100$ for \texttt{BWFLnet}. Finally, the number of initial conditions tested in the multi-start scheme (\Cref{fig:solution_multistart}) was set to $M = 5$.

\subsection{Existing valve configuration results}
\label{sec:results_control}
\begin{figure}[h!]
    \centering
    \captionsetup{justification=centering}
    \subfloat[\label{fig:controlProb_CDFvel_a}\text{\texttt{Pescara}}]{
        \includegraphics[width=0.6\textwidth]{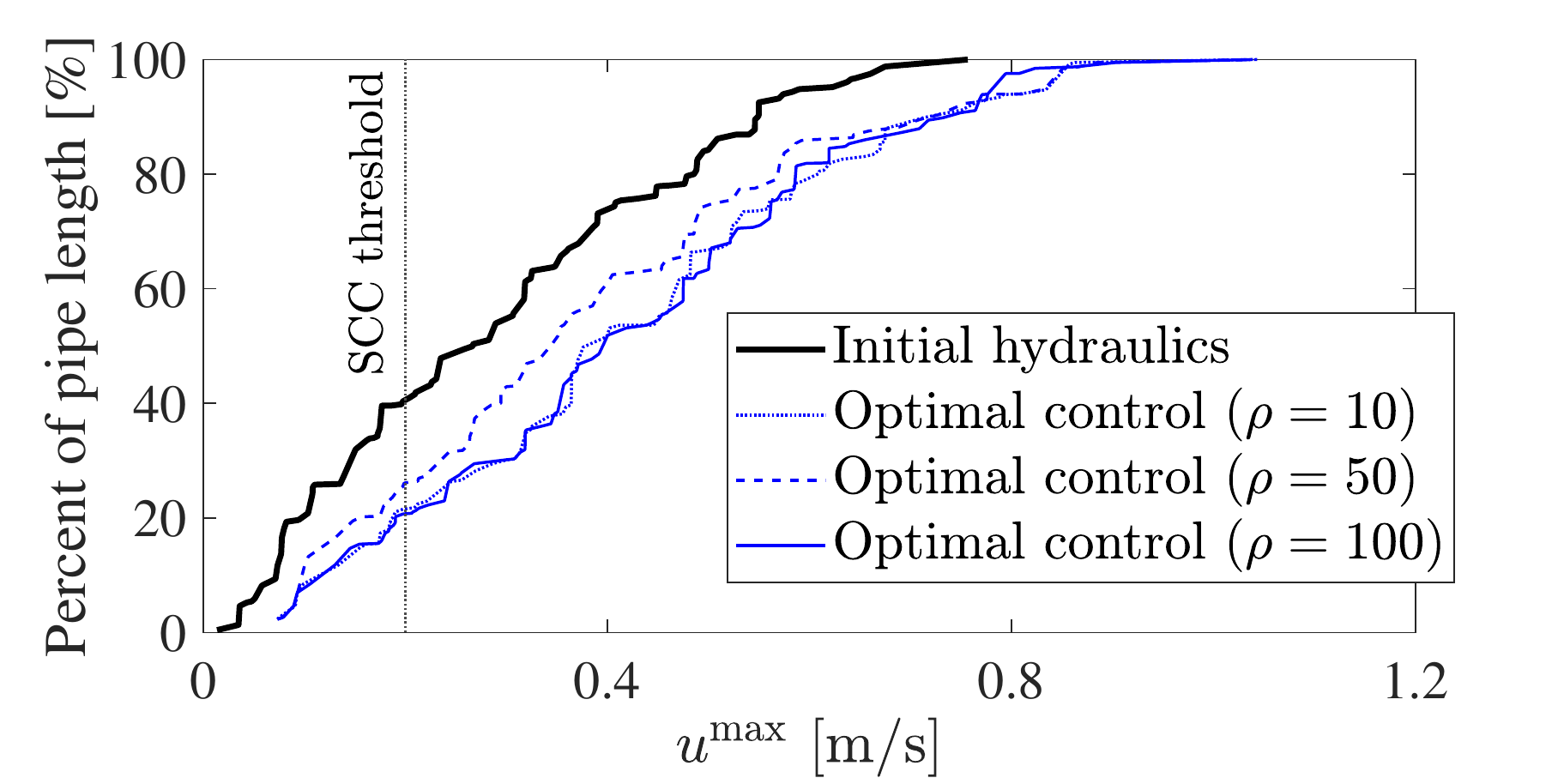}}
    \vspace{0.025cm}   
    \subfloat[\label{fig:controlProb_CDFvel_b}\text{\texttt{Modena}}]{
        \includegraphics[width=0.6\textwidth]{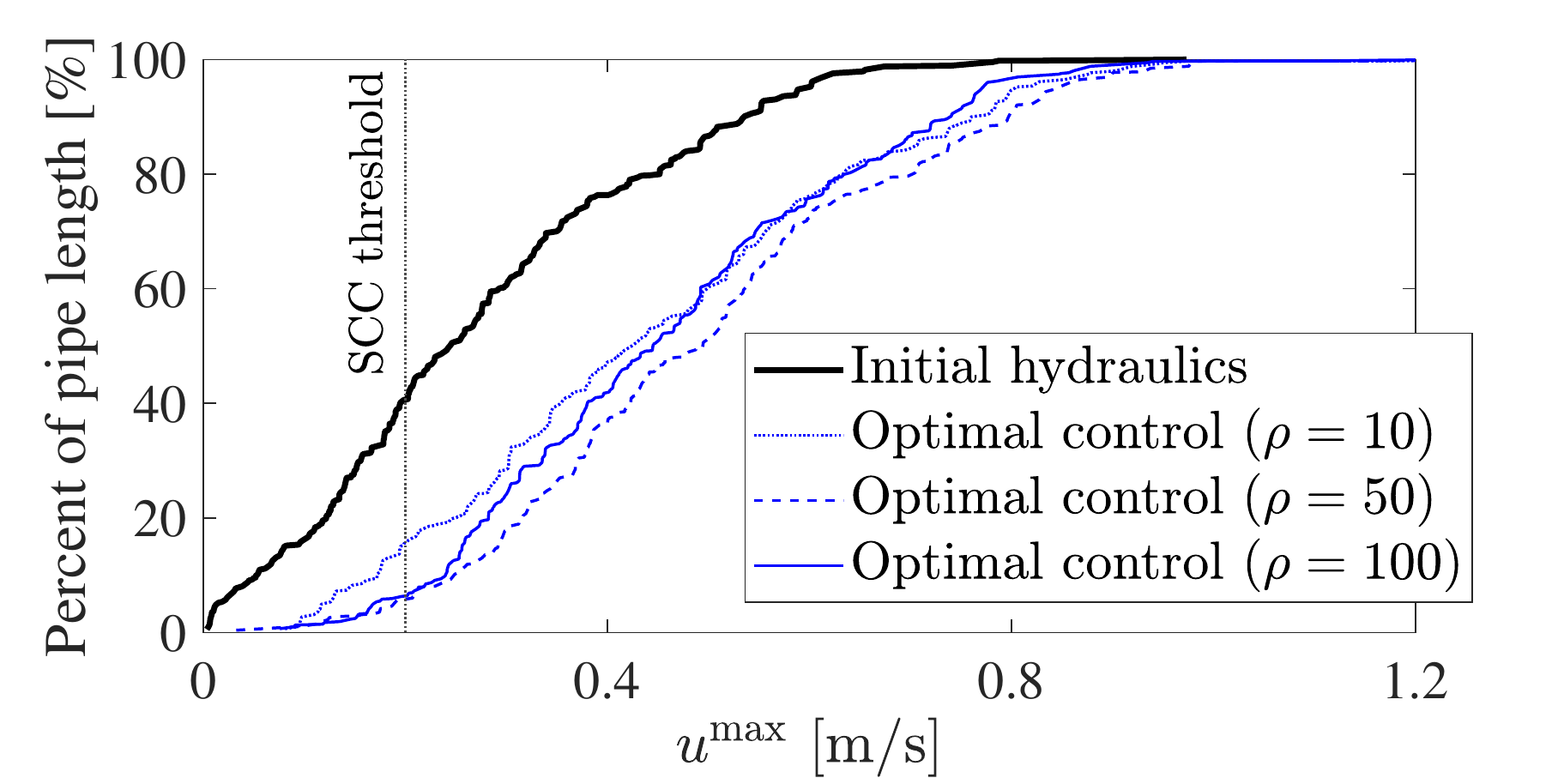}}
    \vspace{0.025cm}    
    \subfloat[\label{fig:controlProb_CDFvel_c}\text{\texttt{BWFLnet}}]{
        \includegraphics[width=0.6\textwidth]{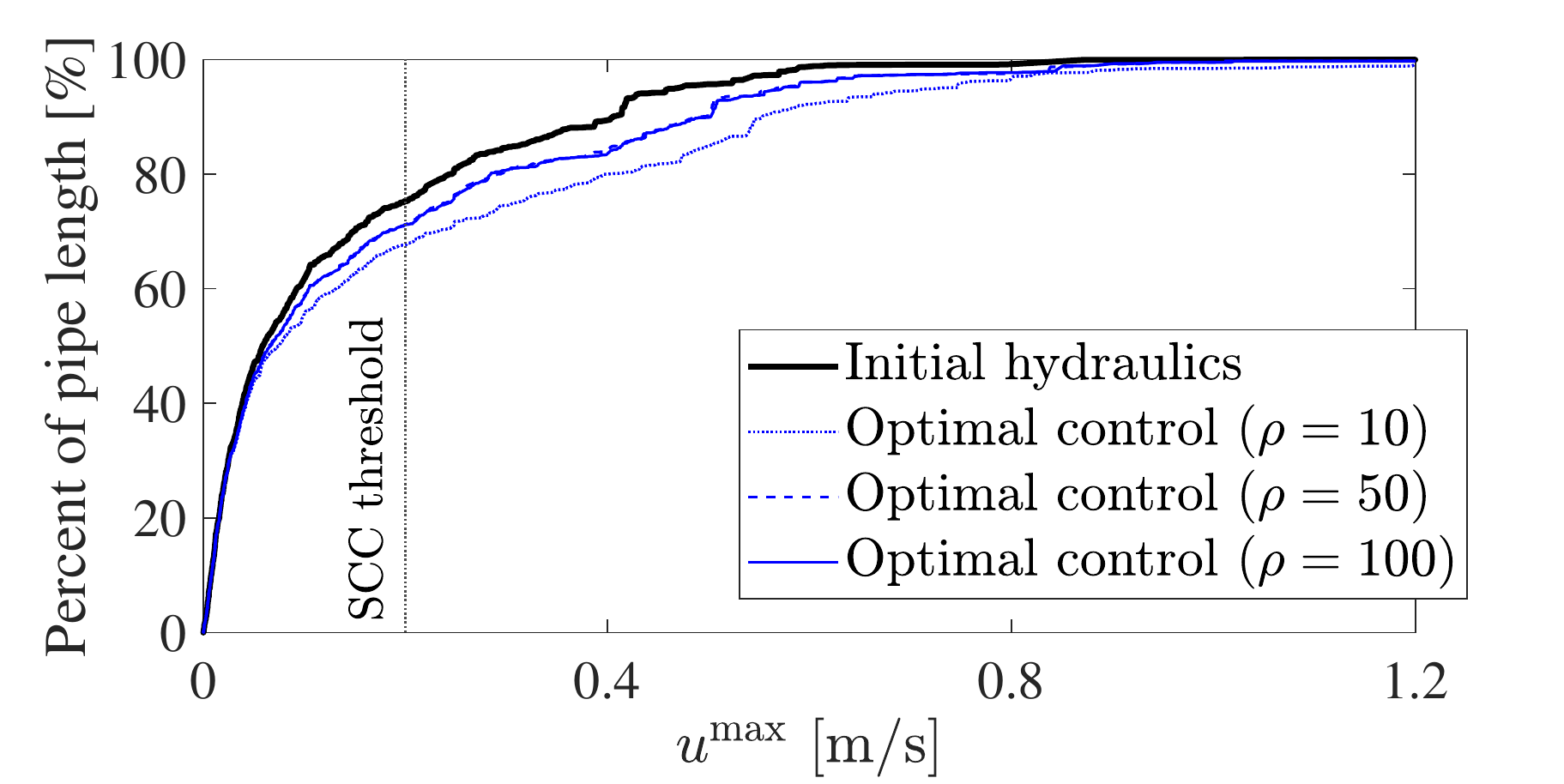}}
    
    \caption{Cumulative distribution plots of maximum pipe flow velocities for control and no control scenarios (with $\rho$ sensitivity)}
    \label{fig:controlProb_CDFvel}
\end{figure}

We first investigated improvements to SCC performance achieved by optimizing the operational settings of existing pressure control valve configurations. This was considered a logical initial step as the implementation SCC-specific controls requires no additional capital investment and can thus be readily applied. As discussed previously, the theoretical case study networks, \texttt{Pescara} and \texttt{Modena}, were assumed to have PRVs assigned at each link leaving a source node (e.g. reservoir). On the other hand, \texttt{BWFLnet} had its existing operational valve configuration assigned. A subproblem of Problem \eqref{eq:MINLP_problem} was then formulated where binary variables $z$ were fixed for the known control valve locations and directional variables $v^+$ and $v^-$ were set free for networks with existing DBV locations (e.g. \texttt{BWFLnet}). Valve control settings and DBV directions were then optimized over the four peak demand periods using the multi-start strategy described in \Cref{sec:multistart_solver}. The results are reported in \Cref{fig:controlProb_CDFvel} through cumulative distribution plots of maximum flow velocities, weighted by pipe length. These plots offer a comparison between network SCC conditions with and without optimized controls. Considering $f_{\widetilde{\text{SCC}}}$ is an approximation of the actual self-cleaning problem, we also include a sensitivity analysis on the sigmoid curvature parameter $\rho$ to test changes in relative performance. Moreover, \Cref{table:controlProb_results} summarizes the control problem results obtained with both the multi-start solver and the GA implementation described in \Cref{sec:results_setup}. We set a maximum time limit of one hour for the off-the-shelf GA solver, which was assumed to reflect an upper bound for the practical application of control problems.

\begin{table}[h]
    \centering
    \setlength{\tabcolsep}{16pt}
    \captionsetup{width=0.95\linewidth, skip=6pt}
    \caption{Control problem results for existing valve configurations}
    \label{table:controlProb_results}
    \small
    \begin{tabular}{lrrr}
    \toprule
    & \multicolumn{3}{c}{Case study networks}\\ \cmidrule{2-4}
     & \texttt{Pescara} &  \texttt{Modena} &  \texttt{BWFLnet} \\
    \midrule
    \rowgroup{SCC objective (\%)} \\
    NLP solver & 63.7 & 69.1 & 24.1 \\
    (CPU time, s) & (1.61) & (4.64) & (200) \\
    GA solver &  63.8 & 65.8 & 23.7 \\
    (CPU time, s) & (3,600) & (3,600) & (3,600) \\
    \midrule
    \rowgroup{AZP objective (m)} \\
    Optimal AZP control & 27.4& 18.0 & 36.4 \\
    Optimal SCC control & 32.6 & 23.3 & 52.6 \\
    \bottomrule
    \end{tabular}
\end{table}

First, \Cref{fig:controlProb_CDFvel} suggests that pressure modulation at existing control valves can significantly improve SCC conditions in the tested looped networks. This is highlighted by the numerical results for \texttt{Modena}, where the length of pipe experiencing self-cleaning velocities increased from approximately 60\% to 95\% for peak demand conditions. In contrast, results for \texttt{BWFLnet}, which has limited connectivity due to its branched structure, showed only modest improvements. These results were expected as the reduced inter-connectivity in \texttt{BWFLnet} leads to pipe flows being primarily a function of downstream demands. For this reason, \Cref{sec:results_dfc} presents the results for a design-for-control problem where we consider the installation of additional DBVs and flushing valves, which aim to improve self-cleaning conditions in branched networks. Furthermore, we report modest changes in the flow velocity profiles for the $\rho$ sensitivity analysis (see \Cref{fig:controlProb_CDFvel}). Since no pattern can be discerned between $\rho$ values, we conclude $\rho = 50$ to be suitable for the design-for-control numerical experiments discussed in the subsequent section.

As shown in \Cref{table:controlProb_results}, the multi-start solver computes solutions that are similar to those obtained by the GA, but at a fraction of the computational effort. While the GA consistently reaches the set time-limit of one hour, the longest computational time experienced by the multi-start solver is equal to 200 seconds. This is important as we are considering an optimal control problem, which may need to be solved in near real-time. The much longer CPU time required by the GA is likely due to the fact that it facilitates a search for the globally optimal solution. Thus, for highly nonlinear optimization problems, like that formulated for the current SCC problem, the GA retains a large solution space so as to avoid exclusion of the global optimum. This characteristic has been recognized as a challenge when using GAs for WDN control problems, especially when dealing with problems having a large number and range of free continuous variables \citep{MAIER2014,ULUSOY2022b}. Furthermore, we note that the GA implementation used in this work applies default parameters (see \Cref{sec:results_setup}) and loads an external hydraulic solver (EPANET2.2) at each fitness evaluation. This off-the-shelf implementation is also likely to contribute to the relatively large CPU times. Extending the time limit to nearly 48 hours, the GA converges to a solution with an SCC value of approximately 24.8\% for \texttt{BWFLnet}. This minimal performance gain suggests that the one-hour time limit is sufficient for the SCC control problem comparison.

In addition, \Cref{table:controlProb_results} provides a comparison of average zone pressure (AZP) performance between control settings optimized for AZP and SCC, respectively. Note that the AZP objective function is formulated as a weighted sum of nodal pressures averaged over $n_t$ time steps (\citeauthor{WRIGHT2015}, \citeyear{WRIGHT2015}; Equation 4). In accordance with findings in \citet{ABRAHAM2016}, we found a clear trade-off between the optimization of control settings for the SCC and AZP objectives. For example, a 16.2 m difference in AZP was recorded for \texttt{BWFLnet}. While this trade-off is not discussed further in the current study, we acknowledge the importance of having SCC and AZP objectives (among others) coexist in an overall operational network control strategy. Furthermore, the transition between such operational objectives may provoke large variations in hydraulic conditions (e.g. pressure transients) and should thus be managed accordingly. Future work is needed to investigate these considerations for the implementation of dynamically adaptive control strategies in operational WDNs.

\subsection{Design-for-control results}
\label{sec:results_dfc}
\begin{figure}[h!]
    \vspace{-0.2cm}
    \centering
    \captionsetup{justification=centering}
    \subfloat[\label{fig:DfC_Pescara_a}\text{\texttt{Pescara} (SCC objective)}]{
        \includegraphics[width=0.45\textwidth]{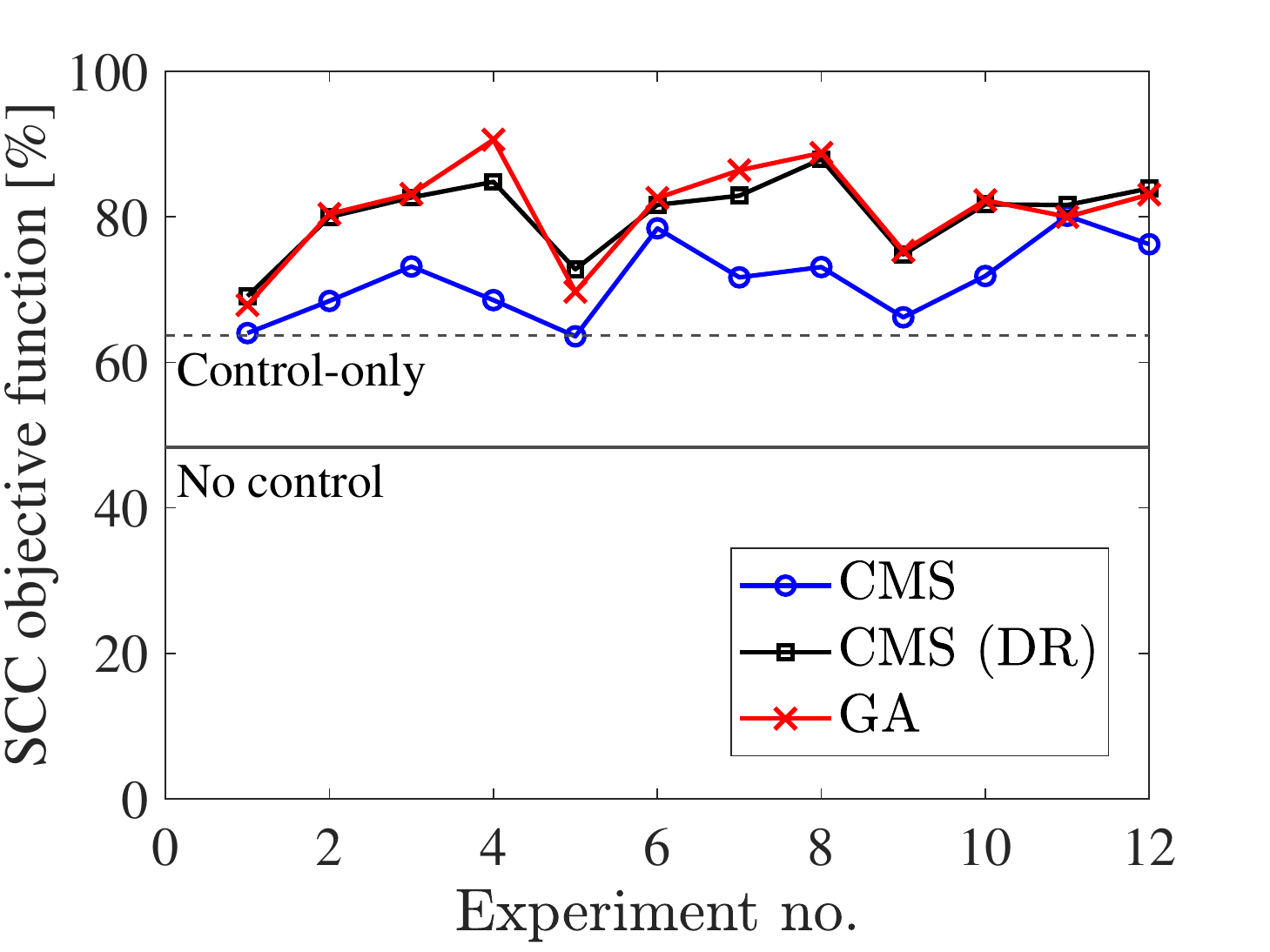}}
    \hspace{0.6cm}
    \subfloat[\label{fig:DfC_Pescara_b}\text{\texttt{Pescara} (CPU time)}]{
        \includegraphics[width=0.45\textwidth]{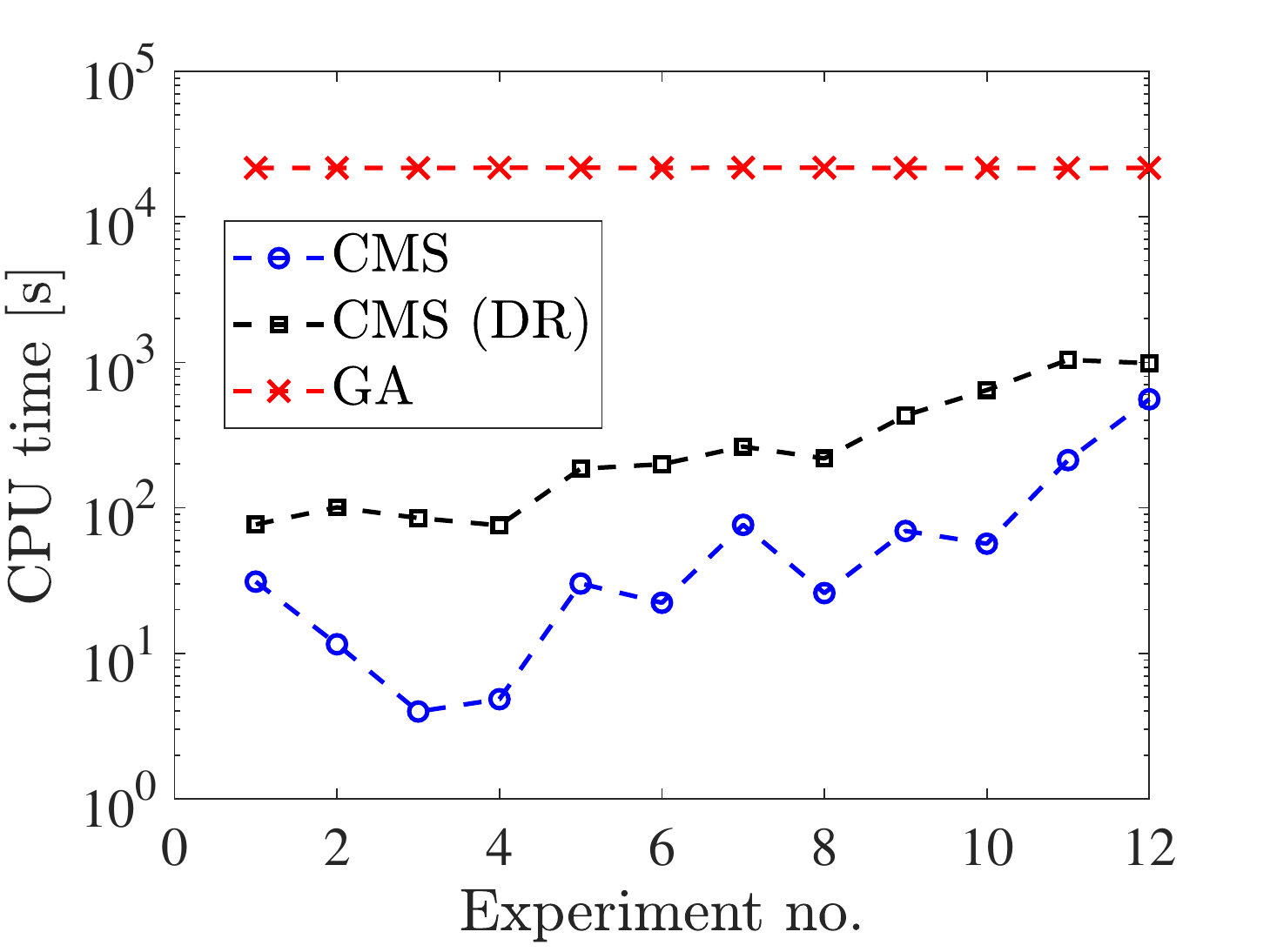}}
    \vspace{0.1cm}    
    \subfloat[\label{fig:DfC_Modena_a}\text{\texttt{Modena} (SCC objective)}]{
        \includegraphics[width=0.45\textwidth]{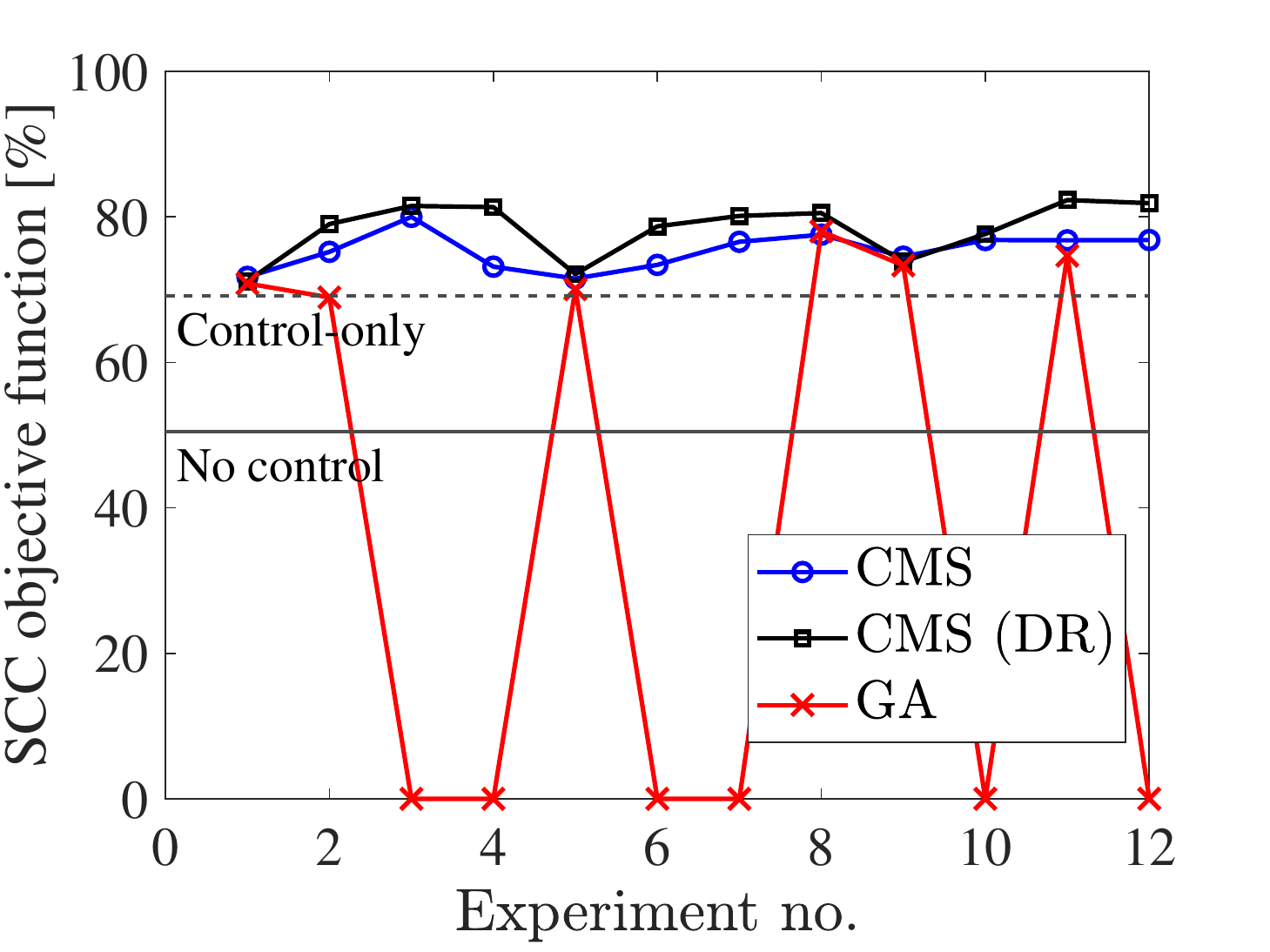}}
    \hspace{0.6cm}
    \subfloat[\label{fig:DfC_Modena_b}\text{\texttt{Modena} (CPU time)}]{
        \includegraphics[width=0.45\textwidth]{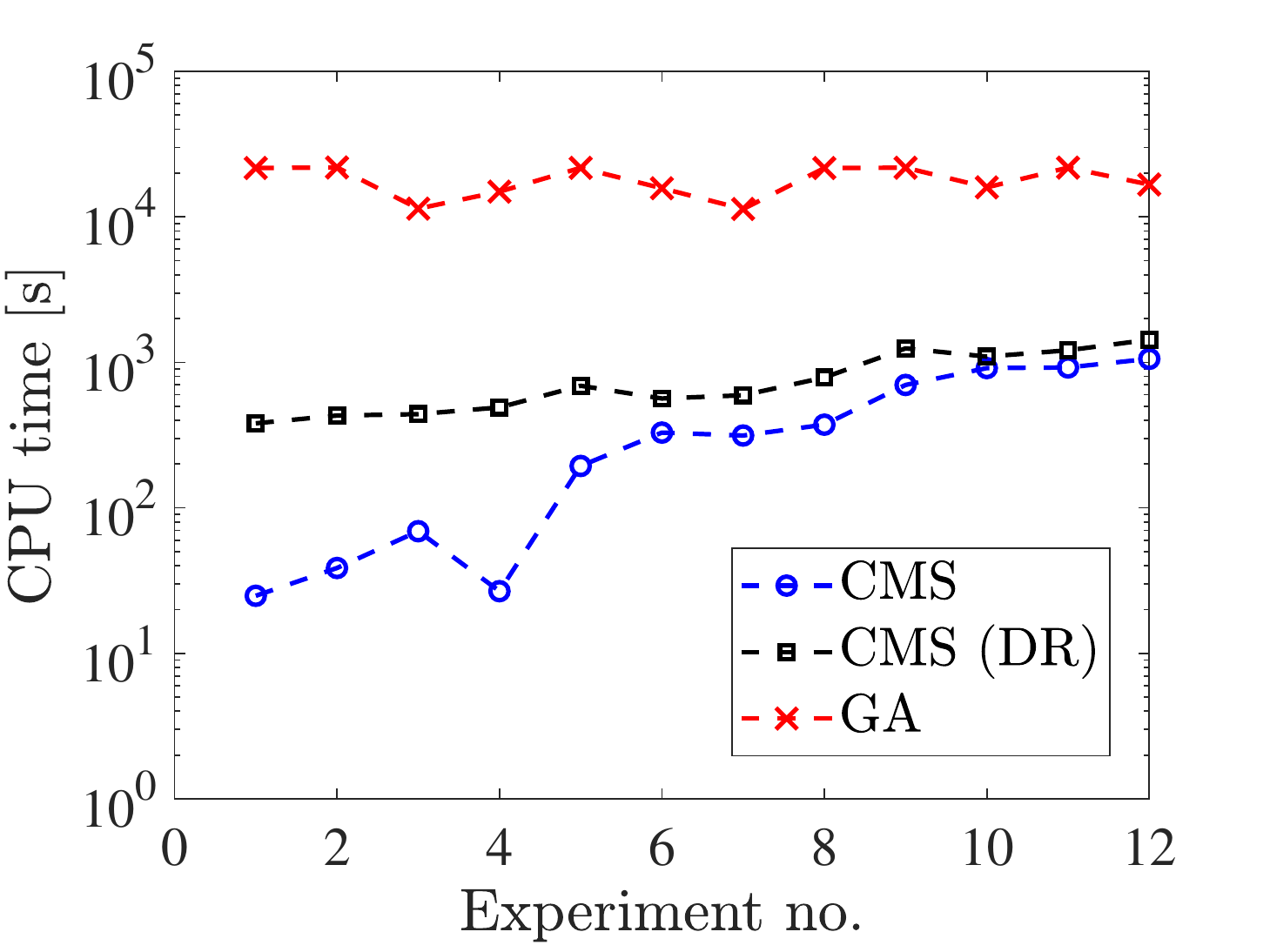}}
    \vspace{0.1cm}    
    \subfloat[\label{fig:DfC_BWFLnet_a}\text{\texttt{BWFLnet} (SCC objective)}]{
        \includegraphics[width=0.45\textwidth]{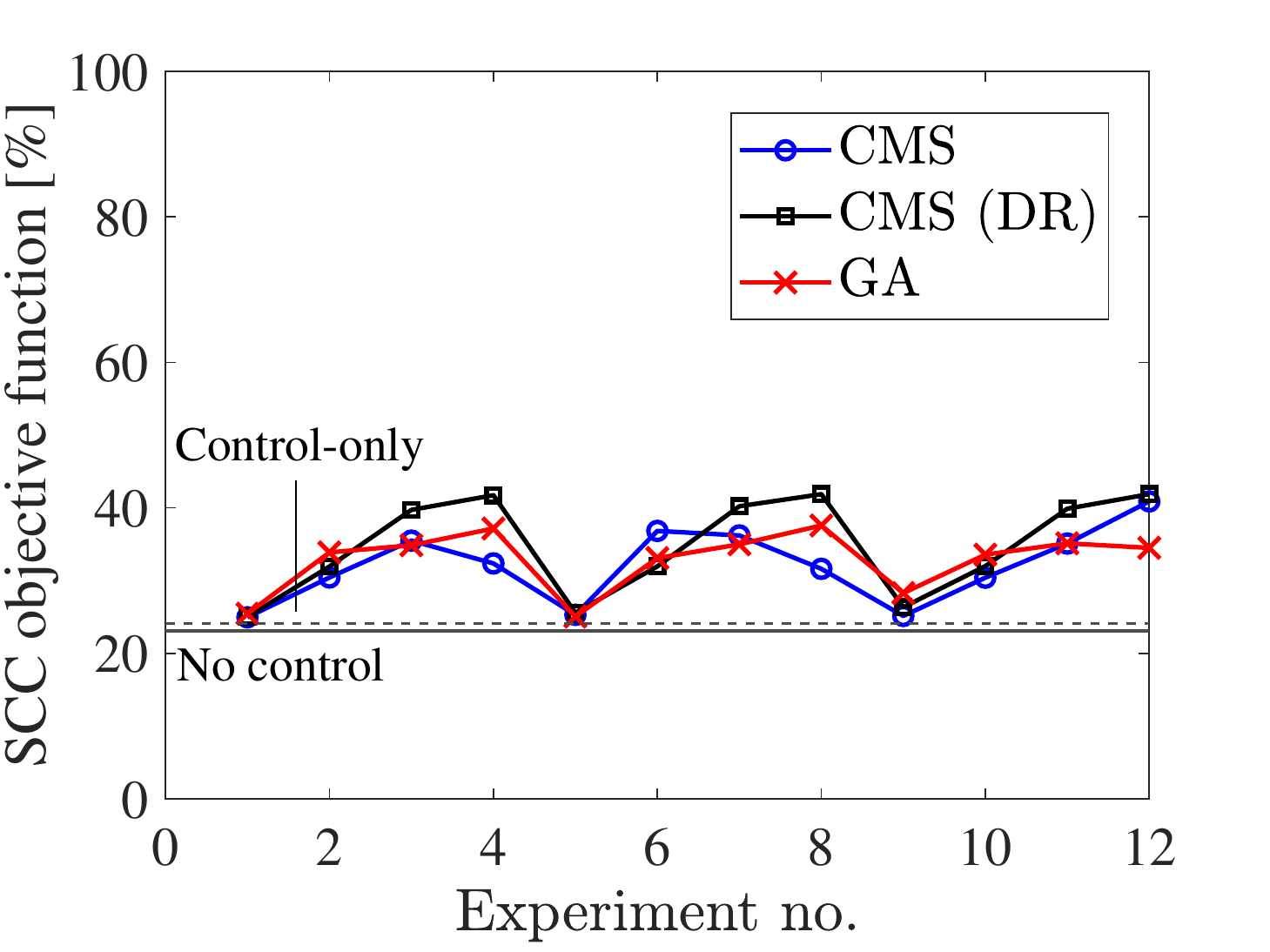}}
    \hspace{0.6cm}
    \subfloat[\label{fig:DfC_BWFLnet_b}\text{\texttt{BWFLnet} (CPU time)}]{
        \includegraphics[width=0.45\textwidth]{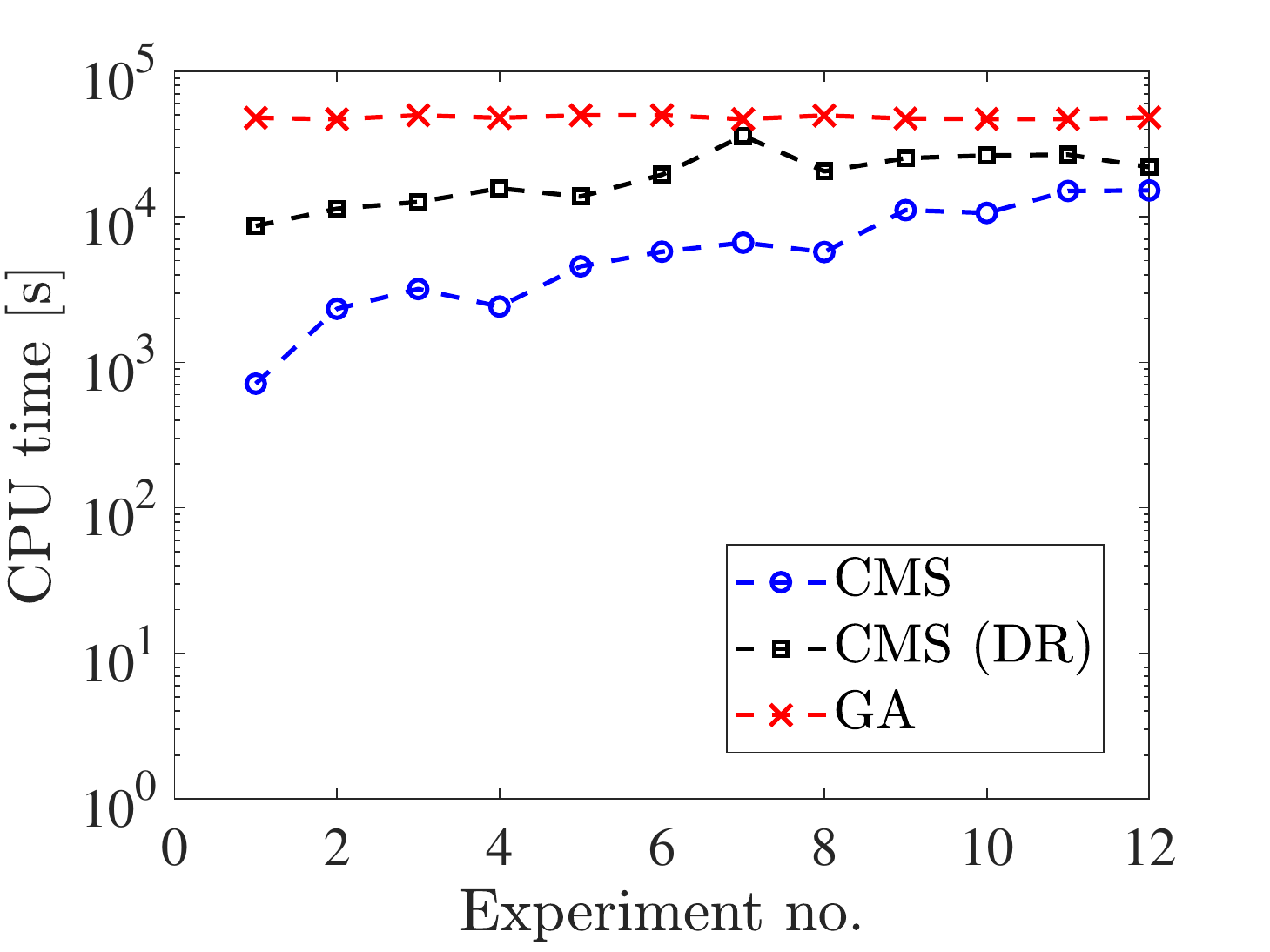}}

    \caption{Design-for-control solutions to Problem \eqref{eq:MINLP_problem} for experiments corresponding to $n_v = 1,\dots,3$ and $n_f = 0,\dots,3$}
    \label{fig:DfC_results}
\end{figure}

This section considers the SCC design-for-control problem for optimal placement and control of new bidirectional DBVs and AFVs. We test a total of 12 valve configurations for each case study network, where Experiment 1 corresponds to $n_v = 1$ and $n_f = 0$ and Experiment 12 corresponds to $n_v = 3$ and $n_f = 3$. As stated in \Cref{sec:results_setup}, each experiment uses hydraulic conditions from four peak demand conditions across a typical diurnal demand profile. Moreover, with the exception of existing DBVs in \texttt{BWFLnet}, the design-for-control numerical experiments include the existing PRV configurations shown in \Cref{fig:netLayouts}, which were used for the optimal control subproblem in \Cref{sec:results_control}. We investigate the solution of Problem \eqref{eq:MINLP_problem} using three methods: (i) the proposed convex multi-start heuristic (CMS); (ii) CMS with domain reduction (as described in \Cref{sec:convex_relax}), referred to as CMS (DR); and (iii) an off-the-shelf GA implementation coupled with the EPANET2.2 hydraulic solver. Note that we set the following time limits for the GA on the basis of that recorded using CMS: 6 hours for \texttt{Pescara} and \texttt{Modena}; and 12 hours for \texttt{BWFLnet}. The SCC objective function results and corresponding CPU times from each solution method are shown in \Cref{fig:DfC_results} for the considered set of numerical experiments.

\begin{figure}[p]
    \centering
    \captionsetup{justification=centering}
    \subfloat[\label{fig:DfC_vel_Pescara_a}\text{\texttt{Pescara} (no control)}]{
        \includegraphics[width=0.43\textwidth]{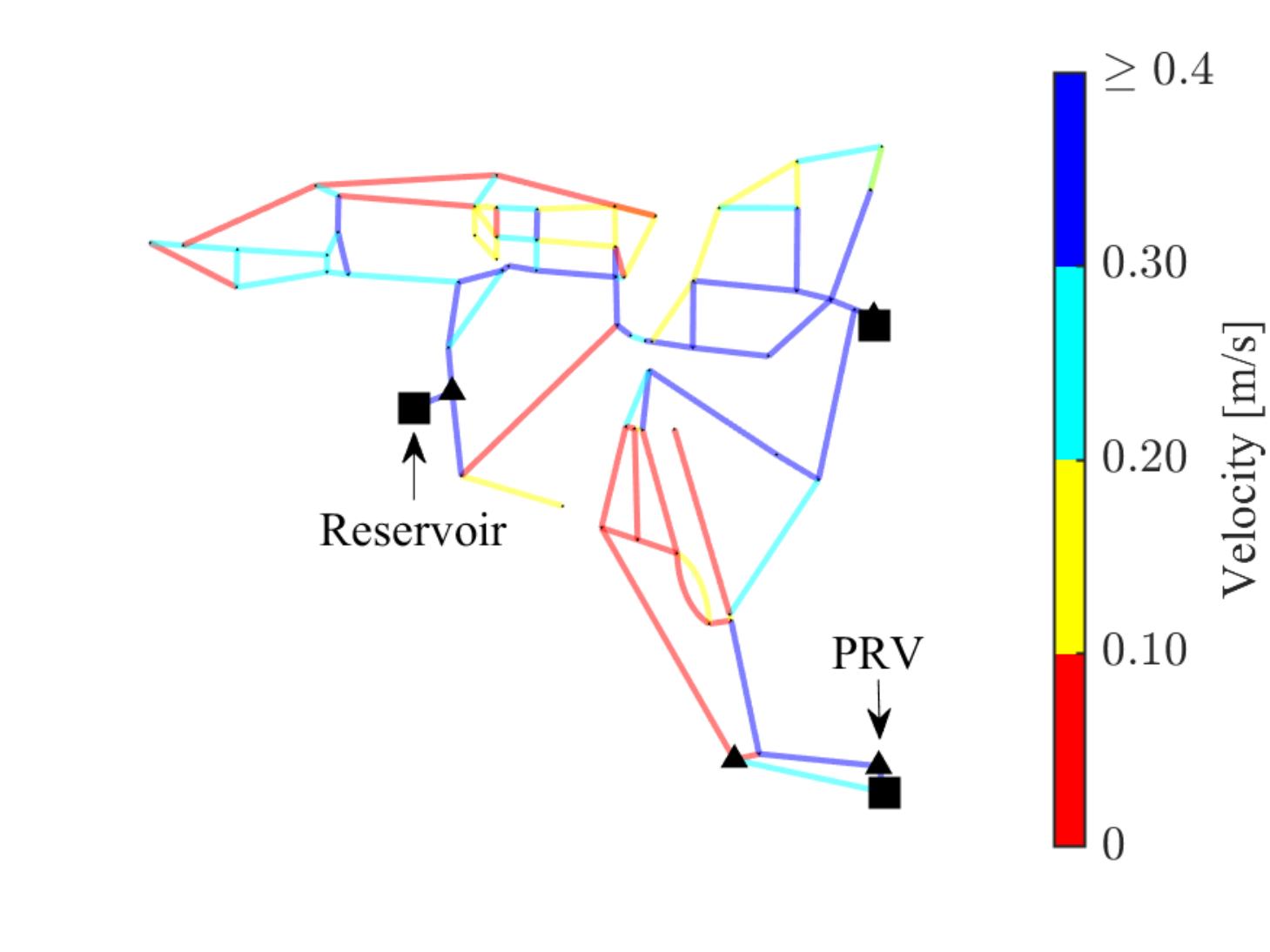}}
    \hspace{0.25cm}
    \subfloat[\label{fig:DfC_vel_Pescara_b}\text{\texttt{Pescara} (optimal design-for-control)}]{
        \includegraphics[width=0.43\textwidth]{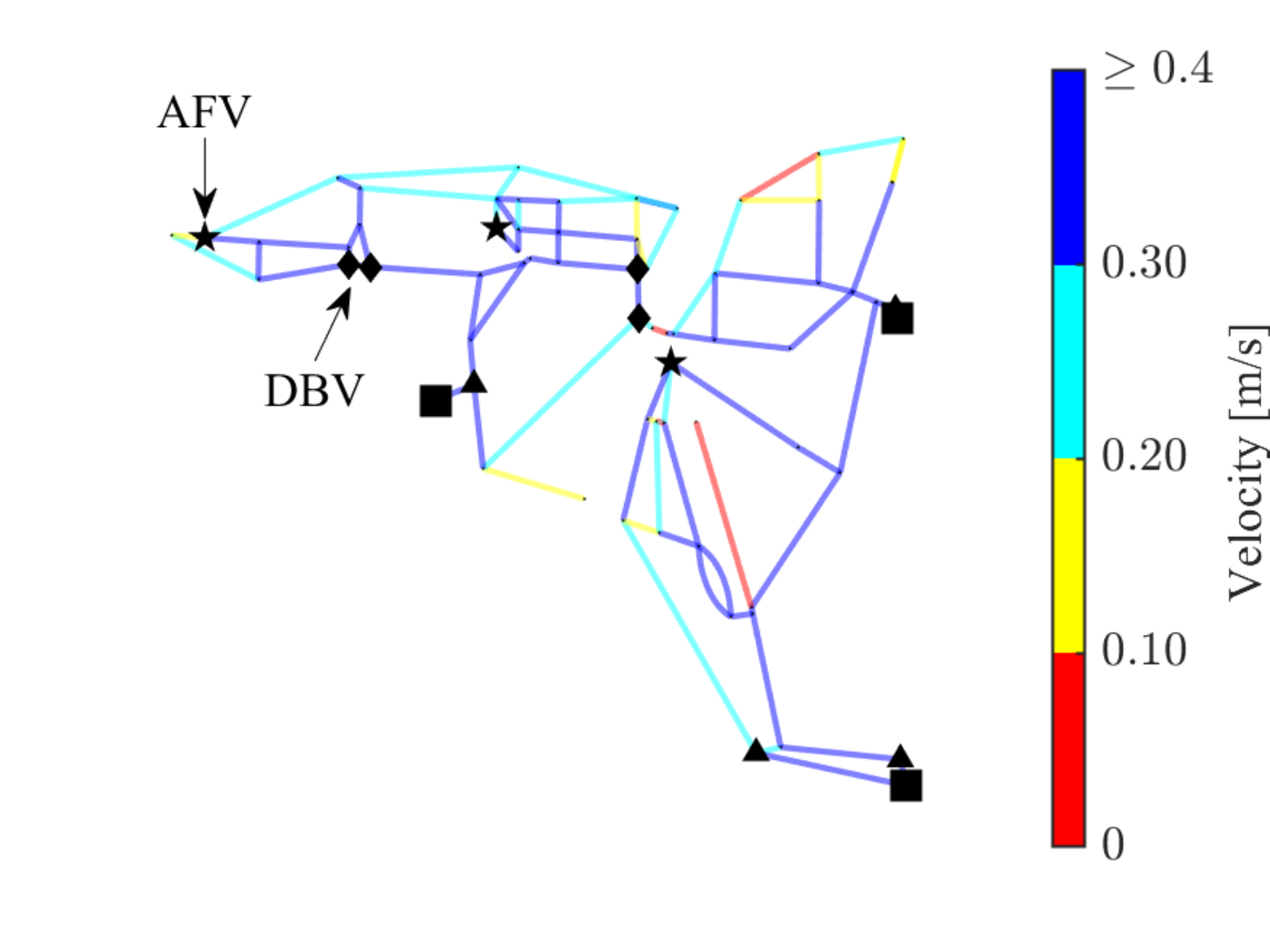}}
    \vspace{0.15cm}    
    \subfloat[\label{fig:DfC_vel_Modena_a}\text{\texttt{Modena} (no control)}]{
        \includegraphics[width=0.43\textwidth]{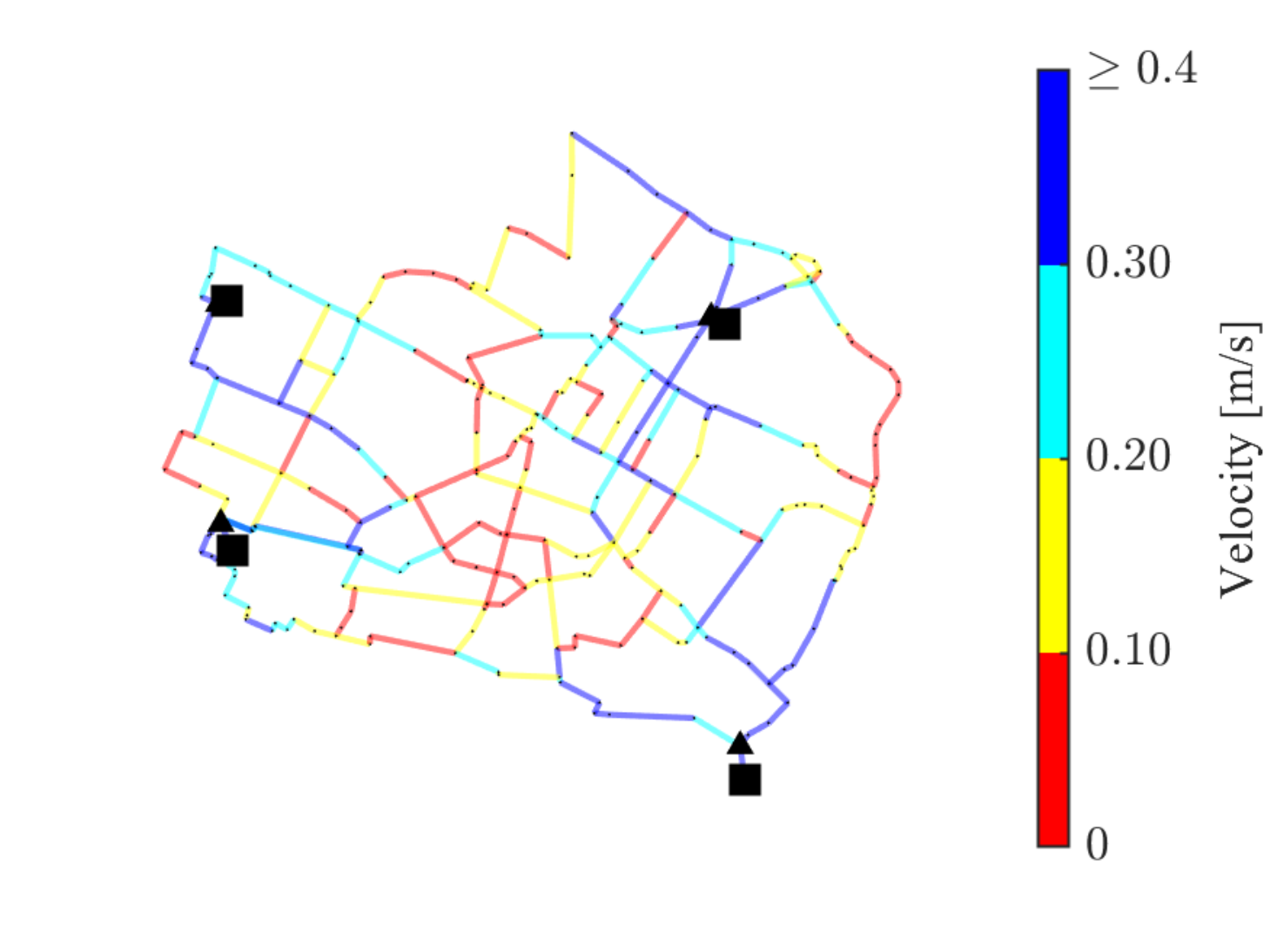}}
    \hspace{0.25cm}
    \subfloat[\label{fig:DfC_vel_Modena_b}\text{\texttt{Modena} (optimal design-for-control)}]{
        \includegraphics[width=0.43\textwidth]{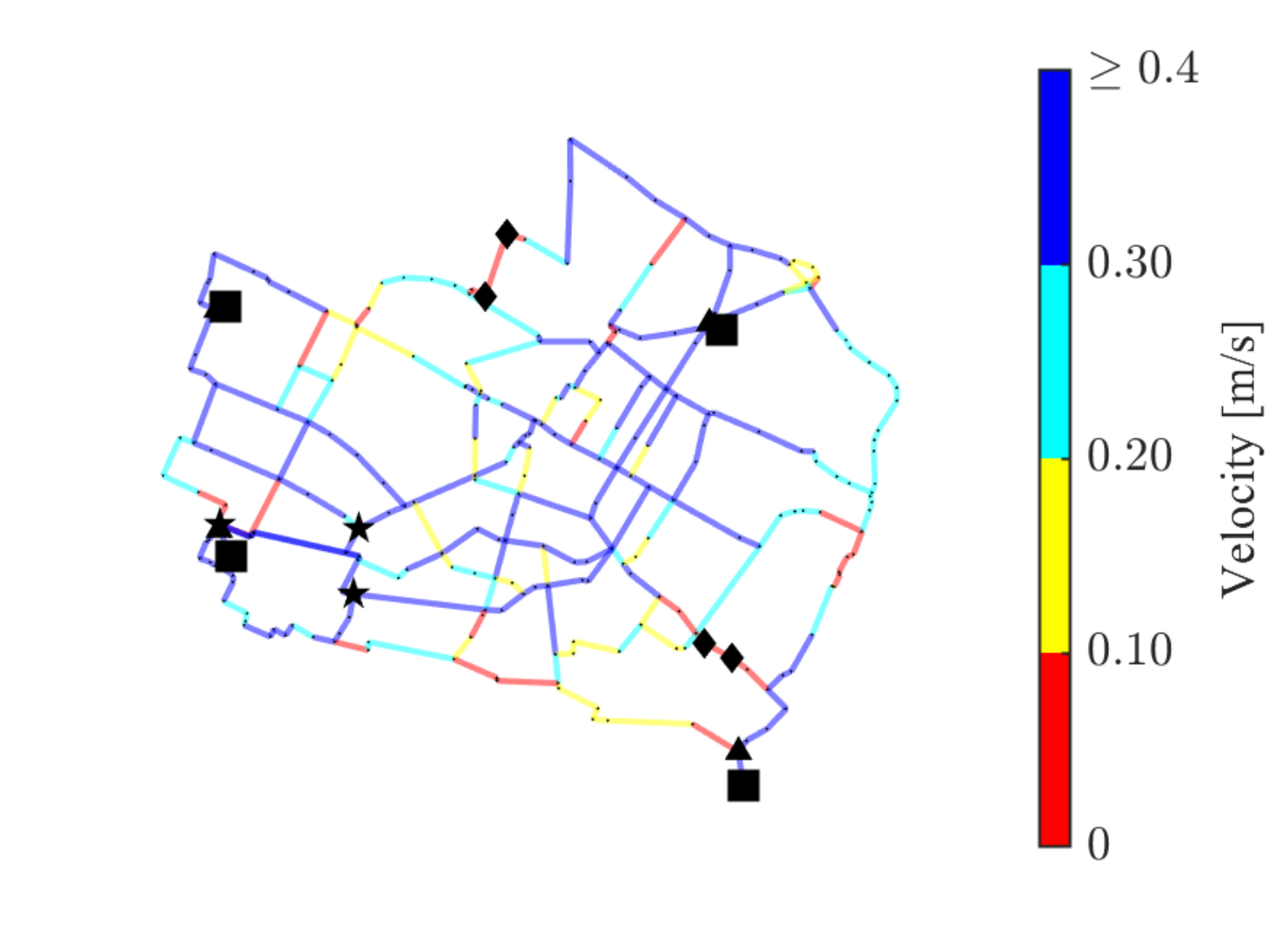}}
    \vspace{0.15cm}    
    \subfloat[\label{fig:DfC_vel_BWFLnet_a}\text{\texttt{BWLFnet} (no control)}]{
        \includegraphics[width=0.46\textwidth]{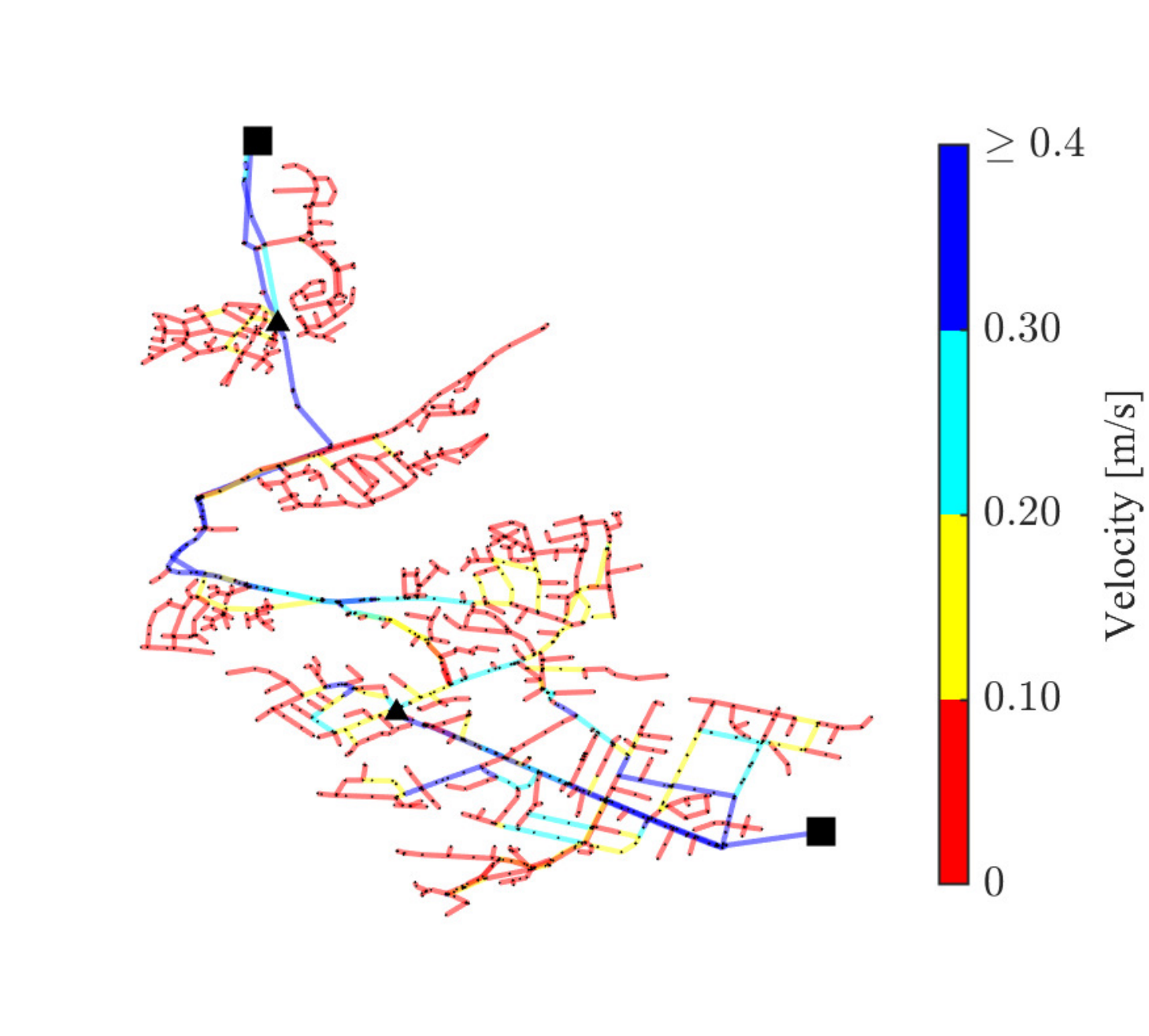}}
    \hspace{0.15cm}
    \subfloat[\label{fig:DfC_vel_BWFLnet_b}\text{\texttt{BWFLnet} (optimal design-for-control)}]{ 
        \includegraphics[width=0.46\textwidth]{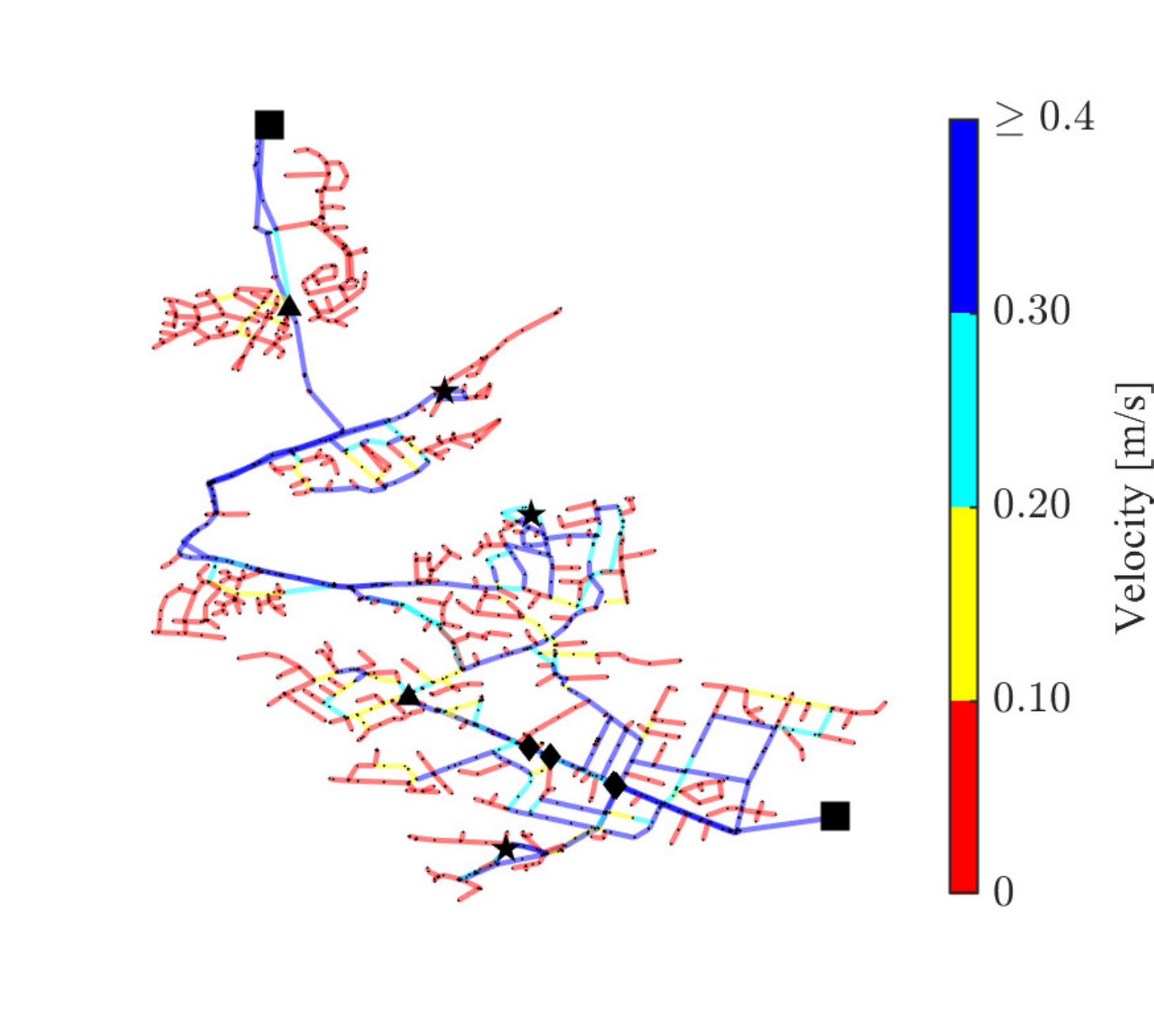}}
    \caption{Improvements to network self-cleaning velocities from optimal design-for-control solutions ($n_v = 2$ and $n_f = 3$)}
    \label{fig:DfC_vel}
\end{figure}

Figures~\ref{fig:DfC_Pescara_a}, \ref{fig:DfC_Modena_a}, and \ref{fig:DfC_BWFLnet_a} show that CMS finds feasible solutions for all numerical experiments. To highlight the improvements in SCC, we include a solid black line to represent initial hydraulic conditions and a dashed black line to report the SCC values obtained by optimally controlling existing valves, as done in \Cref{sec:results_control}. While SCC improvements are experienced for all case studies, we observe significant relative improvements compared with the control-only solution for \texttt{BWFLnet}. As previously discussed in \Cref{sec:results_control}, the SCC improvements for \texttt{BWFLnet} can be attributed to the deployment of AFVs, which facilitate controlled flushing demands at the optimally selected locations. Moreover, the additional DBVs act to further branch the network, thereby enhancing self-cleaning velocities. The increase in SCC is investigated further through network velocity plots shown in \Cref{fig:DfC_vel}. Here, we compare maximum pipe velocities resulting from initial hydraulics (i.e. no control) and, as an example, the design-for-control solution corresponding to $n_v = 2$ and $n_f = 3$. Figures~\ref{fig:DfC_vel_Pescara_b}, \ref{fig:DfC_vel_Modena_b}, and \ref{fig:DfC_vel_BWFLnet_b} clearly indicate a net increase in the overall extent of the network experiencing self-cleaning velocities from optimal valve controls. We observe optimal AFV placement to be near the network periphery in \texttt{Pescara} and \texttt{BWFLnet}. These results are intuitive as flushing demands facilitate self-cleaning velocities in areas that otherwise convey relatively low flow. In contrast, AFVs are concentrated near a source node in \texttt{Modena}. This is a result of existing PRVs modulating inlet pressure to shift the hydraulic balance point (or interface) between source nodes. In particular, the flow rates illustrated in \Cref{fig:DfC_vel_Modena_b} suggest that source nodes near the top of the network supply a greater proportion of the demand. Additionally, we note that the optimal DBV locations, which increase self-cleaning velocities by redistributing flow paths, may result in sub-optimal conditions for controlling network pressure. Since we do not explicitly consider pressure management in the SCC problem formulation, the integration of self-cleaning and pressure management objectives should be a focus of future work. Lastly, \Cref{fig:DfC_vel_BWFLnet_b} highlights the limited influence of control features for the more branched \texttt{BWFLnet} network. Although the placement and operation of AFVs indeed improve self-cleaning hydraulic conditions, these are generally limited to direct upstream flow paths. In the case where self-cleaning velocities are sought for specific areas of the network (e.g. historical complaints), the SCC objective function can be modified to focus on a user-defined subset of pipes.

In comparing the two implementations of CMS, with and without domain reduction, we observe the former to yield better quality feasible solutions. This can be attributed to the tighter relaxations resulting from domain reduction. In fact, tighter relaxations may generate fractional valve placement values closer to the feasible domain, and therefore more likely to result in better quality design-for-control solutions. In addition, the solution of  subproblem~\eqref{eq:LP_problem} with tighter relaxations computes control valve settings $\eta$ that are closer to the feasible set and the global optimum. As described in \Cref{sec:multistart_solver}, these $\eta$ values are used as a starting point for the multi-start solver. Hence, we expect these settings to be advantageous in finding a better solution. Nonetheless, the OBBT algorithm comes at an increase in computational effort. Figures ~\ref{fig:DfC_Pescara_b}, \ref{fig:DfC_Modena_b}, and \ref{fig:DfC_BWFLnet_b} compare CPU times across numerical experiments for each solution method. Although the increase in CPU times are modest for the smaller networks, OBBT does not scale well for the larger operational network. For example, the increase in CPU time reported for \texttt{BWFLnet} was, at minimum, approximately 50\% when OBBT was applied. Even though Problem \eqref{eq:MINLP_problem} is considered a design problem, and is thus solved offline, OBBT can become impractical for large-scale and highly interconnected networks \citep{PECCI2022}. Future work should investigate speed improvements to OBBT or the implementation of tighter relaxations for Problem \eqref{eq:MINLP_problem}.

Finally, we compare design-for-control solutions computed by CMS with an off-the-shelf GA implementation. As described in \Cref{sec:results_setup}, the EPANET2.2 hydraulic solver was employed within the GA's fitness function to test hydraulic feasibility of the generated solutions. The solution was deemed feasible when the EPANET2.2 solver successfully converged and when the entire set of pressure heads across all time steps met the minimum pressure constraint. Otherwise, hydraulically infeasible solutions were discarded. Moreover, pressure set-points at control valves were bounded by the maximum head difference across the network (i.e. reservoir level less minimum node elevation). In \Cref{fig:DfC_Pescara_a}, the GA solution is shown to produce comparable (if not better) SCC results for \texttt{Pescara} than CMS (DR). This is not surprising as the stochastic search procedure employed by the GA is well suited for highly nonlinear optimization problems, such as the current SCC objective function. Although the GA did not satisfy its termination criterion within the prescribed time limit (six hours), the relatively small number of free continuous variables in \texttt{Pescara} (e.g. pressure set-points and flushing rates; see \Cref{table:prob_data}) enabled the computation of good-quality feasible solutions. On the other hand, the GA was unable to find feasible solutions for numerous problem instances using \texttt{Modena}. This was particularly evident for experiments including the design and control of AFVs, which result in higher frictional losses materializing from flushing demands. Note that infeasibility is indicated by an SCC value of zero in \Cref{fig:DfC_Modena_a}. Additional GA experiments revealed such infeasibility issues to be linked to the relatively small range in hydraulically feasible pressures found in \texttt{Modena}. With the minimum regulatory pressure head relaxed from 15 m to 0 m (i.e. $h^{\min}$ set to node elevations), the GA found feasible solutions comparable to CMS (DR). These results highlight the efficacy of the proposed convex heuristic for a wider range of network conditions. For \texttt{BWFLnet}, CMS (DR) finds similar quality design-for-control solutions to the GA. However, as observed in \texttt{Modena}, CMS (DR) is shown to consistently outperform the GA for experiments considering the implementation of AFVs. In light of the comparison between solution methods, we investigated the computational overhead of the GA fitness function, as configured in this study using the EPANET-MATLAB Toolkit. This is a cumbersome process, and one which directly impacts the number of fitness evaluations performed by the GA. By comparison, CMS requires much less overhead from loading external software and computes hydraulic states via an efficient null space solver. \Cref{table:fitness_func_eval} compares the number of EPANET2.2 solver calls (for the 12-hour time limit) with the null space solver calls in CMS. The comparison uses select design-for-control experiments for the \texttt{BWFLnet} case study network.

\begin{table}[h!]
    \centering
    \setlength{\tabcolsep}{16pt}
    \caption{Comparison of the number of EPANET2.2 solver calls and null space hydraulic simulations in CMS for \texttt{BWFLnet}}
    \label{table:fitness_func_eval}
        \begin{tabular}{lrr}
            \toprule
            Experiment no. & EPANET2.2 solver & Null space solver\\
            \midrule
            1 --- $n_v = 1$,  $n_f = 0$ & 11,064 & 4,208\\
            3 --- $n_v = 1$,  $n_f = 2$ & 10,304 & 69,995\\
            5 --- $n_v = 2$,  $n_f = 0$ & 11,064 & 61,274\\
            7 --- $n_v = 2$,  $n_f = 2$ & 10,304 & 177,820\\
            9 --- $n_v = 3$,  $n_f = 0$ & 10,304 & 87,976\\
            11 --- $n_v = 3$,  $n_f = 2$ & 10,304 & 165,671\\
            \bottomrule
        \end{tabular}
    \end{table} 

With the exception of Experiment 1, \Cref{table:fitness_func_eval} reports the number of EPANET2.2 solver calls to be significantly less than the hydraulic simulations recorded in CMS. This identifies the hydraulic analysis as a bottleneck in the overall GA search procedure, particularly since its CPU times are over twice that reported using CMS. The comparison demonstrates the advantages of the proposed CMS algorithm in its ability to efficiently perform hydraulic simulations, which in turn yields more opportunities for CMS to compute solutions to the SCC problem than the GA over a given period.

\section{Conclusions}
\label{sec:conclusions}
Novel operational strategies are needed to minimize the severity and frequency of discolouration incidents in WDNs. In this paper, we propose control and design-for-control strategies for maximizing the self-cleaning capacity (SCC) of a network. This is facilitated through the operational framework of dynamically adaptive networks by extending its control capabilities to enhance water quality.

Building on previous work, we formulate an optimization problem to improve SCC through the optimal placement and operation of pressure control and automatic flushing valves. This results in a nonconvex MINLP problem. Since global solvers become intractable for large problem instances, we propose a tailored convex multi-start heuristic (referred to as CMS) to compute feasible solutions to the SCC design-for-control problem. CMS comprises a convex relaxation and randomized valve sampling procedure, which is demonstrated to effectively handle integer variables modelling valve placement. Given the high nonlinearity of the SCC objective function, CMS also implements a multi-start solver to compute feasible solutions to the SCC control problem, whilst mitigating the occurrence of poor local optima. Compared to an off-the-shelf GA implementation, which employs the EPANET2.2 hydraulic solver, CMS is shown to consistently yield good-quality feasible solutions to the SCC design-for-control problem. We demonstrate the robustness of CMS using networks with varying size and hydraulic complexity, including a large-scale operational network in the UK. Although the GA computed comparable solutions to CMS for most numerical experiments, it yielded worse-quality and, in some instances, infeasible solutions for larger and more complex conditions. Furthermore, we show the proposed multi-start solver to compute fast and scalable solutions to the SCC control problem, highlighting its suitability for application in near real-time control strategies. Altogether, the proposed convex multi-start algorithm represents progress towards the optimal design and control of smart water networks. 

Moving forward, the proposed solution for maximizing SCC needs to be integrated within a wider framework for the design and operation of WDNs with dynamically adaptive control. For example, the objective to minimize average zone pressure and pressure variations should include the periodic transition to a self-cleaning mode, as formulated in this manuscript through the operation of pressure control and automatic flushing valves. From an operational perspective, it is also important that this control framework considers the management of unsteady hydraulics arising from such dynamic controls. In relation to the SCC objective function, additional features should also be investigated to improve self-cleaning performance. These may include the impacts of flow reversals as well as a tailored approach to focus on specific areas of the network (e.g. historical complaints). These considerations, combined with advances in technologies for the operation of dynamically adaptive networks, will enable better control strategies for responding to environmental, financial, and regulatory challenges.


\pagebreak
\renewcommand{\thetable}{A.\arabic{table}}
\renewcommand{\thefigure}{A.\arabic{figure}}
\renewcommand{\theequation}{A.\arabic{equation}}
\renewcommand{\thesection}{A\arabic{section}}
\renewcommand{\thealgorithm}{A\arabic{algorithm}}
\section*{A. Convex relaxations}
\label{sec:A}
This section presents a detailed derivation of the convex relaxations formulated for the nonconvex self-cleaning capacity (SCC) objective function and nonconvex HW head loss model. 

\subsection*{A.1 Nonconvex SCC objective function}
\label{sec:A1}
Recalling Section 3.1, we define piecewise linear relaxations for the following nonconvex sigmoidal inequality constraints
\begin{subequations}
    \label{eq:SCC_obj_relax_a}
    \begin{align}
        \label{eq:SCC_obj_relax_a1}
        &\sigma^+_t \leq \psi^+\bigg(\frac{q_{t}}{A}\bigg), \quad \forall t \in \{1,\ldots,n_t\}\\[3pt]
        \label{eq:SCC_obj_relax_a2}
        & \sigma^-_t \leq \psi^-\bigg(\frac{q_{t}}{A}\bigg), \quad \forall t \in \{1,\ldots,n_t\}
    \end{align}
\end{subequations}
where $\psi^+(\cdot)$ and $\psi^-(\cdot)$ model the positive and negative components of the SCC objective function $f_{\widetilde{\text{SCC}}}$, respectively; $\sigma_t := [\sigma_t^+ \,\sigma_t^-]^T$ is the vector of auxiliary variables introduced to reformulate the SCC objective function; and $q_t \in \mathbbm{R}^{n_p}$ is the vector of link flows. As presented in Section 3.1, the resulting relaxations are written as the following linear inequality constraint:
\begin{equation}
    \label{eq:SCC_obj_relax_c}
    S_tq_t + T_t\sigma_t \le s_t, \quad \forall t \in \{1,\ldots,n_t\} 
\end{equation}
where matrices $S_t := [S_t^+ \, S_t^-]^T$ and $T_t := [T_t^+ \, T_t^-]^T$ and vector $s_t := [s_t^+ \, s_t^-]^T$ depend on flow bounds $q_t^L$ and $q_t^U$ and the sigmoid function parameters.

Before proceeding with the derivation of \eqref{eq:SCC_obj_relax_c}, we note the following. First, we omit time indices $t$ and link indices $j$ for the sake of clarity. Second, we drop the sign of the sigmoid function $\psi(\cdot)$ and the vector of auxiliary variables $\sigma \in \mathbbm{R}^{n_p}$; apart from swapping the upper and lower flow bounds, the mathematical derivation does not change for the positive or negative $\psi(\cdot)$ functions. However, where deemed necessary to clearly explain the formulated relaxations, we use the positive sigmoid function. Finally, we use flow velocities as the domain for $\psi(\cdot)$, which are denoted by $u := \big(\frac{q}{A}\big)$. It follows that the lower and upper domain bounds are denoted by $u^L := \big(\frac{q^L}{A}\big)$ and $u^U := \big(\frac{q^U}{A}\big)$, respectively.

We begin by considering the set
\begin{equation}
    \label{eq:sigmoid_relax_a}
    \bigg\{(u,\sigma) \;\big| \; \sigma = \psi(u), \quad \forall u \in \big[u^L\, u^U\big]\bigg\},
\end{equation}
and its concave envelope given by
\begin{equation}
    \label{eq:sigmoid_relax_b}
    \bigg\{(u,\sigma) \;\big| \; 0 \leq \sigma \leq \hat{\psi}(u), \quad \forall u \in \big[u^L\, u^U\big]\bigg\},
\end{equation}
where $\hat{\psi}(\cdot)$ is a piecewise linear relaxation that depends on bounds $u^L$ and $u^U$ and the tangent point created with the sigmoid function. Following the examples presented in \citet{UDELL2014,UDELL2016}, we discuss the steps to construct $\hat{\psi}(\cdot)$ for various domain scenarios below.

\begin{figure}[p]
    \centering
    \captionsetup{justification=centering}
    \subfloat[\label{fig:sigmoid_relax_1}\text{$u^L < w \leq u^U$}]{
        \includegraphics[width=0.45\textwidth]{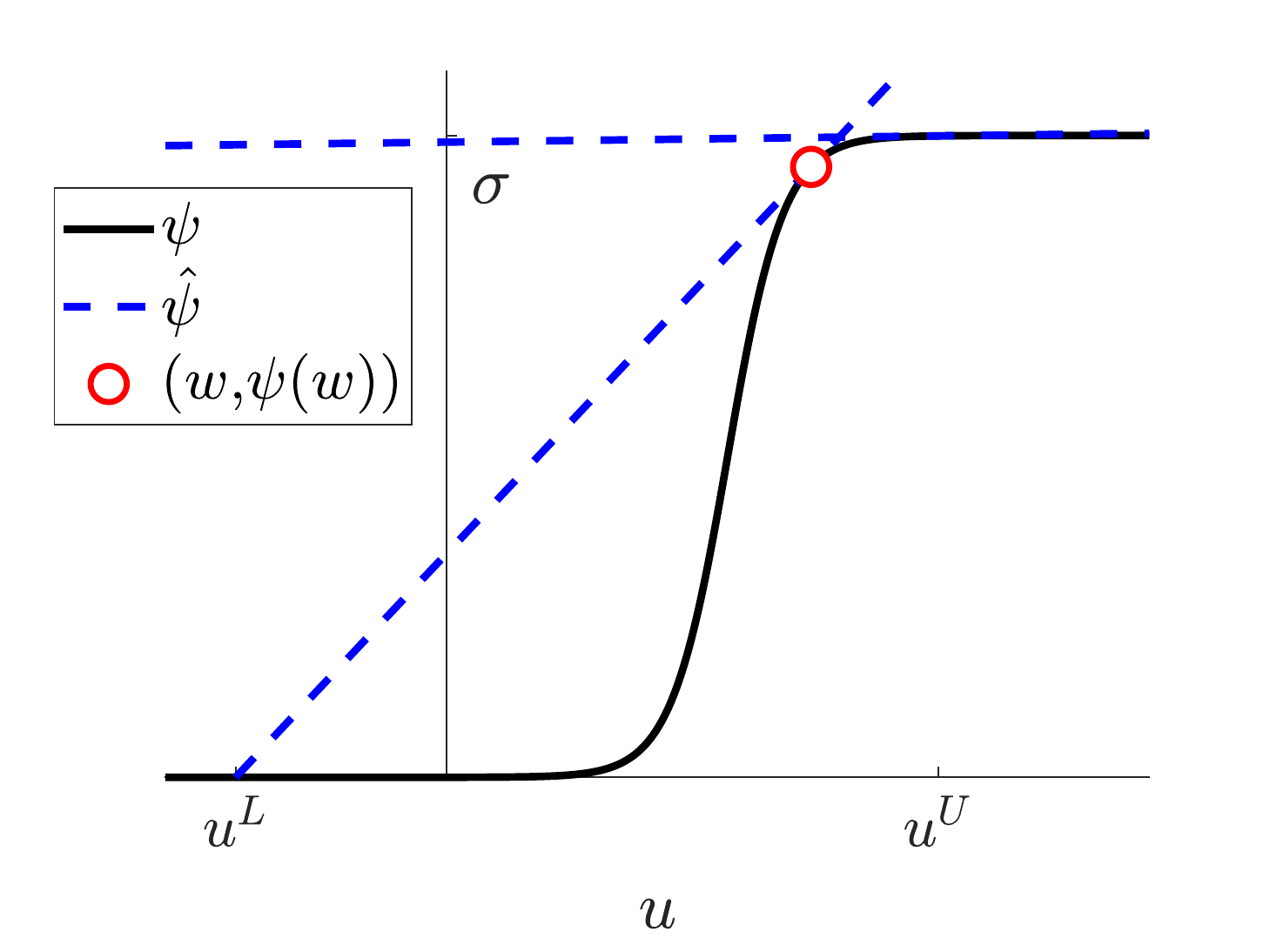}}
    \hspace{0.25cm}
    \subfloat[\label{fig:sigmoid_relax_2}\text{$u^L < u^U \leq w$}]{
        \includegraphics[width=0.45\textwidth]{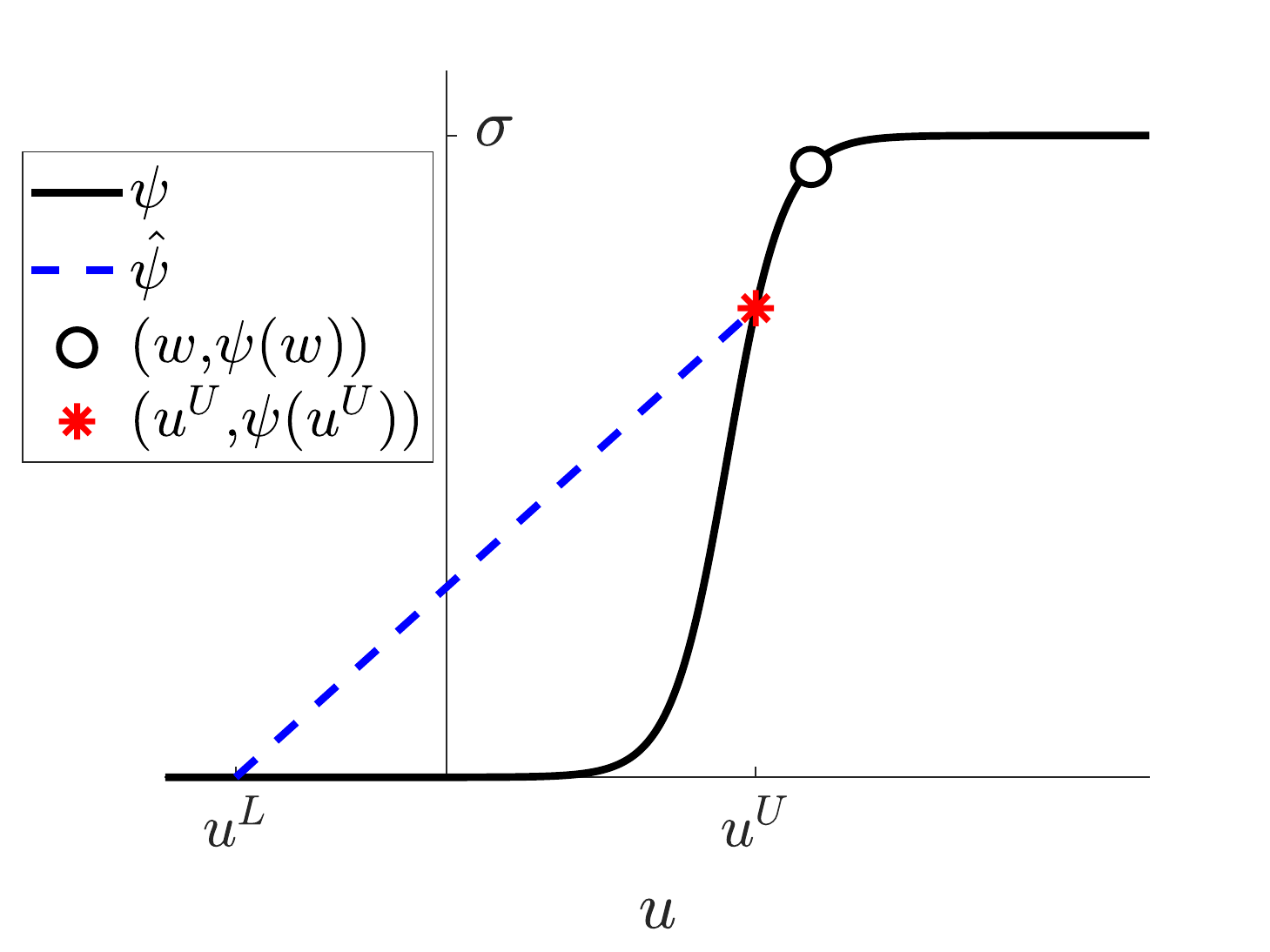}}
    \vspace{1.0cm}    
    \subfloat[\label{fig:sigmoid_relax_3}\text{$u^L = w \leq u^U$}]{
        \includegraphics[width=0.45\textwidth]{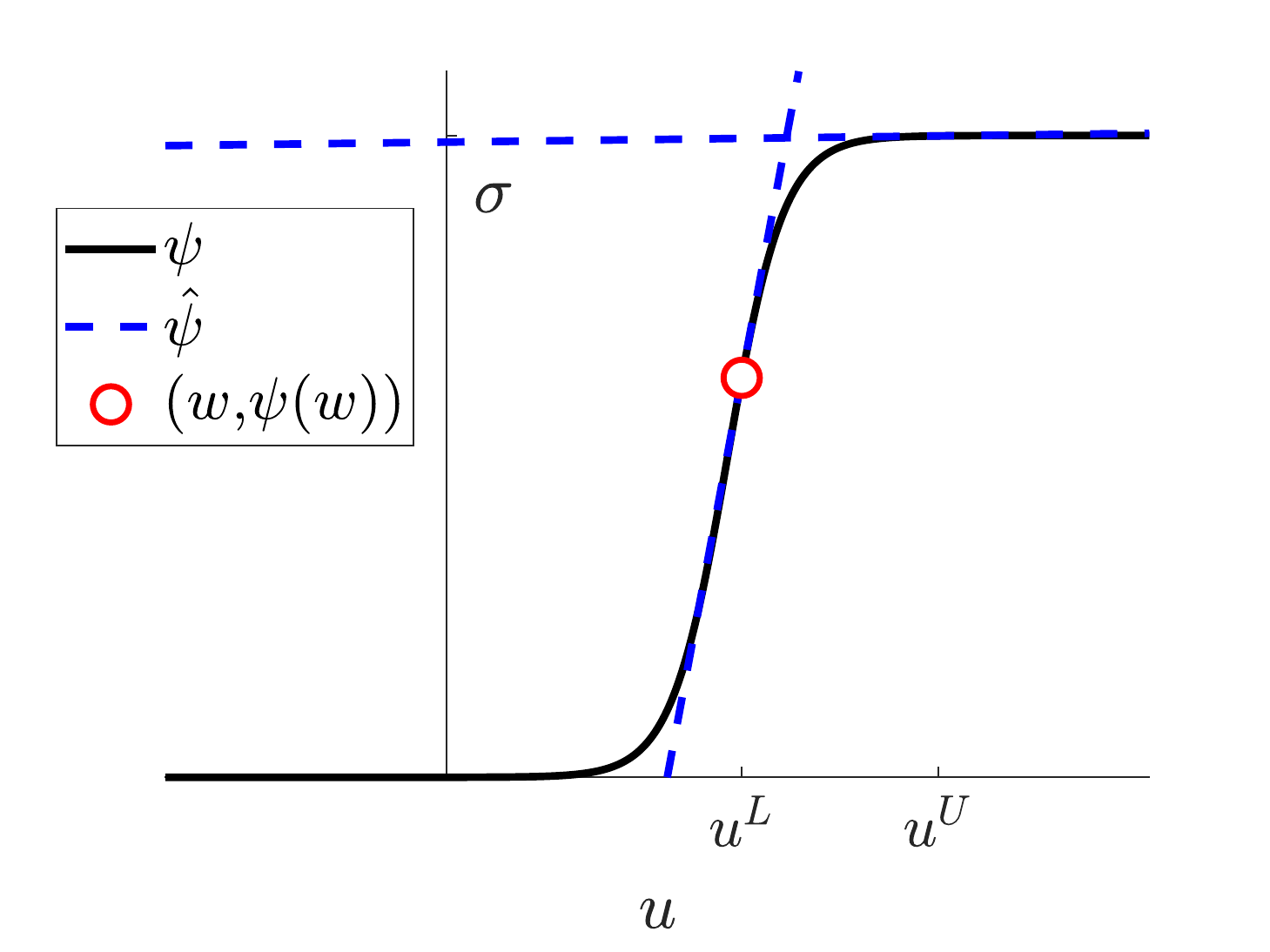}}
    \hspace{0.25cm}
    \subfloat[\label{fig:sigmoid_relax_4}\text{$u^L = w= u^U$}]{
        \includegraphics[width=0.45\textwidth]{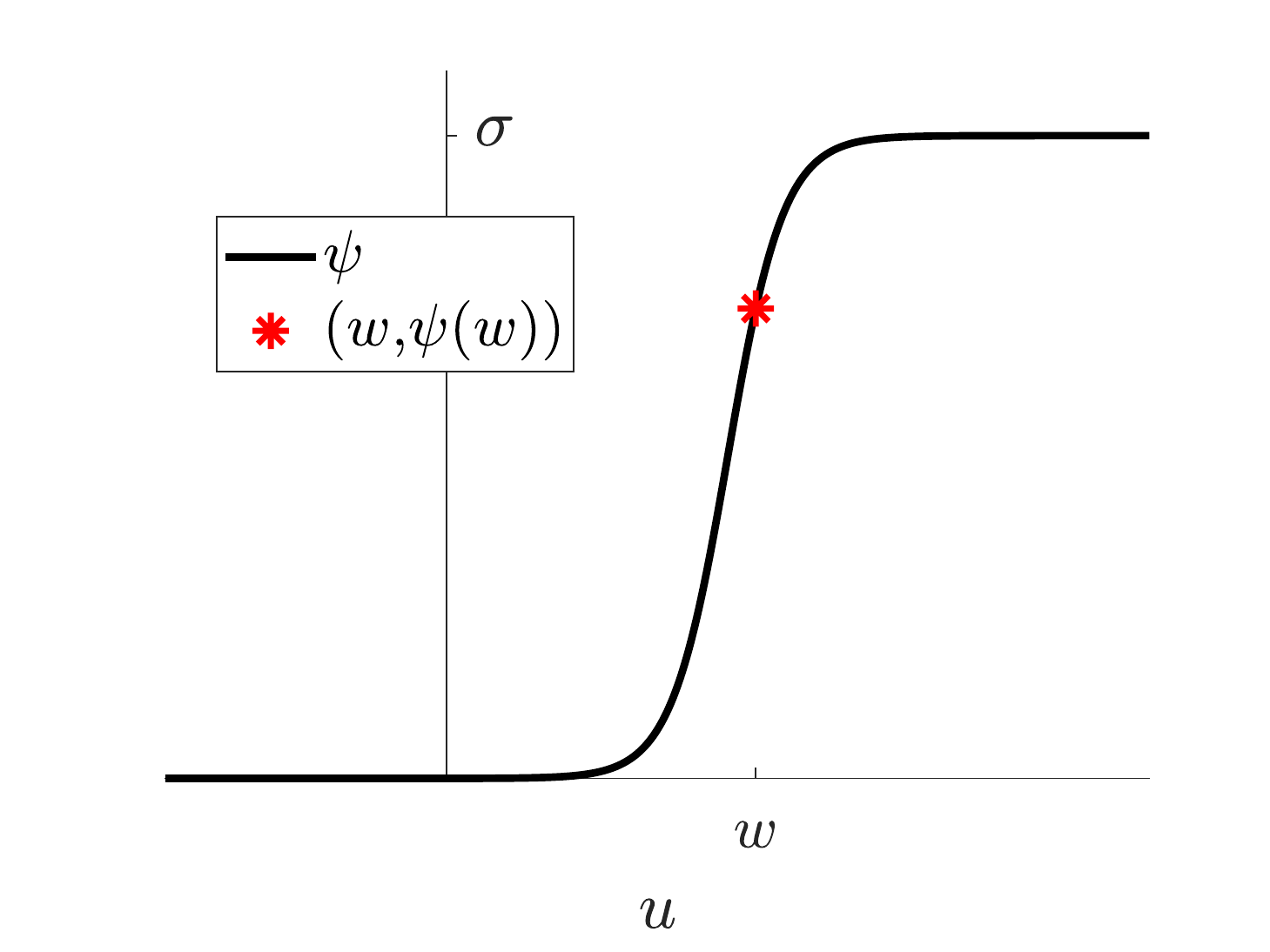}}
    
        \caption{Concave envelopes of nonconvex sigmoidal inequality constraints \eqref{eq:SCC_obj_relax_a}}
    \label{fig:sigmoid_relax}
\end{figure}

First, we look for the point $w$ such that the line from $(u^L,\psi(u^L))$ to $(w,\psi(w))$ is tangent to $\psi(\cdot)$ at $(w,\psi(w))$. The bisection method, an iterative numerical algorithm \citep{MILLER2014}, is used to compute the tangent point $w$ as no analytical expression exists for the equations of intersection in consideration. This follows the steps detailed in \Cref{alg:bisection_sigmoid}.
\begin{algorithm}[h!]
    \caption{Bisection method for nonconvex sigmoid objective function}
    \label{alg:bisection_sigmoid}
    \begin{algorithmic}[1]
    \State \textbf{Input:} $\psi(\cdot)$ model parameters, SCC threshold velocity $u^{\min}$, and flow velocity bounds $u^L$ and $u^U$
    \State \textbf{Output:} unique point $w$ where the line from $(y,\psi(y))$ to $(w,\psi(w))$ is tangent to $\psi(\cdot)$ at $(w,\psi(w))$
    \State \textbf{Define:} for $\psi^+$, $y := u^L$; for $\psi^-$, $y := u^U$; equation of intersection $f(x) := \psi'(x)(x-y) + \psi(y) - \psi(x)$, where $x = w$ at the root of $f(x)$
    \State Initialize: (i) initial lower and upper bounds $x_0$ and $x_1$; and (ii) error tolerance $\epsilon_{\text{tol}}$
    \If {$u^L \geq u^{\min}$}
        \State $w \; \leftarrow \; u^L$
    \Else
        \State $x_2 \leftarrow \frac{x_0+x_1}{2}$
        \While {$|x_1 - x_0| \geq \epsilon_{tol}$}
            \If {$f(x_2) \times f(x_0) \leq 0$}
                \State $x_1 \; \leftarrow \;  x_2$
                \State $x_2 \; \leftarrow \;  \frac{x_0+x_1}{2}$
            \Else
                \State $x_0 \; \leftarrow \;  x_2$
                \State $x_2 \; \leftarrow \;  \frac{x_0+x_1}{2}$
            \EndIf
        \EndWhile
        \State $w \; \leftarrow \;  x_2$
    \EndIf
    \end{algorithmic}
\end{algorithm}
Since the SCC problem is mostly comprised of linear constraints, it is also convenient to construct a linear upper bound to the concave segment of $\psi(\cdot)$. It follows that point $k$ represents the intersection of linear relaxations, at which point the maximum error of $\hat{\psi}(\cdot) - \psi(\cdot)$ occurs. The constructed concave envelopes $\hat{\psi}(\cdot)$ are described below for the four conceivable domain scenarios, with their corresponding plots shown in \Cref{fig:sigmoid_relax}.

If $u^L < w \leq u^U$, the concave envelope $\hat{\psi}(u)$ is defined as the following piecewise linear function (see \Cref{fig:sigmoid_relax_1}):
\begin{equation} 
    \label{eq:sigmoid_relax_c1}
    \begin{alignedat}{3}
        &\hat{\psi}(u) =
        \begin{cases} 
            \psi(w) + \psi'(w)(u-w), & u^L \leq u \leq k \\[3pt]
            \psi(u^U) + \psi'(u^U)(u-u^U), & k \leq u \leq u^U
        \end{cases}
    \end{alignedat}
\end{equation}

\noindent
When $u^L < u^U \leq w$, the concave envelope $\hat{\psi}(u)$ is defined a the following linear function (see \Cref{fig:sigmoid_relax_2}):
\begin{equation} 
    \label{eq:sigmoid_relax_c2}
    \begin{alignedat}{3}
        &\hat{\psi}(u) = \frac{\psi(u^U) - \psi(u^L)}{u^U - u^L}(u-u^L), \quad u^L \leq u \leq u^U
    \end{alignedat}
\end{equation}

\noindent
When $u^L =w$ and $w \leq u^U$, the concave envelope $\hat{\psi}(u)$ is defined as the following piecewise linear function (see \Cref{fig:sigmoid_relax_3}):
\begin{equation} 
    \label{eq:sigmoid_relax_c3}
    \begin{alignedat}{3}
        &\hat{\psi}(u) =
        \begin{cases} 
            \psi(u^L) + \psi'(u^L)(u-u^L), & u^L \leq u \leq k \\[3pt]
            \psi(u^U) + \psi'(u^U)(u-u^U), & k \leq u \leq u^U
        \end{cases}
    \end{alignedat}
\end{equation}

\noindent
Lastly, if $u^L = w = u^U$, the concave envelope $\hat{\psi}(u)$ is defined as the point (see \Cref{fig:sigmoid_relax_4}):
\begin{equation} 
    \label{eq:sigmoid_relax_c4}
    \begin{alignedat}{3}
        &\hat{\psi}(u) = \psi(w)
    \end{alignedat}
\end{equation}

\noindent
Sigmoid functions $\psi^+$ and $\psi^-$ in \eqref{eq:SCC_obj_relax_a} are then replaced with the appropriate concave envelopes $\hat{\psi}(\cdot)$ described in \eqref{eq:sigmoid_relax_c1} - \eqref{eq:sigmoid_relax_c4}. After rearranging terms, and substituting flow velocity $u$ with $\big(\frac{q}{A}\big)$ to align with the SCC problem variables, the linear inequality constraint shown in \eqref{eq:SCC_obj_relax_c} is formulated.

\vspace{0.1cm}

\subsection*{A2. Nonconvex HW head loss constraints}
\label{sec:A2}
Building on the material presented in \citet{PECCI2019} (Appendix 1), we formulate polyhedral relaxations for the HW head loss model. Recall the energy conservation equality constraint
\begin{equation}
    \label{eq:nonconvex_head_loss}
    \theta_t - \phi(q_t) = 0, \quad \forall t \in \{1,\ldots,n_t\}
\end{equation}
where $\phi(\cdot)$ models the head loss materialized across network links $j \in \{1,\dots,n_p\}$; $\theta_t \in \mathbbm{R}^{n_p}$ is the vector of auxiliary variables introduced to isolate the nonconvex head loss term; and $q_t \in \mathbbm{R}^{n_p}$ is the vector of link flows. The resulting polyhedral relaxations are written as the following linear inequality constraint:
\begin{equation}
    \label{eq:HW_relax_a}
    R_tq_t + E_t\theta_t \le r_t, \quad \forall t \in \{1,\ldots,n_t\} 
\end{equation}
where matrices $R_t$ and $E_t$ and vector $r_t$ depend on flow bounds $q_t^L$ and $q_t^U$ and the HW model parameters. We follow a similar procedure to that detailed in \citet{LIBERTI2003} for monomials of odd degree and later adapted in \citet{PECCI2019} for the quadratic approximation head loss model. The contribution of this work concerns the derivation of polyhedral relaxations for the HW head loss model, which has a fractional exponent dictating the relationship between $q$ and $\phi(\cdot)$.  

Omitting indices $t$ and $j$, we begin by considering the set
\begin{equation}
    \label{eq:HW_relax_b1}
    \bigg\{(q,\theta) \;\big| \; \theta = \phi(q), \quad \forall q \in \big[q^L\, q^U\big]\bigg\},
\end{equation}
and its convex relaxation given by
\begin{equation}
    \label{eq:HW_relax_b2}
    \bigg\{(q,\theta) \;\big| \; \text{$\underline{\phi}$}(q) \leq \theta \leq \bar{\phi}(q), \quad \forall q \in \big[q^L\, q^U\big]\bigg\},
\end{equation}
where $\underline{\phi}(\cdot)$ and $\bar{\phi}(\cdot)$ are piecewise linear functions which formulate lower and upper relaxations of $\phi(\cdot)$, respectively. These functions depend on flow bounds $q^L$ and $q^U$ and the tangent points created with $\phi(\cdot)$. In the following, we present the derivations of $\underline{\phi}(\cdot)$ and $\bar{\phi}(\cdot)$ for the complete range of domain scenarios.

\begin{figure}[p]
    \centering
    \captionsetup{justification=centering}
    \subfloat[\label{fig:HW_relax_1}\text{$q^L < \bar{z} < 0 < \underline{z} < q^U$}]{
        \includegraphics[width=0.45\textwidth]{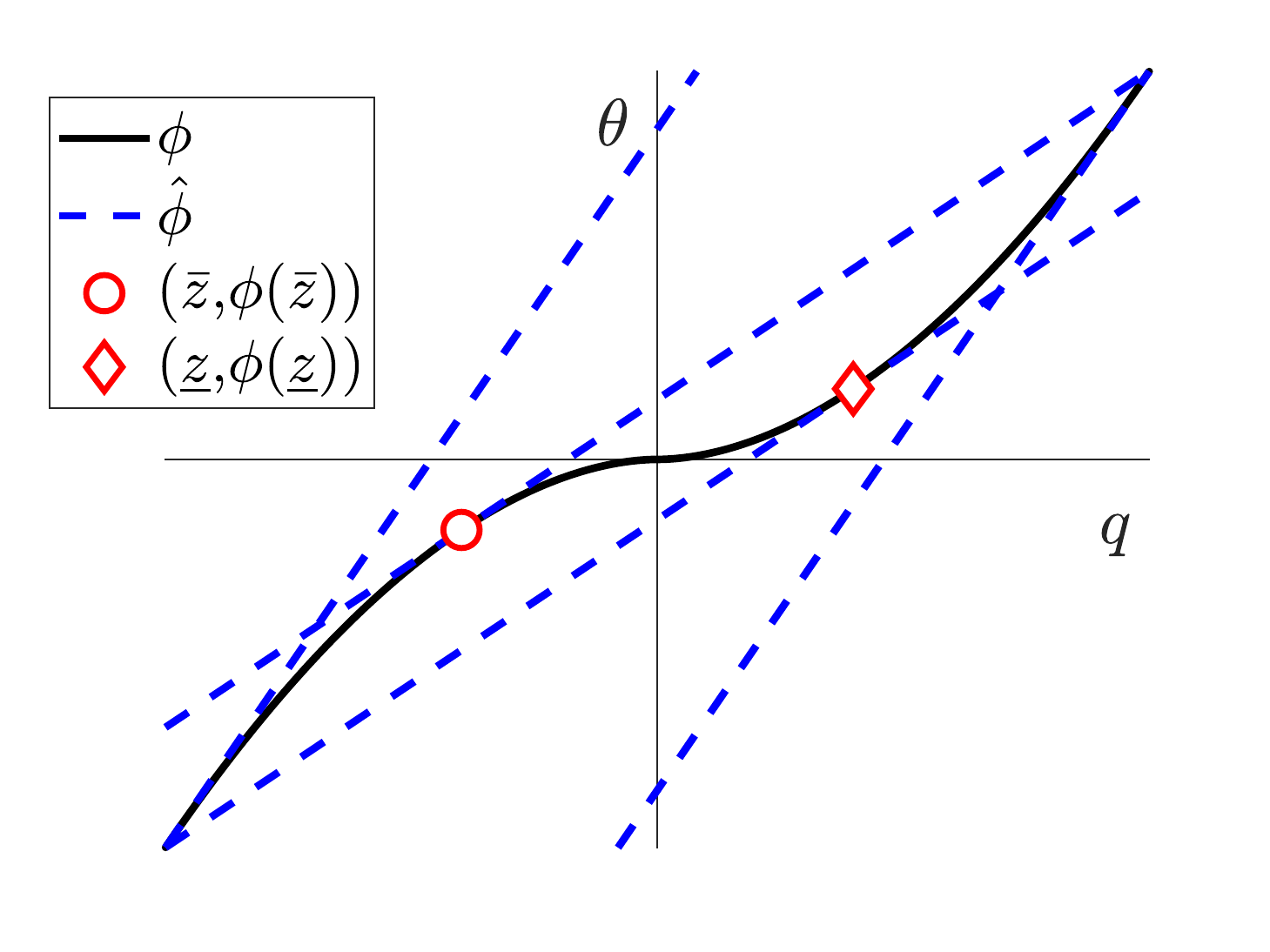}}
    \hspace{0.1cm}
    \subfloat[\label{fig:HW_relax_2}\text{$\bar{z} \leq q^L < 0 < \underline{z} < q^U$}]{
        \includegraphics[width=0.45\textwidth]{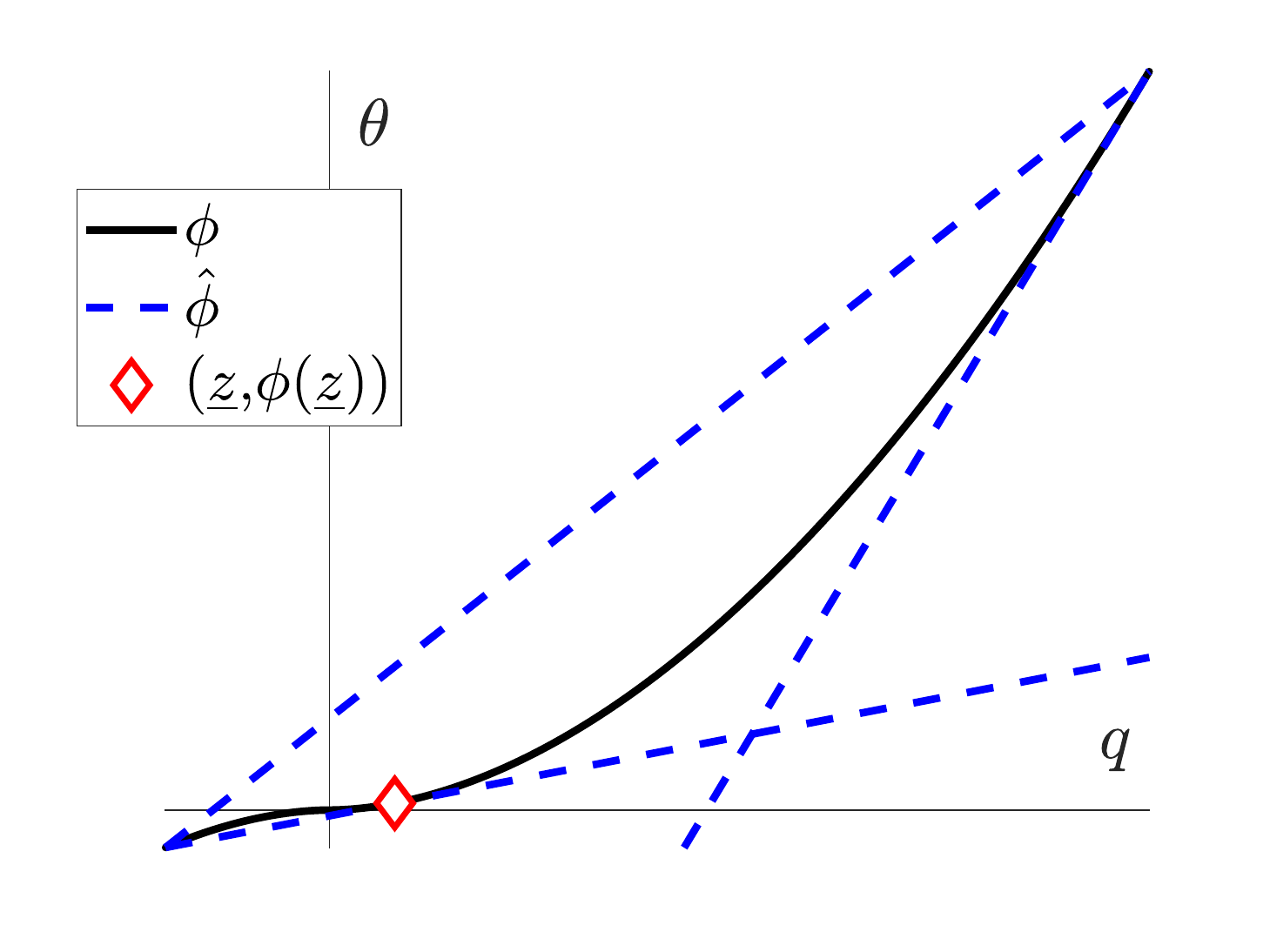}}
    \vspace{0.25cm}
    \subfloat[\label{fig:HW_relax_3}\text{$q^L < \bar{z} < 0 < q^U \leq \underline{z}$}]{
        \includegraphics[width=0.45\textwidth]{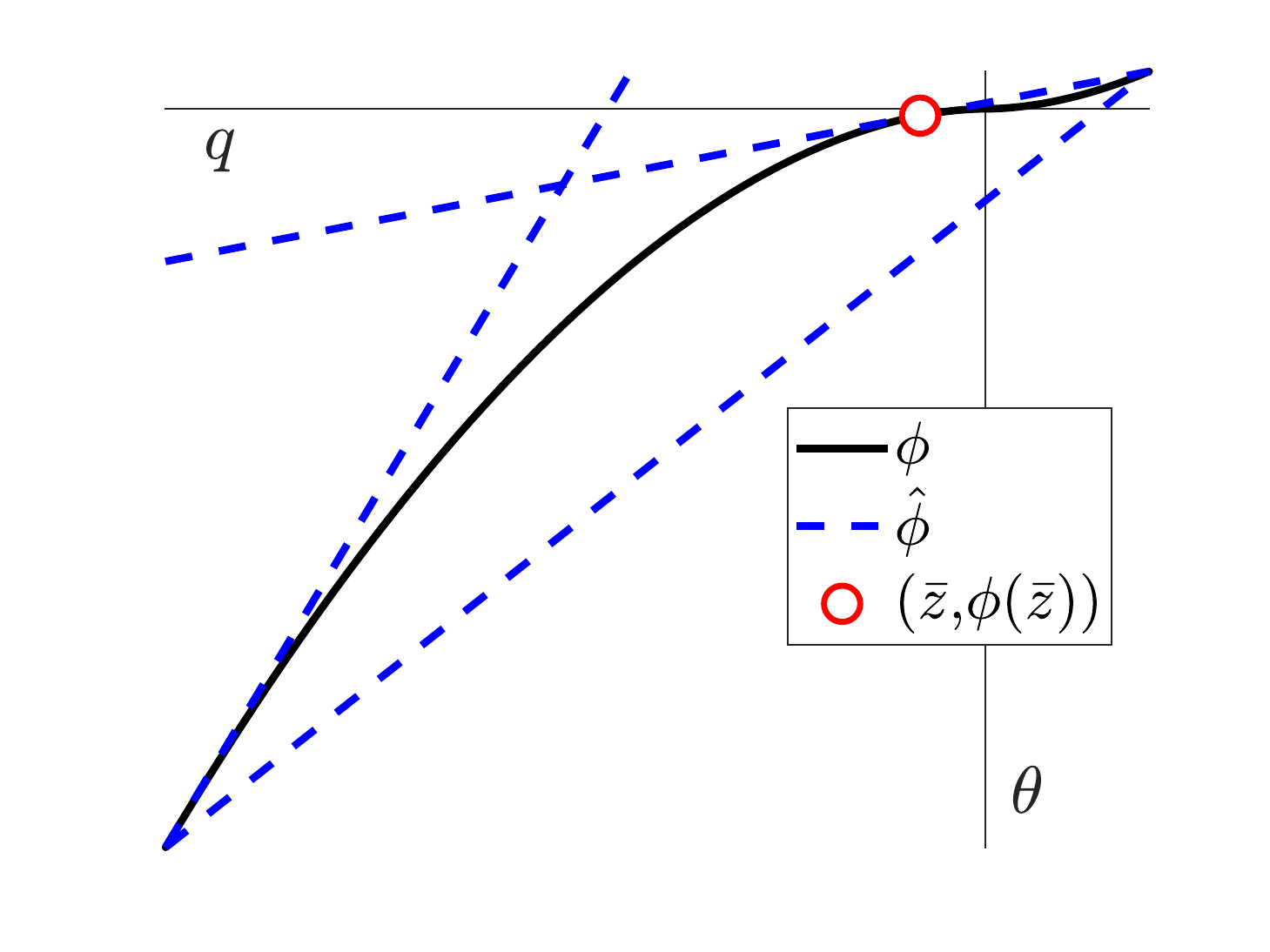}}
    \hspace{0.1cm}    
    \subfloat[\label{fig:HW_relax_4}\text{$0 \leq q^L < q^U$}]{
        \includegraphics[width=0.45\textwidth]{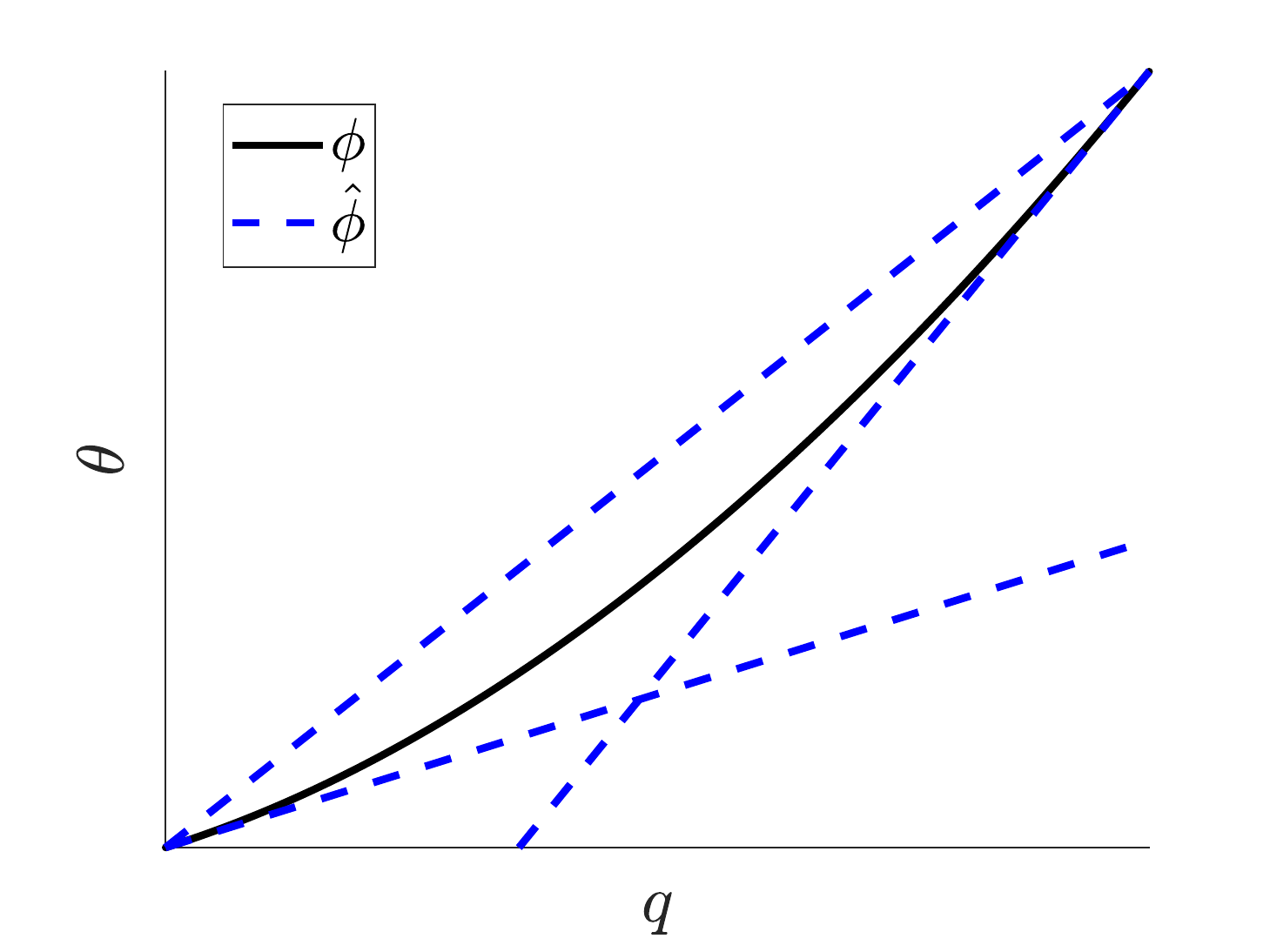}}
    \vspace{0.25cm}
    \subfloat[\label{fig:HW_relax_5}\text{$q^L < q^U \leq 0$}]{
        \includegraphics[width=0.45\textwidth]{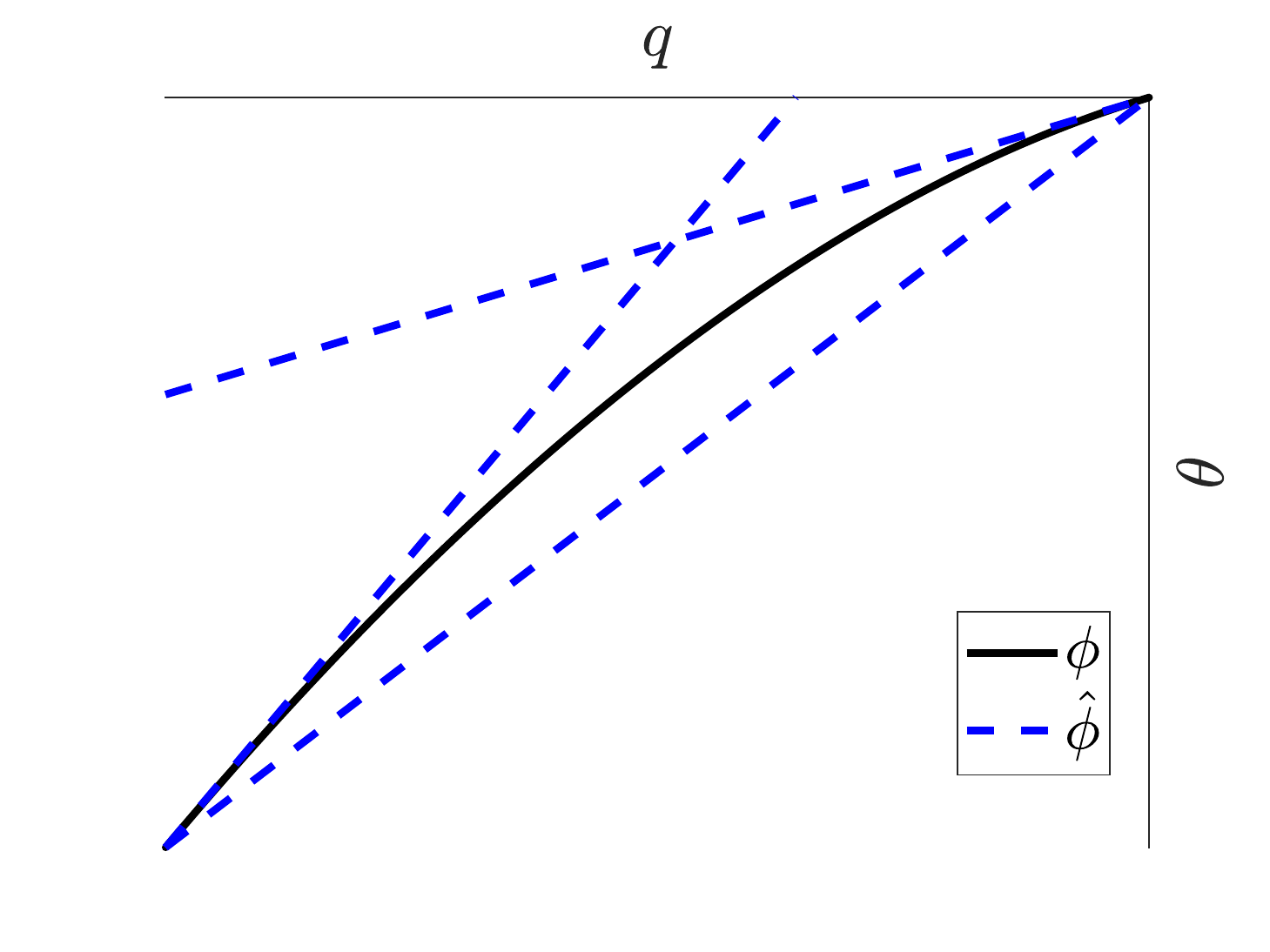}}
    
        \caption{Polyhedral relaxations of nonconvex HW head loss constraints \eqref{eq:nonconvex_head_loss}}
    \label{fig:HW_relax}
\end{figure}

We start by looking for the point $\underline{z}$ such that the line from $(q^L,\phi(q^L))$ to $(\underline{z},\phi(\underline{z}))$ is tangent to $\phi(\cdot)$ at $(\underline{z},\phi(\underline{z}))$. In contrast to the quadratic approximation method derived in \citet{PECCI2019}, there is no analytical expression for finding root $\underline{z}$ to the intersection between $\phi(\cdot)$ and the tangent line at $(\underline{z},\phi(\underline{z}))$. Thus, the bisection method is used to find $\underline{z}$ within an appropriate numerical tolerance. While similar to \Cref{alg:bisection_sigmoid}, we present the exact steps taken to compute $\underline{z}$ in \Cref{alg:bisection_HW} for the sake of completeness. Analogously, the unique point $\bar{z}$ is computed to find the line from $(q^U,\phi(q^U))$ to $(\bar{z},\phi(\bar{z}))$. The resulting set of linear relaxations for $\phi(\cdot)$ are described below, with the corresponding plots presented in \Cref{fig:HW_relax}. Note that we define the set of polyhedral relaxations as $\hat{\phi}(\cdot) := \{\underline{\phi} \cup \bar{\phi}\}$. Moreover, let $\underline{k}$ and $\bar{k}$ denote the points of intersection between the linear functions in $\underline{\phi}(\cdot)$ and $\bar{\phi}(\cdot)$, respectively.

\begin{algorithm}[h!]
    \caption{Bisection method for nonconvex HW head loss model}
    \label{alg:bisection_HW}
    \begin{algorithmic}[1]
    \State \textbf{Input:} $\phi(\cdot)$ model parameters and flow bounds $q^L$ and $q^U$
    \State \textbf{Output:} unique point $z$ where the line from $(y,\phi(y))$ to $(z,\phi(z))$ is tangent to $\phi(\cdot)$ at $(z,\phi(z))$
    \State \textbf{Define:} for $\underline{\phi}$, $z := \underline{z}$, $y_1 := q^L$, and $y_2 := q^U$; for $\bar{\phi}$, $z := \bar{z}$, $y_1 := q^U$, and $y_2 := q^L$; equation of intersection $f(x) := \phi'(x)(x-y_1) + \phi(y_1) - \phi(x)$, where $x = z$ at the root of $f(x)$
    \State Initialize: (i) initial lower and upper bounds $x_0$ and $x_1$; and (ii) error tolerance $\epsilon_{\text{tol}}$
    \If {$f(y_2) \times f(0) \leq 0$}
        \State $x_2 \leftarrow \frac{x_0+x_1}{2}$
        \While {$|x_1 - x_0| \geq \epsilon_{tol}$}
            \If {$f(x_2) \times f(x_0) \leq 0$}
                \State $x_1 \; \leftarrow \;  x_2$
                \State $x_2 \; \leftarrow \;  \frac{x_0+x_1}{2}$
            \Else
                \State $x_0 \; \leftarrow \;  x_2$
                \State $x_2 \; \leftarrow \;  \frac{x_0+x_1}{2}$
            \EndIf
        \EndWhile
        \State $z \; \leftarrow \; x_2$
    \Else
        \State $z \; \leftarrow \;$  NaN
    \EndIf
    \end{algorithmic}
\end{algorithm}

If $q^L < \bar{z} < 0 < \underline{z} < q^U$, the upper and lower bounds on $\theta$ are defined as the following piecewise linear functions (see \Cref{fig:HW_relax_1}):
\begin{equation} 
    \label{eq:HW_relax_c1}
    \begin{alignedat}{3}
        &\bar{\phi}(q) =
        \begin{cases} 
            \phi(q^L) + \phi'(q^L)(q-q^L), & q^L \leq q \leq \bar{k} \\[3pt]
            \phi(q^U) + \phi'(\bar{z})(q-q^U), & \bar{k} \leq q \leq q^U
        \end{cases}\\[9pt]
        &\underline{\phi}(q) =
        \begin{cases} 
            \phi(q^L) + \phi'(\underline{z})(q-q^L), & q^L \leq q \leq \underline{k} \\[3pt]
            \phi(q^U) + \phi'(q^U)(q-q^U), & \underline{k} \leq q \leq q^U
        \end{cases}
    \end{alignedat}
\end{equation}

\noindent
When $\bar{z} \leq q^L < 0 < \underline{z} < q^U$, $\bar{\phi}(\cdot)$ and $\underline{\phi}(\cdot)$ are defined as follows (see \Cref{fig:HW_relax_2}):
\begin{equation} 
    \label{eq:HW_relax_c2}
    \begin{alignedat}{3}
        &\bar{\phi}(q) = \frac{\phi(q^U) - \phi(q^L)}{q^U - q^L}(q-q^L), \quad q^L \leq q \leq q^U\\[9pt]
        &\underline{\phi}(q) =
        \begin{cases} 
            \phi(q^L) + \phi'(\underline{z})(q-q^L), & q^L \leq q \leq \underline{k} \\[3pt]
            \phi(q^U) + \phi'(q^U)(q-q^U), & \underline{k} \leq q \leq q^U
        \end{cases}
    \end{alignedat}
\end{equation}

\noindent
When $q^L < \bar{z} < 0 < q^U \leq \underline{z}$, $\bar{\phi}(\cdot)$ and $\underline{\phi}(\cdot)$ are defined as follows (see \Cref{fig:HW_relax_3}):
\begin{equation} 
    \label{eq:HW_relax_c3}
    \begin{alignedat}{3}
        &\bar{\phi}(q) =
        \begin{cases} 
            \phi(q^L) + \phi'(q^L)(q-q^L), & q^L \leq q \leq \bar{k} \\[3pt]
            \phi(q^U) + \phi'(\bar{z})(q-q^U), & \bar{k} \leq q \leq q^U
        \end{cases}\\[9pt]
        &\underline{\phi}(q) = \frac{\phi(q^U) - \phi(q^L)}{q^U - q^L}(q-q^L), \quad q^L \leq q \leq q^U
    \end{alignedat}
\end{equation}

\noindent
When $0 \leq q^L < q^U$, $\bar{\phi}(\cdot)$ and $\underline{\phi}(\cdot)$ are restricted to the positive domain and defined as follows (see \Cref{fig:HW_relax_4}):
\begin{equation} 
    \label{eq:HW_relax_c4}
    \begin{alignedat}{3}
        &\bar{\phi}(q) = \frac{\phi(q^U) - \phi(q^L)}{q^U - q^L}(q-q^L), \quad q^L \leq q \leq q^U\\[9pt]
        &\underline{\phi}(q) =
        \begin{cases} 
            \phi(q^L) + \phi'(q^L)(q-q^L), & q^L \leq q \leq \underline{k} \\[3pt]
            \phi(q^U) + \phi'(q^U)(q-q^U), & \underline{k} \leq q \leq q^U
        \end{cases}
    \end{alignedat}
\end{equation}

\noindent
Lastly, if $q^L < q^U \leq 0$, $\bar{\phi}(\cdot)$ and $\underline{\phi}(\cdot)$ are restricted to the negative domain and defined as follows (see \Cref{fig:HW_relax_5}):
\begin{equation} 
    \label{eq:HW_relax_c5}
    \begin{alignedat}{3}
        &\bar{\phi}(q) =
        \begin{cases} 
            \phi(q^L) + \phi'(q^L)(q-q^L), & q^L \leq q \leq \bar{k} \\[3pt]
            \phi(q^U) + \phi'(q^U)(q-q^U), & \bar{k} \leq q \leq q^U
        \end{cases}\\[9pt]
        &\underline{\phi}(q) = \frac{\phi(q^U) - \phi(q^L)}{q^U - q^L}(q-q^L), \quad q^L \leq q \leq q^U
    \end{alignedat}
\end{equation}

\noindent
The HW head loss model $\phi(\cdot)$ in \eqref{eq:nonconvex_head_loss} is then replaced with the appropriate polyhedral relaxation $\hat{\phi}(\cdot)$ described in \eqref{eq:HW_relax_c1} - \eqref{eq:HW_relax_c5}. After rearranging terms, the linear inequality constraint in \eqref{eq:HW_relax_a} is formulated.

\vspace{0.1cm}

\pagebreak
\setcounter{figure}{0}
\setcounter{equation}{0}
\setcounter{algorithm}{0}
\renewcommand{\thetable}{B.\arabic{table}}
\renewcommand{\thefigure}{B.\arabic{figure}}
\renewcommand{\theequation}{B.\arabic{equation}}
\renewcommand{\thesection}{B\arabic{section}}
\renewcommand{\thealgorithm}{B\arabic{algorithm}}
\section*{B. Optimization-based bound tightening algorithm}
\label{sec:B}

In Section 3.1, we implement an optimization-based bound tightening (OBBT) algorithm to strengthen the convex relaxations formulated for Subproblem (LP). Following the supplementary material from \citet{PECCI2022}, we tighten polyhedral relaxations for the nonconvex constraints by reducing the flow variable domain for link indices corresponding to the set of core links $\mathcal{C}$. For each OBBT iteration, we solve two optimization problems to independently maximize and minimize flow $q_{j,t}$ for all links $j \in \mathcal{C}$ and hydraulic time steps $t \in \{1,\dots,n_t\}$. These are represented by the linear program described in Subproblem (LP), with objective functions set to both minimize and maximize $q_{j,t}$. We therefore solve $2n_t|\mathcal{C}|$ linear programs for each OBBT iteration. The algorithm is terminated once bound tightening progress diminishes between successive iterations, which is defined by a set tolerance on the maximum flow domain diameter. Moreover, we exclude forest links from the optimization scheme since their flows are defined \textit{a priori} by the aggregate of downstream demands and the maximum flow defined for automatic flushing valves. Psuedocode for the OBBT algorithm is described in \Cref{alg:obbt_algorithm}.

\begin{algorithm}
    \caption{OBBT}
    \label{alg:obbt_algorithm}
    \begin{algorithmic}[1]
    \State \textbf{Input:} network data and initial flow bounds $q^L$ and $q^U$
    \State \textbf{Output:} tightened flow bounds $q^L$ and $q^U$
    \State Initialize $k_{\max}$, $\epsilon_{\text{tol}}$, and set $k = 1$, $\epsilon = 0$
    \State Initialize empty set $\texttt{diam} \;\leftarrow\; \emptyset$
    \State Compute initial maximum flow domain diameter: $\texttt{diam}(1) = \text{max}_{t \in \{1,\dots,n_t\}} \text{max}_{j \in \mathcal{C}}(q_{j,t}^U - q_{j,t}^L)$  
    \While{$\epsilon_{\text{tol}} \geq \epsilon$ \textbf{and} $k \leq k_{\max}$}
        \For{$t \in \{1,\dots,n_t\}$}
            \For{$j \in \mathcal{C}$}
            \State Solve Subproblem (LP) with objective function set to minimize $q_{t,j}$    
            \State Update lower flow bound $q_{t,j}^L$
            \State Solve Subproblem (LP) with objective function set to maximize $q_{t,j}$
            \State Update upper flow bound $q_{t,j}^U$
            \EndFor
        \EndFor
        \State Compute $\texttt{diam}(k+1) = \text{max}_{t \in \{1,\dots,n_t\}} \text{max}_{j \in \mathcal{C}}(q_{j,t}^U - q_{j,t}^L)$
        \State $\epsilon \;\leftarrow\; \frac{\texttt{diam}(k+1)}{\texttt{diam}(k)}$ 
        \State $k \;\leftarrow\; k+1$
    \EndWhile
    \end{algorithmic}
\end{algorithm}

\pagebreak
\setcounter{figure}{0}
\setcounter{equation}{0}
\setcounter{algorithm}{0}
\renewcommand{\thetable}{C.\arabic{table}}
\renewcommand{\thefigure}{C.\arabic{figure}}
\renewcommand{\theequation}{C.\arabic{equation}}
\renewcommand{\thesection}{C\arabic{section}}
\renewcommand{\thealgorithm}{C\arabic{algorithm}}
\section*{C. Strictly feasible sequential convex programming algorithm}
\label{sec:C}

In Section 3.3.1, we implement a strictly feasible sequential convex programming (SFSCP) solver to compute feasible solutions to the nonlinear programming (NLP) control problem. This follows the implementations in \citet{WRIGHT2015} and \citet{ABRAHAM2016}, with the addition of the continuous decision variable $\alpha$ to model the operation of flushing valves. We present the SFSCP solver psuedocode in \Cref{alg:sfscp_solver}, which is adapted from \citet{WRIGHT2015}. The variables used in \Cref{alg:sfscp_solver} are defined as follows: $f$ denotes the continuous sum of sigmoids SCC objective function $f_{\widetilde{\text{SCC}}}$; $\eta$ and $\alpha$ model operational settings for pressure control and flushing valves, respectively; $d\eta$ and $d\alpha$ denote step directions for the aforementioned continuous variables; $q$ and $h$ represent hydraulic states computed for each control configuration; and $\mathcal{F}$ is the set defining the hydraulically feasible region.

\begin{algorithm}
    \caption{SFSCP solver}
    \label{alg:sfscp_solver}
    \begin{algorithmic}[1]
    \State \textbf{Input:} network data, fixed valve configuration, and feasible starting point $x_0$
    \State \textbf{Output:} $f^*$ and corresponding operational control settings $\eta^*$ and $\alpha^*$ 
    \State Initialize $\epsilon_{\text{tol}}$, $k_{\max}$, and set $k = 1$
    \State Set $x_k \;\leftarrow\; x^0$
    \State Compute initial hydraulic states $q_k, h_k$, and objective function $f_k$
    \While {$\frac{|f_k - f_{k+1}|}{|f_k|} > \epsilon_{\text{tol}}$ \textbf{and} $k \leq k_{\max}$}
        \State $\eta_{k+1}, \alpha_{k+1} \;\leftarrow\;$ NLP control problem with first-order Taylor approximations of nonlinear SCC objective function and energy conservation constraints
        \State $d\eta_{k}, d\alpha_{k} \;\leftarrow\; \eta_{k+1} - \eta_k, \alpha_{k+1} - \alpha_k$
        \State $\beta = 1$
        \State $q_{k+1}, h_{k+1} \; \leftarrow \;$ hydraulic simulation with controls $\eta_{k} + \beta d\eta_k, \alpha_{k} + \beta d\alpha_k$
        \State $f_{k+1} \; \leftarrow \;$ recompute objective function
        \While {$f_k > f_{k+1}$ \textbf{or} $x_{k+1} \notin \mathcal{F}$}
            \State $\beta \;\leftarrow\; \beta/2$
            \State $\eta_{k+1}, \alpha_{k+1} \;\leftarrow\; \eta_{k} + \beta d\eta_k, \alpha_{k} + \beta d\alpha_k$
            \State $q_{k+1}, h_{k+1} \; \leftarrow \;$ hydraulic simulation with controls $\eta_{k} + \beta d\eta_k, \alpha_{k} + \beta d\alpha_k$
            \State $f_{k+1} \; \leftarrow \;$ recompute objective function
        \EndWhile
        \State $k \;\leftarrow\; k+1$
        \State $\eta_k, \alpha_k \; \leftarrow \; \eta_{k+1}, \alpha_{k+1}$
        \State $f_k \;\leftarrow\; f_{k+1}$
    \EndWhile
    \end{algorithmic}
\end{algorithm}

\pagebreak
\setcounter{figure}{0}
\setcounter{equation}{0}
\setcounter{algorithm}{0}
\renewcommand{\thetable}{D.\arabic{table}}
\renewcommand{\thefigure}{D.\arabic{figure}}
\renewcommand{\theequation}{D.\arabic{equation}}
\renewcommand{\thesection}{D\arabic{section}}
\section*{D. NLP solver performance comparison}
\label{sec:D}

\begin{figure}[h!]
    \centering
    \includegraphics[width=0.55\linewidth]{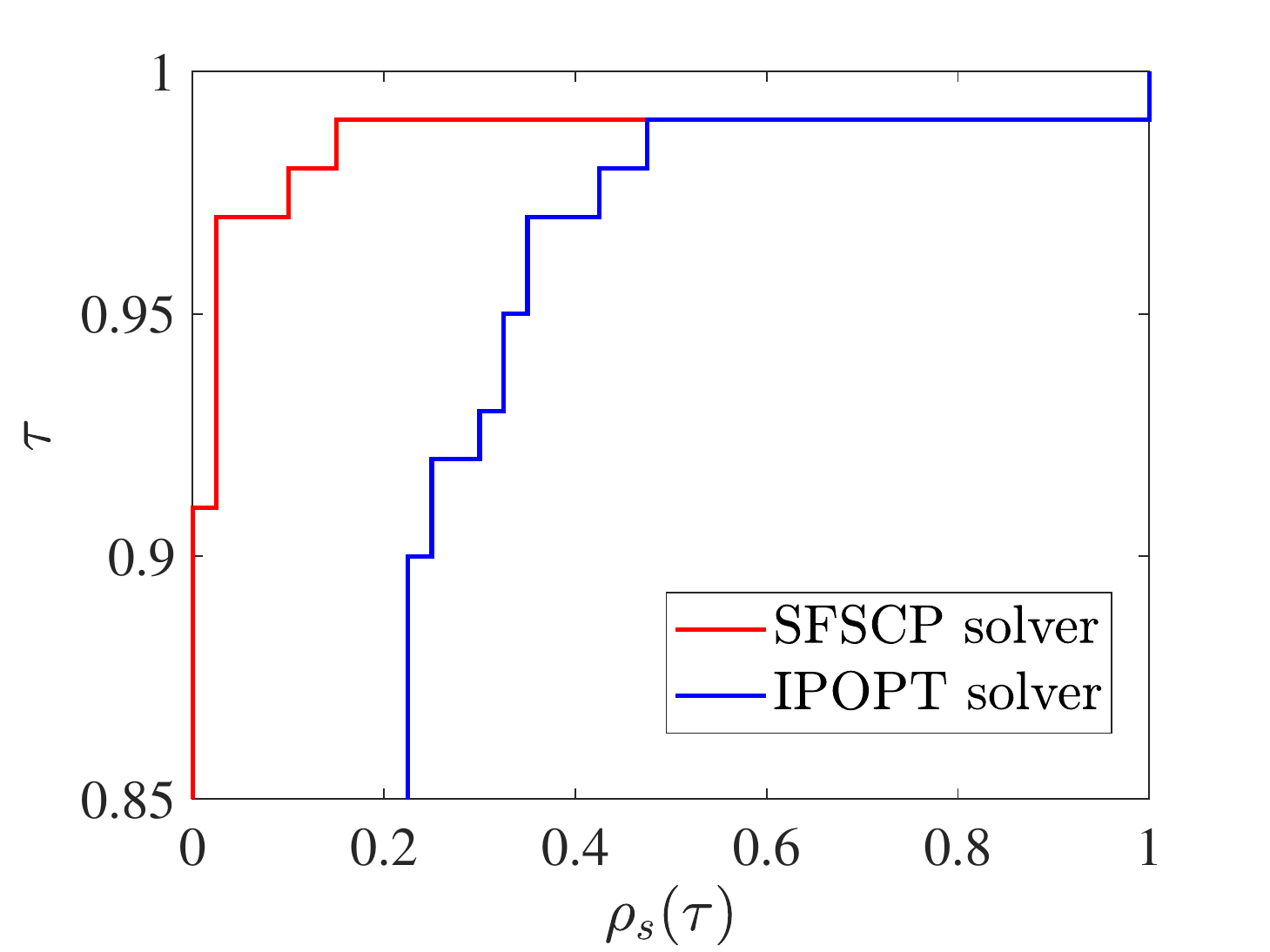}
    \caption{SFSCP and IPOPT solver performance profiles}
    \label{fig:solver_performance_cdf}
\end{figure}

We use the SFSCP solver in this work for multiple reasons. First, it guarantees strict hydraulic feasibility for each optimization step and shows fast convergence properties in previously studied WDN control problems, as demonstrated in \citet{WRIGHT2015}. Second, since state-of-the-art nonlinear optimization solvers (e.g. IPOPT) require second-order derivatives, it avoids the complications associated with the nonsmooth HW head loss model, which has an unbounded hessian at the origin. Potential complications are exacerbated by the highly nonlinear SCC objective function used in this study. In order to justify application of the SFSCP solver, though, we conduct a performance comparison with solutions computed by IPOPT. Following \citet{DOLAN2002}, we derive performance profiles to objectively compare the two solvers. Let $\mathcal{X} := \{1,\dots,n_x\}$ and $\mathcal{S} := \big\{s_\text{SFSCP}, s_\text{IPOPT}\big\}$ be the set of test experiments and solvers, respectively, where $n_x$ is the number of experiments tested. Since we are interested in comparing the SCC objective function values, we set the performance metric $f_{x,s}$ equal to $f_{\widetilde{\text{SCC}}}$ for each experiment $x$ produced by solver $s$. We then compare solver performance for all experiments $x \in \mathcal{X}$ by normalizing $f_{x,s}$ to the best computed solution for a particular experiment by any solver $s \in \mathcal{S}$. This comparison is expressed as follows:
\begin{equation}
    \label{eq:solver_performance_compute}
    r_{x,s} = \frac{f_{x,s}}{\min{\big\{f_{x,s} \; | \; s \in \mathcal{S}\big\}}},
\end{equation}
where $r_{x,s}$ is the performance ratio of solver $s$ for experiment $x$. Note that, if the solver produces an error or fails to find a feasible solution to the SCC control problem, we assign $f_{x,s} = +\infty$. An overall performance assessment of solver $s \in \mathcal{S}$ is then obtained by deriving the cumulative distribution function
\begin{equation}
    \label{eq:solver_performance_cdf}
    \rho_s(\tau) = \frac{1}{n_x}\text{size}\big\{x \in \mathcal{X} \; | \; r_{x,s} \leq \tau\big\},
\end{equation}
where $\rho_s(\tau)$ represents the percentage of experiments within a factor $\tau$ of the best solution for solver $s \in \mathcal{S}$. \Cref{fig:solver_performance_cdf} presents the SFSCP and IPOPT solver performance profiles for $n_x = 36$ design-for-control experiments. In addition to existing unidirectional pressure reducing valves (PRVs), these experiments involved a combination of $n_v = 1,\dots,3$ bidirectional dynamic boundary valves (DBVs) and $n_f = 0,\dots,3$ automatic flushing valves (AFVs), as described in Section 4.1 of the manuscript.

From \Cref{fig:solver_performance_cdf}, we first observe that the SFSCP solver finds feasible solutions for all experiments, while IPOPT returns infeasibility errors for just over 20\% of the tested cases. This is illustrated by the distance along the x-axis from the origin to the vertical segment of the IPOPT solver (blue) plot. When looking at individual experiments, we identified a few features which may be contributing to IPOPT's infeasibility errors. First, infeasibility occurs almost entirely for the larger case study network (e.g. \texttt{BWFLnet}), of which also has a greater range in operational demands and control features. Second, since good performance was recorded for the feasibility restoration problem (see Section 3.3.2 of the manuscript), we believe the highly nonlinear SCC objective function, in combination with the unbounded hessian of the HW model, to be a key factor in the encountered infeasibility errors. In any case, \Cref{fig:solver_performance_cdf} shows the SFSCP solver to be within 10\% of the best solution for all tested experiments, which is comparable to the IPOPT solver for experiments when no feasibility issues were present. This suggests that the SFSCP solver is appropriate for the current SCC design-for-problem study.

\section*{Acknowledgments}
This work was supported by EPSRC (EP/P004229/1, Dynamically Adaptive and Resilient Water Supply Networks for a Sustainable Future), Bristol Water Plc, Analytical Technology (ATi), and Imperial College London's Department of Civil and Environmental Engineering Skempton Scholarship.

\bibliography{mybibfile.bib}

\end{document}